\documentclass[a4paper,11pt,oneside,reqno]{amsart}
\usepackage{silence}
\allowdisplaybreaks
\usepackage[utf8]{inputenc}
\usepackage[pdftex]{graphicx}
\usepackage{amssymb}
\usepackage{esint}
\usepackage{amsmath}
\usepackage{amsfonts}
\usepackage{amsthm}
\usepackage{bm}
\usepackage{mathrsfs}
\usepackage{mathtools}
\usepackage{comment}
\mathtoolsset{showonlyrefs=true}
\usepackage{appendix}
\usepackage{enumerate}
\usepackage[a4paper, total={6.6in, 9.7in}]{geometry}
\usepackage[stretch=30,shrink=30]{microtype}
\IfFileExists{stmaryrd.sty}{\usepackage{stmaryrd}}{%
    \newcommand{\llbracket}{\mathopen{\text{\normalfont\fontfamily{lmr}\selectfont\textlbrackdbl}}}%
    \newcommand{\rrbracket}{\mathclose{\text{\normalfont\fontfamily{lmr}\selectfont\textrbrackdbl}}}%
}
\usepackage{xcolor}
\definecolor{antiquefuchsia}{rgb}{0.57, 0.36, 0.51}
\definecolor{azure}{rgb}{0.0, 0.5, 1.0}

\usepackage[backref=page, colorlinks = true, linkcolor = black, citecolor = antiquefuchsia]{hyperref}

\renewcommand*{\backref}[1]{}
\renewcommand*{\backrefalt}[4]{%
    \ifcase #1 (Not cited.)%
    \or        (Cited on page~#2.)%
    \else      (Cited on pages~#2.)%
    \fi}

\makeatletter
\def\th@plain{%
	\thm@notefont{}
	\itshape 
}
\def\th@definition{%
	\thm@notefont{}
	\normalfont 
}
\makeatother

\numberwithin{equation}{section}

\newtheorem{theorem}{Theorem}[section]
\newtheorem{lemma}[theorem]{Lemma}
\newtheorem{proposition}[theorem]{Proposition}
\newtheorem{corollary}[theorem]{Corollary}

\theoremstyle{definition}

\newtheorem{remark}[theorem]{Remark}

\renewcommand{\b}{\mathbb{B}}

\newcommand{\Hh}{\mathscr{H}}

\DeclareMathOperator{\spt}{spt}

\DeclareMathOperator{\vol}{vol}

\DeclareMathOperator{\tr}{tr}

\DeclareMathOperator{\Sing}{Sing}

\newcommand{\Reg}{\mathrm{Reg}}

\newcommand{\U}{\operatorname{U}}

\newcommand{\cp}{\mathbb{CP}}

\newcommand{\spn}{\mathrm{span}}

\hypersetup{bookmarksdepth=2,
    pdftitle={Partial Regularity of Weakly Pseudoholomorphic Maps Into Projective Algebraic Varieties},
    pdfauthor={Riccardo Caniato}}

\title[Sharp Partial Regularity of Weakly Pseudoholomorphic Maps]{Sharp Partial Regularity of Weakly Pseudoholomorphic Maps Into Projective Algebraic Varieties}

\author[Riccardo Caniato]{Riccardo Caniato}
\address{Warwick Mathematics Institute, Zeeman Building, University of Warwick, Coventry, CV4 7AL, United Kingdom}
\email{riccardo.caniato@warwick.ac.uk}
\date{\today} 

\begin{document}
\begin{abstract}
    For every $m\geq2$, we prove that a weakly pseudoholomorphic $W^{1,2}_{\mathrm{loc}}$ map from a real-analytic almost complex $2m$-manifold into a projective algebraic variety has a countably $(2m-4)$-rectifiable singular set. When $m=2$, the singular set is shown to be locally finite. We first establish both conclusions on smooth almost complex domains under the additional assumption that the pullback of the Fubini--Study form is distributionally closed. We then show that this conservation law follows from the Cauchy--Riemann equation alone on real-analytic domains.
\end{abstract}
\maketitle
\tableofcontents
\section{Introduction}
\label{Section: Introduction}
For $m\ge 2$, let $(M^{2m},J)$ be a smooth almost complex manifold equipped with a $J$-compatible Riemannian metric $g$, and consider a closed almost K\"ahler target manifold $(N,j_N,\omega_N)$.\footnote{\,By \textit{almost K\"ahler} manifold, we mean that the symplectic form $\omega_N$ is compatible with the almost complex structure $j_N$, and therefore induces the Riemannian metric $h:=\omega_N(\,\cdot\,,j_N\,\cdot\,)$.} A map $u\in W^{1,2}_{\mathrm{loc}}(M,N)$ is called \emph{weakly pseudoholomorphic} if
\begin{align}
\label{Equation: Intro pseudoholomorphic equation}
    du\circ J=j_N\circ du \qquad\mbox{ a.e. on }M.
\end{align}
The notion of pseudoholomorphic map connects the geometry of almost complex manifolds with the variational theory for the Dirichlet energy and calibrated geometry. Weakly pseudoholomorphic maps of finite Dirichlet energy on Riemann surfaces, known as \emph{pseudoholomorphic curves}, are always smooth. This remarkable regularity property underlies, for instance, Gromov's foundational work on symplectic rigidity; see \cite{GromovPseudoholomorphicCurves}. Since weakly pseudoholomorphic maps on domains of dimension at least four may exhibit nonremovable singularities despite having finite Dirichlet energy, understanding the size and geometric structure of their singular sets is a natural problem in geometric analysis.

We call $x\in M$ a \emph{regular point} of a map $u\in W^{1,2}_{\mathrm{loc}}(M,N)$ if $u$ admits a smooth representative in a neighbourhood of $x$. The open set $\Reg\,u\subset M$ formed by all such points is called the \emph{regular set} of the map, while $\Sing u:=M\smallsetminus\Reg\,u$ is referred to as its \emph{singular set}. Any point $x\in\Sing u$ is called a \emph{singular point}, or a \emph{singularity}, of $u$. Throughout, Hausdorff dimensions and Hausdorff measures are computed using $g$. Since smooth Riemannian metrics are locally equivalent, the conclusions below are independent of this choice. 

For a Sobolev map $u\in W^{1,2}_{\mathrm{loc}}(M,N)$, the size of its singular set is closely related to whether the pullback of the target symplectic form is distributionally closed. Indeed, for a smooth map $u$ we trivially have $d(u^*\omega_N)=u^*(d\omega_N)=0$, whereas if $u$ is just $W^{1,2}_{\mathrm{loc}}$ the pullback $u^*\omega_N$ is still a well-defined locally integrable 2-form, but closedness may fail on $\Sing u$. For instance, consider the finite-energy map $v(x)=x/|x|$ from the unit ball in $\mathbb R^3$ to the unit sphere $\mathbb S^2$. If $\omega_{\mathbb S^2}$ denotes the standard area form on the 2-sphere, then $v^*\omega_{\mathbb S^2}$ is closed off the origin and has integral $4\pi$ over every 2-sphere centred there. Consequently, $\ast d(v^*\omega_{\mathbb S^2})=4\pi\delta_0$ in the sense of distributions (see e.g. \cite{BrezisCoronLieb}). We are then led to consider the geometric conservation law
\begin{align}
\label{Equation: Intro closedness assumption}
    d(u^*\omega_N)=0 \qquad\text{in the sense of distributions on }M.
\end{align}
To reinforce our point that failure of \eqref{Equation: Intro closedness assumption} leads to obstructions to regularity, note that
\begin{align*}
    d(u^*\omega_N)\neq 0 \qquad\Rightarrow\qquad \Hh^{2m-3}(\Sing u)>0.
\end{align*}
Indeed, suppose that $\Hh^{2m-3}(\Sing u)=0$. Fix three coordinate directions in a chart. Since the projection of $\Sing u$ onto the complementary $(2m-3)$-plane has measure zero, almost every three-dimensional slice parallel to the chosen directions misses $\Sing u$. The restriction of $u^*\omega_N$ to each such slice is smooth and closed. Integrating by parts on the slices and applying Fubini's theorem gives the vanishing of the corresponding component of $d(u^*\omega_N)$ in distributions. Varying the three directions gives \eqref{Equation: Intro closedness assumption}.

To study the regularity of a map $u$ satisfying \eqref{Equation: Intro pseudoholomorphic equation} and \eqref{Equation: Intro closedness assumption}, we first consider the case in which the domain is also almost K\"ahler. Since its fundamental form $\Omega:=g(J\,\cdot\,,\,\cdot\,)$ is closed, $u$ can be shown to be \emph{stationary harmonic}, i.e. critical for the Dirichlet energy under compactly supported outer and inner variations \cite{RiviereTian,CaniatoRiviere}. By the partial regularity theorem of Evans \cite{EvansStationaryHarmonicMaps} for sphere targets and its extension by Bethuel \cite{BethuelStationaryHarmonicMaps} to general smooth targets,\footnote{\,By exploiting the antisymmetric structure of the harmonic map equation, Rivi\`ere--Struwe \cite{RiviereStruwe} obtained the corresponding partial H\"older regularity result for $C^2$ targets.} it follows that 
\begin{align}\label{Equation: Intro preliminary bound}
    \Hh^{2m-2}(\Sing u)=0.
\end{align}
On a general almost complex domain, $\Omega$ need not be closed but the same regularity estimate can still be obtained. Indeed, the current $T_u$ defined by $u^*\omega_N$ is semicalibrated by $\Omega^{m-1}/(m-1)!$, and its mass measure is the Dirichlet energy measure $\frac12|du|^2\,d\vol_g$. Since \eqref{Equation: Intro closedness assumption} makes $T_u$ a cycle, the almost-monotonicity formula for its mass applies to the rescaled Dirichlet energy of $u$ \cite{PumbergerRiviere,CaniatoRiviere}. Together with Wang's $\varepsilon$-regularity theorem \cite{WangJHolomorphic}, this implies \eqref{Equation: Intro preliminary bound}.

When both $J$ and $j_N$ are integrable, complex analytic geometry allows us to obtain a stronger bound. Indeed, by Sobolev slicing and regularity in complex dimension one, the restriction of a weakly holomorphic map $u\in W^{1,2}_{\mathrm{loc}}(M,N)$ to almost every coordinate disc in a holomorphic chart has a holomorphic representative. Shiffman's theorem \cite[Theorem~1]{ShiffmanSeparatelyMeromorphic} therefore implies that $u$ agrees almost everywhere with a meromorphic map (see also \cite{UhlenbeckYauSubbundles,Popovici} for the projective
case). Thus, $\Sing u$ coincides with the indeterminacy locus of this map. Since this locus is an analytic set of complex codimension at least two, it has locally finite $\Hh^{2m-4}$-measure (see \cite[Theorem~3.1.1]{King} and \cite[Proposition~1.6]{HarveyShiffman}).

If $J$ is nonintegrable, new difficulties arise due to the lack of complex coordinates. Nevertheless, for $m=2$, Rivi\`ere--Tian \cite{RiviereTian} proved that a weakly pseudoholomorphic map into a projective algebraic variety has isolated singularities if it is locally a strong $W^{1,2}$ limit of smooth maps.\footnote{\,For strong approximation in the ambient $\cp^n$, this hypothesis is equivalent to \eqref{Equation: Intro closedness assumption} for the Fubini--Study form \cite{BethuelCoronDemengelHelein}.} Although their proof uses positivity of intersections in a way specific to real dimension four, it is natural to ask if a similar form of improved regularity can be obtained on general almost complex domains. In this paper, we give a positive answer to this question. 

\subsection{Main results}

We define a \emph{complex projective algebraic variety} $X\subset\cp^n$ as the common zero set of finitely many homogeneous complex polynomials in $n+1$ variables. As this allows $X$ to be singular or reducible, we regard maps into $X$ as maps into $\cp^n$ taking values in $X$ almost everywhere. In particular, we interpret smoothness and the pseudoholomorphic equation in the ambient projective space. We write $j$ for its standard complex structure and $\omega$ for its Fubini--Study form, normalised by
\begin{align*}
    \int_{\cp^1}\iota^*\omega=1,
\end{align*}
where $\iota:\cp^1\hookrightarrow\cp^n$ is any projective linear embedding. As all projective lines are homologous, this normalisation is independent of the chosen line.

We say that $(M,J)$ is \emph{real analytic} if $M$ admits a real-analytic atlas in which the coefficients of $J$ are analytic. In particular, this holds for the nonintegrable structures studied by Holt--Zhang on the Kodaira--Thurston manifold \cite[Section~2]{HoltZhang} and by Di Scala--Vezzoni on the Iwasawa manifold \cite[Theorem~4.3]{DiScalaVezzoni}. Indeed, in both cases $J$ descends from an invariant tensor on a real-analytic Lie group. Another example is the Eells--Salamon structure on the twistor space of an oriented Riemannian four-manifold \cite{EellsSalamon}. For a real-analytic base metric, this structure is real analytic because it is defined in terms of the horizontal distribution of the Levi-Civita connection and the complex structure on the twistor fibres.

\begin{theorem}
    \label{Theorem: Real analytic partial regularity}
    Let $(M^{2m},J)$ be a real-analytic almost complex manifold without boundary, with $m\geq2$, and let $X\subset\cp^n$ be a complex projective algebraic variety, with $n\geq1$. If $u\in W^{1,2}_{\mathrm{loc}}(M,X)$ is weakly pseudoholomorphic, then $\Sing u$ is closed and countably $(2m-4)$-rectifiable. In particular, 
    \begin{align}\label{Equation: Intro Haussdorf dimension bound}
        \dim_{\Hh}\Sing u\leq2m-4
    \end{align}
    and $\Hh^{2m-4}\,\!\llcorner\Sing u$ is $\sigma$-finite. If $m=2$, then $\Sing u$ is locally finite.
\end{theorem}
Thus, for maps from real-analytic domains into projective algebraic varieties, satisfying the Cauchy--Riemann equation alone forces a very strong form of regularity. This should be compared with the fact that weakly harmonic maps on domains of dimension larger than $2$ can be discontinuous everywhere (see \cite{RiviereDiscontinuousHarmonicMaps}), and even imposing stronger variational assumptions leads to worse bounds than \eqref{Equation: Intro Haussdorf dimension bound}. For instance, the radial map $v$ above, extended by constancy in $2m-3$ further variables, solves the harmonic map equation and locally minimises the Dirichlet energy  \cite{BrezisCoronLieb}. Yet, its singular set has real codimension three, attaining the bound of Schoen--Uhlenbeck \cite{SchoenUhlenbeck}.
\begin{remark}
    The dimension bound is sharp even for maps into $\cp^1$. Indeed, the map
    \begin{align*}
        u(z_1,\ldots,z_m)=[z_1:z_2]
    \end{align*}
    belongs to $W^{1,2}_{\mathrm{loc}}(\mathbb C^m,\cp^1)$, is weakly holomorphic, and has singular set precisely the $(2m-4)$-dimensional plane $\{z_1=z_2=0\}$ (see Proposition~\ref{Proposition: Sharpness of codimension four}). One can also realise this example on a nonintegrable domain. On the Iwasawa manifold above, two members of the invariant $(1,0)$-coframe are closed \cite[Eq.~(4.2)]{DiScalaVezzoni}. Their local primitives $z_1,z_2$, chosen to vanish at a given point, define a pseudoholomorphic submersion $f$ to $\mathbb C^2$. Composing $f$ with $(w_1,w_2)\mapsto[w_1:w_2]$ gives, after assigning a value on $f^{-1}(0)$, a weakly pseudoholomorphic map whose singular set is precisely the fibre over the origin, of real codimension four.
\end{remark}
For the maps in Theorem~\ref{Theorem: Real analytic partial regularity} one can also describe the blow-up at each singular point. Indeed, Theorem~\ref{Theorem: Automatic closedness} gives $d(u^*\omega)=0$, and since $[\omega]$ generates $H^2_{\mathrm{dR}}(\cp^n;\mathbb R)$, Lemma~\ref{Lemma: Sobolev pullbacks} gives $d(u^*\beta)=0$ for every smooth closed $2$-form $\beta$ on $\cp^n$. The criterion of Bethuel--Coron--Demengel--H\'elein \cite{BethuelCoronDemengelHelein} therefore gives local strong $W^{1,2}$ approximation by smooth $\cp^n$-valued maps. It follows from the uniqueness theorems of Caniato--Rivi\`ere \cite{CaniatoRiviere} that at every singular point both the tangent map and the tangent cone to the current defined by $u^*\omega$ are unique. By the energy identity \eqref{Equation: Energy identity}, the rescaled Dirichlet energy measures then converge weakly to the mass measure of this cone.

Since the tangent map is unique in our setting, a singular point belongs to the $k$th tangent-map stratum precisely when its tangent map has at most $k$ independent directions of translation invariance. Naber--Valtorta \cite{NaberValtortaHarmonicMaps} proved that these strata are rectifiable for stationary harmonic maps and that, for energy minimisers, the full singular set is countably $(2m-3)$-rectifiable with locally finite $\Hh^{2m-3}$-measure. For stationary pseudoholomorphic maps, Guo--Liu--Wang \cite[Theorem~1.3 and Lemma~3.4]{GuoLiuWang} obtained sharp tubular estimates for the quantitative strata and proved rectifiability of the corresponding tangent-map strata. To see why the latter can miss singular points, consider $u(z_1,z_2)=[z_1:z_2^2]$. Although $u$ is singular at the origin, its rescalings $u(rz_1,rz_2)=[z_1:r z_2^2]$ converge in $L^2_{\mathrm{loc}}$ to the constant map $[1:0]$ as $r\rightarrow0$. As this tangent map is invariant under every translation, the origin lies outside every lower-dimensional tangent-map stratum. Rectifiability of these strata alone therefore cannot be used to estimate the size of the set of singular points with constant tangent map, whereas Theorem~\ref{Theorem: Real analytic partial regularity} proves that the full singular set is countably $(2m-4)$-rectifiable.

\smallskip
The proof of Theorem~\ref{Theorem: Real analytic partial regularity} rests on two independent results. First, in Sections~\ref{Section: Preliminaries}--\ref{Section: Partial regularity and sharpness}, we establish partial regularity on smooth almost complex domains treating \eqref{Equation: Intro closedness assumption} as an additional assumption.

\begin{theorem}
    \label{Theorem: Intro main regularity}
    Let $(M^{2m},J)$ be a smooth almost complex manifold without boundary, with $m\geq2$, and let $X\subset\cp^n$ be a complex projective algebraic variety, with $n\geq1$. If $u\in W^{1,2}_{\mathrm{loc}}(M,X)$ is weakly pseudoholomorphic and satisfies \eqref{Equation: Intro closedness assumption}, then $\Sing u$ is closed and countably $(2m-4)$-rectifiable. In particular, 
    \begin{align*}
        \dim_{\Hh}\Sing u\leq2m-4
    \end{align*}
    and $\Hh^{2m-4}\,\!\llcorner\Sing u$ is $\sigma$-finite. If $m=2$, then $\Sing u$ is locally finite.
\end{theorem}
We then show in Section~\ref{Section: Automatic closedness} that the conservation law \eqref{Equation: Intro closedness assumption} is automatically satisfied if $J$ is real-analytic.
\begin{theorem}
    \label{Theorem: Automatic closedness}
    Let $(M^{2m},J)$ be a real-analytic almost complex manifold without
    boundary, with $m\geq1$, and let $n\geq1$. Every weakly
    pseudoholomorphic map $u\in W^{1,2}_{\mathrm{loc}}(M,\cp^n)$
    satisfies
    \begin{align*}
        d(u^*\omega)=0
        \qquad\text{in the sense of distributions on }M.
    \end{align*}
\end{theorem}
%

\subsection{Strategy of the proofs}

The first step in the proof of Theorem~\ref{Theorem: Intro main regularity}
is to characterise the singular points of $u$ by their energy density.
To this end, set
\begin{align*}
    \Omega(X,Y):=g(JX,Y),\qquad
    \Phi:=\frac{\Omega^{m-1}}{(m-1)!},\qquad
    \mu_u:=\frac12|du|^2\,d\vol_g .
\end{align*}
By the energy identity we have $\mu_u=u^*\omega\wedge\Phi$. Since the
current defined by $u^*\omega$ is a cycle, almost monotonicity gives
the existence of the density
\begin{align*}
    \theta_u(x):=\lim_{r\rightarrow0}r^{2-2m}\mu_u(B_r(x)).
\end{align*}
Note that $\theta_u$ vanishes at every regular point. Conversely, small-energy
regularity gives a positive lower bound for $\theta_u$ at singular
points, uniform on compact subsets of $M$. Thus, $\Sing u=\{x\in M:\theta_u(x)>0\}$, and it suffices to estimate the
size of the positive-density set of $\mu_u$.

To study this set, we disintegrate $\mu_u$ as an average over the dual
projective space $A=(\cp^n)^*$, equipped with its unitary-invariant
probability measure $\lambda$. For almost every hyperplane $H\in A$,
the inverse image
\begin{align*}
    F_H:=\{x\in\Reg\,u\mbox{ : }u(x)\in H\}
\end{align*}
is a smooth pseudoholomorphic hypersurface of $\Reg\,u$. Hence, integration over $F_H$ with its complex orientation gives
a cycle on $\Reg\,u$. To extend it to a cycle on $M$, we must
control its mass near $\Sing u$ and show that its extension acquires
no boundary there.

For this purpose we use a weighted incidence formula. For each smooth
real function $f$ on $A$, integration along the fibres of the
projective incidence manifold gives a smooth closed $2$-form $\beta_f$
on $\cp^n$ satisfying
\begin{align*}
    \beta_1=\omega,\qquad
    \beta_f=\left(\int_A f\,d\lambda\right)\omega+d\eta_f
\end{align*}
for some smooth $1$-form $\eta_f$. On $\Reg\,u$, where $u$ is smooth,
the corresponding incidence formula reads
\begin{align}\label{Equation: Intro weighted incidence}
    \int_A f(H)\left(\int_{F_H}\zeta\right)d\lambda(H)
    =\int_{\Reg\,u}u^*\beta_f\wedge\zeta
\end{align}
for every smooth compactly supported $(2m-2)$-form $\zeta$ there.
Taking $f=1$ and $\zeta=\chi\Phi$, we see that averaging the mass measures of the $F_H$ over $H$ gives $\mu_u$ on $\Reg\,u$. As $\mu_u$ is locally finite across $\Sing u$, almost every $F_H$ has locally finite mass there as well. Thus, integration over $F_H$ defines a current $T_H$ on $M$, and extension by zero gives the weighted identity on all of $M$.

To see that these extensions are cycles, we use the freedom to vary
$f$. Indeed, the formula for $\beta_f$, together with
\eqref{Equation: Intro closedness assumption} and
Lemma~\ref{Lemma: Sobolev pullbacks}, gives $d(u^*\beta_f)=0$. Thus, taking $\zeta=d\xi$ for a smooth compactly supported $(2m-3)$-form $\xi$ on $M$, we see that the weighted integral on the left vanishes for every $f$. It follows that $\partial T_H=0$ for almost every $H$, and the unweighted mass identity becomes
\begin{align*}
    \mu_u=\int_A\|T_H\|\,d\lambda(H).
\end{align*}

All these cycles are semicalibrated by the same form $\Phi$, since
$\Phi$ has comass one and restricts to the volume form on each
oriented complex hyperplane. If $\Phi$ is closed, it is a calibration
in the sense of Harvey--Lawson \cite{HarveyLawson} and the cycles $T_H$ minimise mass in homology \cite{MorganCalibrations1988,MorganCalibrations1990}.
Here we use the common semicalibration to rule out mass cancellation in the sums $T_H+T_{H'}$. For almost every pair $(H,H')$, the mass measures of $T_H$ and $T_{H'}$ are also mutually singular. These two facts imply that
\begin{align*}
    \spt T_H\cap\spt T_{H'}\subset\Sing(T_H+T_{H'})
    \qquad\text{for a.e. }(H,H').
\end{align*}
Indeed, if the sum were regular at a point of both supports, the
constancy theorem would force both cycles to have positive constant
multiplicity on the same local sheet, contradicting mutual singularity.

It remains to bound the singular set of the sum. For area-minimising
currents, the codimension-two bound is due to Almgren
\cite{AlmgrenBigRegularity}, with a later proof by De Lellis--Spadaro
\cite{DeLellisSpadaroRegularityI,DeLellisSpadaroRegularityII,DeLellisSpadaroRegularityIII}.
Spolaor's extension to semicalibrated currents \cite{Spolaor} applies
to $T_H+T_{H'}$. Since the sum has dimension $2m-2$, the pairwise
intersection has Hausdorff dimension at most $2m-4$. When $m=2$ we
obtain more, since the sum is a positive $(1,1)$ integral cycle and
every smooth almost complex structure in real dimension four admits
locally compatible symplectic forms \cite[Lemma~A.1]{RiviereTian}
(see also \cite[Theorem~1]{LejmiLocalSymplectic}). Thus the theorem
of Rivi\`ere--Tian \cite{RiviereTianCurrents} gives isolated
singularities for the sum and hence local finiteness of the
intersection. Note that the theorem of De Lellis--Spadaro--Spolaor
for semicalibrated integral $2$-currents
\cite{DeLellisSpadaroSpolaorTangentCones,DeLellisSpadaroSpolaorI,DeLellisSpadaroSpolaorII,DeLellisSpadaroSpolaor}
gives the same conclusion and includes the case of special
Legendrian cycles in $\mathbb S^5$ treated by Bellettini--Rivi\`ere
\cite{BellettiniRiviereLegendrian}.

To pass from these pairwise estimates to one for $\Sing u$, consider
the probability that the support of a cycle contains a given point, i.e.
\begin{align*}
    \mathfrak p(x):=\lambda\big(\{H\in A\mbox{ : }x\in\spt T_H\}\big).
\end{align*}
Uniform almost monotonicity for the cycles allows us to pass to
the density limit under the integral defining $\mu_u$. Thus
$\theta_u(x)$ is the average of the densities of the cycles at $x$.
As each density vanishes off the corresponding support, this gives
$\mathfrak p(x)>0$ wherever $\theta_u(x)>0$. The same estimate gives
a uniform positive lower bound for $\mathfrak p$ on compact portions
of $\Sing u$. For two independently chosen hyperplanes
the probability that both supports contain $x$ is $\mathfrak p(x)^2$.
Thus, for any finite Borel measure $\nu$ on $M$, Tonelli's theorem gives
\begin{align*}
    \int_{A\times A}\nu(\spt T_H\cap\spt T_{H'})\,d(\lambda\otimes\lambda)(H,H')
    =\int_M\mathfrak p(x)^2\,d\nu(x).
\end{align*}
If $\dim_{\Hh}\Sing u>2m-4$, Frostman's lemma would give a nonzero finite
measure $\nu$ on $\Sing u$ vanishing on all sets of dimension at
most $2m-4$. The left-hand side would then vanish, whereas the
right-hand side is positive. In dimension four, suppose a compact set
contains a sequence of distinct singular points $x_j$. Almost every
pairwise intersection contains only finitely many of them by
Rivi\`ere--Tian \cite{RiviereTianCurrents}. Dominated convergence
therefore gives $\mathfrak p(x_j)^2\to0$, contrary to the uniform
positive lower bound on that compact set.

For $m>2$ we use the incidence identity once more to prove
rectifiability of $\Sing u$. In the area-minimising case,
rectifiability in codimension two was proved by
De Lellis--Skorobogatova \cite{DeLellisSkorobogatovaI,DeLellisSkorobogatovaII}
and De Lellis--Minter--Skorobogatova \cite{DeLellisMinterSkorobogatovaIII}.
Krummel--Wickramasekera \cite{KrummelWickramasekeraI,KrummelWickramasekeraII}
developed an independent approach to these singularities. The
extension to semicalibrated currents due to
Minter--Parise--Skorobogatova--Spolaor \cite{MinterPariseSkorobogatovaSpolaor}
gives countable $(2m-4)$-rectifiability of almost every pairwise
intersection. If $F\subset\Sing u$ were a purely $(2m-4)$-unrectifiable Borel set
of positive finite $\Hh^{2m-4}$-measure, taking
$\nu=\Hh^{2m-4}\,\!\llcorner F$ in the incidence identity would again
give zero on the left and a positive integral on the right. By
De Lellis--Fleschler \cite[Theorem~1.1]{DeLellisFleschler}, the
absence of such a set implies that $\Hh^{2m-4}\,\!\llcorner\Sing u$ is
$\sigma$-finite. The rectifiable--purely unrectifiable decomposition on each
finite-measure piece then gives the asserted rectifiability. The
underlying averaging argument applies to general families of
semicalibrated cycles, as formulated in Theorem~\ref{Theorem: Averaged cycles}.

For the proof of the conservation law, we first take the target to
be $\cp^1$ and proceed by induction on $m$. Where $J$ is integrable,
almost every level current of $u$ has the same boundary as the current
defined by $u^*\omega$ \cite{AlbertiBaldoOrlandi}. In particular,
the difference of two such currents is closed. The tangent planes of
the levels are complex hyperplanes, so Shiffman's codimension-one
theorem \cite{ShiffmanCodimensionOne} identifies the difference with
a signed holomorphic chain. Since the level currents have the
positive complex orientation and are carried by disjoint levels,
they are precisely the positive and negative parts of this chain.
Both parts are closed. Their common boundary therefore vanishes,
giving $d(u^*\omega)=0$.

Where $J$ is nonintegrable, $du$ annihilates the image of the
Nijenhuis tensor of $J$, since the target complex structure is
integrable. Hence, the module of germs of real-analytic vector fields annihilated by $du$ is nonzero. Since this module is finitely generated, we may choose a finite family of local generators. Wherever one of these fields
is nonzero, $u$ factors locally through a weakly pseudoholomorphic map
on an almost complex domain of real dimension $2m-2$. By induction,
$u^*\omega$ is therefore closed off the common zero set of the fields.
To describe this zero set, we use real analyticity to continue arcs
in it to families of pseudoholomorphic discs contained in the set.
It follows that the strata of codimension two are almost complex
and those of codimension three lie locally in almost complex
hypersurfaces.
Elkhadhra's removability theorem \cite{ElkhadhraPluripolar,ElkhadhraPoincareLelong}
then applies across both kinds of strata, since $u^*\omega$ defines
a positive current. Thus $d(u^*\omega)$ is supported on the remaining strata of
codimension at least four. Since $u^*\omega$ is a locally integrable $2$-form, the support
estimate of Lemma~\ref{Lemma: High codimension support} now gives
$d(u^*\omega)=0$.

To pass from $\cp^1$ to $\cp^n$, we compose $u$ with the meromorphic projections $\cp^n\rightarrow\cp^{n-1}$ considered in \cite[Appendix~A]{CaniatoRiviere}. Almost every composition is weakly pseudoholomorphic, so its pullback of the Fubini--Study form is closed by induction on $n$. Averaging these pullbacks recovers $u^*\omega$ in $L^1_{\mathrm{loc}}$ and gives closedness for the
original map.
\subsection{Further questions}
For $m>2$, it remains open whether $\Sing u$ has locally finite $\Hh^{2m-4}$-measure, as conjectured in \cite{RiviereTian,RiviereTianCurrents,CaniatoRiviere}, even for real-analytic $J$. Proposition~\ref{Proposition: Conditional local finiteness} shows that it would suffice to prove local finiteness of the same measure on $\spt T_H\cap\spt T_{H'}$ for almost every pair. The bounds may depend on the pair and need not be integrable in $(H,H')$.

A second question is whether the conservation law continues to hold when $J$ is only smooth. For projective targets, an affirmative answer would remove the real-analyticity assumption from Theorem~\ref{Theorem: Real analytic partial regularity}, since Theorem~\ref{Theorem: Intro main regularity} already applies on smooth domains. 

Concerning the extension of the main results beyond projective targets, our proof uses the fact that the inverse image of almost every projective hyperplane defines a cycle on $M$ and that averaging the mass measures of these cycles recovers the Dirichlet energy measure of $u$. It is therefore natural to ask which other target manifolds admit a family of hypersurfaces with the same property.

\subsection*{On the use of AI}
The author used ChatGPT-6 Astra solely to improve the clarity and style of the prose and to assist in identifying potential errors or gaps in the mathematical arguments. The research ideas and all the proofs were developed without the use of AI tools. All mathematical arguments were independently checked by the author, who takes full responsibility for the content of this paper.

\subsection*{Acknowledgements}
We would like to thank Tristan Rivi\`ere for pointing out this problem to the author, and for the many enlightening discussions we had over the years concerning the subject matter of this paper.

We are also thankful to Bal\'azs M\'ark B\'ek\'esi, Mattia Freguglia, and Alessandro Pigati for stimulating the author's thoughts on the content of Section 6.

Finally, we are grateful to Costante Bellettini, Davide Parise, Yannick Sire, and Luca Spolaor for several interesting conversations on questions related to this work.

%

\section{Preliminaries}
\label{Section: Preliminaries}

Throughout Sections~\ref{Section: Preliminaries}--\ref{Section: Partial regularity and sharpness}, $(M^{2m},J)$ is a smooth almost complex manifold equipped with a smooth $J$-compatible Riemannian metric $g$. We endow $M$ with the orientation induced by $J$ and define
\begin{align*}
    \Omega(X,Y):=g(JX,Y),
    \qquad
    \Phi:=\frac{\Omega^{m-1}}{(m-1)!}.
\end{align*}
We retain the notation $j$ and $\omega$ for the standard complex structure and the Fubini--Study form on $\cp^n$, normalised as in the introduction. We write $h:=\omega(\,\cdot\,,j\,\cdot\,)$ for the metric induced on $\cp^n$ by $j$ and $\omega$. Norms of differentials are computed using $g$ and $h$, and $B_r(x)$ denotes the geodesic ball in $(M,g)$.

We define the space $W^{1,2}_{\mathrm{loc}}(M,\cp^n)$ using a fixed smooth isometric embedding of $(\cp^n,h)$ into a Euclidean space. For a map $u$ in this class, the approximate differential defines the pullback $u^*\omega$ as an $L^1_{\mathrm{loc}}$ 2-form. Let $u\in W^{1,2}_{\mathrm{loc}}(M,\cp^n)$ be such that
\begin{align}
    \label{Equation: Assumptions on the map}
    \begin{split}
        du\circ J&=j\circ du
        \qquad\mbox{ a.e. on }M,\\
        d(u^*\omega)&=0
        \qquad\qquad\mbox{ in the sense of distributions on } M.
    \end{split}
\end{align}
We associate with $u$ its energy measure $\mu_u$ and its rescaled energy $\mathscr{E}_u$, given respectively by
\begin{align*}
    \mu_u:=\frac12|du|^2\,d\vol_g,
    \qquad
    \mathscr{E}_u(x,r):=r^{2-2m}\mu_u(B_r(x)).
\end{align*}
All geometric constants are chosen uniformly on compact subsets of $M$.

\subsection{The energy identity and Sobolev pullbacks}
\label{Subsection: Energy and Sobolev pullbacks}

The next lemma records the classical Wirtinger inequality and the associated energy identity. The comass statement is given in \cite[Section~II.3, Example~I]{HarveyLawson}. The energy identity for weakly pseudoholomorphic maps and the corresponding inequality for arbitrary Sobolev maps can be found in \cite[Corollary~2.1 and Lemma~3.1]{CaniatoRiviere}. We include the pointwise computation to specify our normalisations and the exact defect term. 
\begin{lemma}
    \label{Lemma: Energy identity}
    The form $\Phi$ has comass one and semicalibrates every
    $J$-holomorphic $(2m-2)$-plane with its complex orientation.
    Moreover, for every $v\in W^{1,2}_{\mathrm{loc}}(M,\cp^n)$,
    \begin{align}
        \frac12|d v|^2\,d\vol_g-v^*\omega\wedge\Phi
        =\frac14|d v\circ J-j\circ dv|^2\,d\vol_g
        \qquad\text{a.e. on }M.
        \label{Equation: Wirtinger energy defect}
    \end{align}
    In particular, every weakly pseudoholomorphic map $u$ satisfies
    \begin{align}
        u^*\omega\wedge\Phi
        =\frac12|d u|^2\,d\vol_g.
        \label{Equation: Energy identity}
    \end{align}
\end{lemma}
\begin{proof}
    Fix $x\in M$. Since $g$ is $J$-compatible, there is a
    positively oriented $g$-orthonormal basis
    $\{e_1,Je_1,\ldots,e_m,Je_m\}$ of $T_xM$.
    Denote its dual basis by
    $\{\varepsilon^1,\eta^1,\ldots,\varepsilon^m,\eta^m\}$.
    By the definition of $\Omega$, we have
    \begin{align*}
        \Omega_x
        &=\sum_{a=1}^m\varepsilon^a\wedge\eta^a,\\
        (d\vol_g)_x
        &=(\varepsilon^1\wedge\eta^1)\wedge\cdots
        \wedge(\varepsilon^m\wedge\eta^m),\\
        \Phi_x
        &=\sum_{a=1}^m
        \bigwedge_{\substack{1\leq b\leq m\\b\neq a}}
        (\varepsilon^b\wedge\eta^b).
    \end{align*}
    In the last expression the factors are ordered by increasing
    $b$. Expanding the wedge product shows that every 2-form
    $\zeta\in\Lambda^2T_x^*M$ satisfies
    \begin{align}
        \zeta\wedge\Phi_x
        =\left(\sum_{a=1}^m\zeta(e_a,Je_a)\right)
        (d\vol_g)_x.
        \label{Equation: Contraction with the fundamental form}
    \end{align}
    Indeed, in the product with the term of $\Phi_x$ omitting
    the $a$-th pair, only the coefficient of
    $\varepsilon^a\wedge\eta^a$ in $\zeta$ contributes.

    The same expressions show that $\Phi_x=*\Omega_x$, where
    the Hodge star is taken with respect to $g$ and the complex
    orientation. Let $\xi$ be an oriented unit simple
    $(2m-2)$-vector, and choose an orthonormal basis $p,q$ of
    its orthogonal complement such that $\xi\wedge p\wedge q$
    has the complex orientation. The definition of the Hodge
    star then gives
    \begin{align*}
        |\Phi_x(\xi)|=|\Omega_x(p,q)|
        =|g_x(Jp,q)|\leq1.
    \end{align*}
    If the plane represented by $\xi$ is $J$-invariant and
    carries its complex orientation, its orthogonal complement
    is also $J$-invariant and has the complementary complex
    orientation. We can therefore choose $q=Jp$, obtaining
    $\Phi_x(\xi)=g_x(Jp,Jp)=1$.
    This proves both the comass and semicalibration assertions.

    To prove the energy identity, fix a point $x$ where
    $v$ has an approximate differential
    \begin{align*}
        d v_x:T_xM\rightarrow T_{v(x)}\cp^n.
    \end{align*}
    At such a point, the pullback is the alternating bilinear
    form defined by
    \begin{align*}
        (v^*\omega)_x(X,Y)
        :=\omega_{v(x)}(d v_xX,d v_xY)
        \qquad(X,Y\in T_xM).
    \end{align*}
    Since $|v^*\omega|\leq C|d v|^2$, this formula defines a locally
    integrable 2-form. In the following calculation all
    domain quantities are evaluated at $x$, and all target
    quantities at $v(x)$.

    Applying \eqref{Equation: Contraction with the fundamental form}
    to $\zeta=(v^*\omega)_x$ gives
    \begin{align*}
        (v^*\omega\wedge\Phi)_x
        =\left(\sum_{a=1}^m
        \omega(d v_xe_a,d v_xJe_a)\right)(d\vol_g)_x.
    \end{align*}
    Since $j$ is an isometry for $h$ and
    $\omega(V,W)=h(jV,W)$, for each $a$ we have
    \begin{align*}
        |d v_xJe_a-jd v_xe_a|_h^2
        &=|d v_xJe_a|_h^2+|d v_xe_a|_h^2
        -2h(d v_xJe_a,jd v_xe_a)\\
        &=|d v_xJe_a|_h^2+|d v_xe_a|_h^2
        -2\omega(d v_xe_a,d v_xJe_a).
    \end{align*}
    Summing over $a$ and using the Hilbert--Schmidt norm of
    the differential yields
    \begin{align*}
        \frac12|d v_x|^2
        -\sum_{a=1}^m\omega(d v_xe_a,d v_xJe_a)
        =\frac12\sum_{a=1}^m|d v_xJe_a-jd v_xe_a|_h^2.
    \end{align*}
    Finally, the relations $J^2=-\mathrm{Id}$ and
    $j^2=-\mathrm{Id}$ give
    \begin{align*}
        (d v_x\circ J-j\circ d v_x)(Je_a)
        =-j(d v_x\circ J-j\circ d v_x)(e_a).
    \end{align*}
    Thus the defect has equal norms on $e_a$ and $Je_a$.
    Taking its squared norm on the full orthonormal basis gives
    \begin{align*}
        |d v_x\circ J-j\circ d v_x|^2
        =2\sum_{a=1}^m|d v_xJe_a-jd v_xe_a|_h^2.
    \end{align*}
    Combining the preceding identities and multiplying by
    $(d\vol_g)_x$ proves
    \eqref{Equation: Wirtinger energy defect} at almost every
    point. For a weakly pseudoholomorphic map $u$, the defect
    $d u\circ J-j\circ d u$ vanishes almost everywhere, giving
    \eqref{Equation: Energy identity}.
\end{proof}
The following identity is a standard chain rule for Sobolev pullbacks; see \cite[Lemma~4.1(2)]{HajlaszSchikorraTyson}. We recall its short proof here for completeness.
\begin{lemma}
    \label{Lemma: Sobolev pullbacks}
    Let $M$ be a smooth manifold, $N$ a compact smooth manifold,
    $v\in W^{1,2}_{\mathrm{loc}}(M,N)$ and $\eta$ a smooth 1-form on $N$. Then
    \begin{align} \label{Equation: Sobolev pullback identity}
        d(v^*\eta)=v^*(d\eta) \qquad\text{in the sense of distributions on }\,M.
    \end{align}
    Consequently, if $N=\cp^n$ and $d(v^*\omega)=0$, then
    $d(v^*\beta)=0$ for every smooth closed 2-form $\beta$ on
    $\cp^n$.
\end{lemma}
\begin{proof}
    Work in a relatively compact coordinate domain and fix a smooth
    embedding $N\subset\mathbb R^L$. Extend $\eta$ to a smooth
    compactly supported 1-form $\widetilde\eta$ on $\mathbb R^L$.
    Choose smooth Euclidean-valued maps $v_\ell$ converging strongly
    to $v$ in $W^{1,2}$ and, after taking a subsequence, almost
    everywhere. The boundedness of the coefficients of
    $\widetilde\eta$ and $d\widetilde\eta$ gives
    \begin{align*}
        v_\ell^*\widetilde\eta&\rightarrow v^*\eta\hspace{12.6mm}\text{ in }L^2,\\
        v_\ell^*(d\widetilde\eta)&\rightarrow v^*(d\eta)\qquad
        \text{ in }L^1.
    \end{align*}
    Indeed, in the second convergence, terms containing
    $d v_\ell-d v$ tend to zero by the Cauchy--Schwarz inequality.
    The remaining term tends to zero almost everywhere and is
    bounded by a constant times $|d v|^2$, so dominated convergence
    applies. The first convergence follows in the same way. Passing to the
    limit in the smooth pullback identity proves
    \eqref{Equation: Sobolev pullback identity}.

    Finally, $H^2_{\mathrm{dR}}(\cp^n;\mathbb R)$ is generated by
    $[\omega]$. Thus $\beta=c\omega+d\eta$ for some $c\in\mathbb R$
    and some smooth 1-form $\eta$. Therefore
    \begin{align*}
        d(v^*\beta)=c\,d(v^*\omega)+d^2(v^*\eta)=0.
    \end{align*}
	The statement follows.
\end{proof}

\subsection{Almost monotonicity and initial regularity}
\label{Subsection: Initial regularity}
The following almost-monotonicity formula can be proved by routine computations; see for instance \cite[Proposition~2.1 and Remark~2.2]{CaniatoRiviere}, where it is deduced from the almost-monotonicity formula for semicalibrated normal cycles in \cite[Proposition~1]{PumbergerRiviere}. In this argument, local strong approximability is used only to ensure that $d(u^*\omega)=0$, which is already part of \eqref{Equation: Assumptions on the map}; no closedness of $\Omega$ is required. By applying the formula in normal coordinates centred at each point, we can choose the geometric constants uniformly on compact subsets. The correction factor $e^{Ar}(1+Ar)$ appearing in the cited formula gives the exponential factor below after increasing the constant, since
\begin{align*}
    e^{A(R-r)}\frac{1+AR}{1+Ar}\leq e^{2A(R-r)}
    \qquad\forall\,0<r<R.
\end{align*}
\begin{lemma}
    \label{Lemma: Almost monotonicity of the energy}
    Let $K\subset M$ be compact. There exist $r_K>0$ and $C_K>0$,
    depending only on the geometric data near $K$, such that every
    map satisfying \eqref{Equation: Assumptions on the map} obeys
    \begin{align}\label{Equation: Energy almost monotonicity}
        \mathscr{E}_u(x,r)\leq e^{C_K(R-r)}\mathscr{E}_u(x,R) \qquad\forall\,x\in K,\,\forall\,0<r<R<r_K.
    \end{align}
    In particular,
    \begin{align*}
        \theta_u(x):=\lim_{r\rightarrow0}\mathscr{E}_u(x,r)
    \end{align*}
    exists and is finite at every point of $M$.
\end{lemma}
The $\varepsilon$-regularity statement below is a consequence of Wang's regularity theorem for weakly pseudoholomorphic maps with small Morrey norm \cite{WangJHolomorphic}, combined with Lemma~\ref{Lemma: Almost monotonicity of the energy}. See also \cite[Section~3]{CaniatoRiviere} for the corresponding regularity statement under local strong approximability. The same result can also be proved directly by using the Jacobian structure of the pseudoholomorphic equation and the compensated-compactness estimate of Coifman--Lions--Meyer--Semmes \cite{CoifmanLionsMeyerSemmes}, following the energy-decay argument in \cite[Lemma~3.1]{BellettiniTian}.
\begin{proposition}
    \label{Proposition: Initial regularity}
    Let $u$ satisfy \eqref{Equation: Assumptions on the map}. Then, $\Sing u$ is closed and
    \begin{align}\label{Equation: Initial singular set estimate}
        \Sing u=\{x\in M\mbox{ : }\theta_u(x)>0\}, \qquad \Hh^{2m-2}(\Sing u)=0.
    \end{align}
    Moreover, for every compact $K\subset M$, there exist
    $\varepsilon_K,r_K>0$, depending only on the geometric data
    near $K$, such that
    \begin{align}
        \mathscr{E}_u(x,r)<\varepsilon_K
        &\quad\Rightarrow\quad
        u\text{ has a smooth representative on }B_{r/2}(x),
        \label{Equation: Epsilon regularity}\\
        \theta_u(x)&\geq\varepsilon_K
        \qquad\text{for every }x\in \Sing u\cap K,
        \label{Equation: Singular density threshold}
    \end{align}
    whenever $x\in K$ and $0<r<r_K$ in
    \eqref{Equation: Epsilon regularity}.
\end{proposition}

\begin{proof}
    Apply Lemma~\ref{Lemma: Almost monotonicity of the energy} on a
    compact neighbourhood of $K$. After decreasing $r_K$, for every
    $y\in B_{r/2}(x)$ and $0<s\leq r/2$ we have
    \begin{align*}
        \mathscr{E}_u(y,s)
        \leq e^{C_Kr}\mathscr{E}_u(y,r/2)
        \leq 2^{2m-2}e^{C_Kr}\mathscr{E}_u(x,r),
    \end{align*}
    since $B_{r/2}(y)\subset B_r(x)$. Thus, if
    $\mathscr{E}_u(x,r)<\varepsilon_K$ with $\varepsilon_K$ sufficiently
    small, Wang's small-Morrey-norm theorem \cite{WangJHolomorphic} gives
    smoothness on $B_{r/2}(x)$. At a singular point $x\in K$ we therefore
    have $\mathscr{E}_u(x,r)\geq\varepsilon_K$ for
    every sufficiently small $r$, whereas at a regular point the
    energy density vanishes by smoothness. This proves the density
    characterisation and the lower bound on $\Sing u\cap K$.

    Lebesgue differentiation gives $\theta_u=0$ almost everywhere,
    so $\Sing u$ has zero Riemannian volume. Since $\mu_u$ is absolutely
    continuous with respect to volume, $\mu_u(\Sing u)=0$.
    The density lower bound on $\Sing u\cap K$ and the $5r$-covering lemma
    now give $\Hh^{2m-2}(\Sing u\cap K)=0$.
\end{proof}

\subsection{Semicalibrated cycles}
\label{Subsection: Semicalibrated cycles}

Let $(V,g_V)$ be a smooth Riemannian manifold. We denote by $\mathbf I_{q,\mathrm{loc}}(V)$ the set of locally integral $q$-currents, by $\|T\|$ the mass measure of a current $T\in\mathbf I_{q,\mathrm{loc}}(V)$, and by $\vec T$ its oriented unit approximate tangent $q$-vector. We say that $T\in\mathbf I_{q,\mathrm{loc}}(V)$ is a \emph{cycle} if  $\partial T=0$, i.e.
\begin{align*}
	\partial T(\eta)\coloneq T(d\eta)=0
\end{align*}
for every smooth and compactly supported $q-1$-form $\eta$ on $V$.  A smooth $q$-form $\varphi$ is a \emph{semicalibration} if its comass is at most one. We say that $\varphi$ semicalibrates $T\in\mathbf I_{q,\mathrm{loc}}(V)$ if
\begin{align}\label{Equation: Prelim semicalibration}
    \varphi_x(\vec T(x))=1 \qquad\text{ for }\|T\|\text{-a.e. } x\in V.
\end{align}
Note that, even when $\varphi$ is not closed, \eqref{Equation: Prelim semicalibration} implies that
\begin{align}
    T(\chi\varphi)=\int_V\chi\,d\|T\|
    \qquad\forall\,\chi\in C_c^\infty(V).
    \label{Equation: Semicalibrated mass identity}
\end{align}
Moreover, two currents $T_1,T_2\in\mathbf I_{q,\mathrm{loc}}(V)$ semicalibrated by the same form satisfy
\begin{align}\label{Equation: Additivity of semicalibrated mass}
    \|T_1+T_2\|=\|T_1\|+\|T_2\|,
\end{align}
and their sum has the same semicalibration. Indeed, the triangle inequality gives $\|T_1+T_2\|\leq\|T_1\|+\|T_2\|$. Testing the sum on $\chi\varphi$ for nonnegative $\chi$ gives the reverse inequality by the comass bound. Equality also gives the semicalibration of the sum.

The following statement is a standard consequence of the almost-monotonicity formula for semicalibrated cycles; see \cite[Proposition~1]{PumbergerRiviere} and \cite[Lemma~0.4 and the proof of Lemma~1.3]{Spolaor}. The constants depend only on the local geometry and the fixed semicalibration, and can therefore be chosen uniformly on compact subsets, independently of the current. The domination estimate follows by comparing the mass ratio at each smaller radius with that at a fixed radius $R$ such that $B_R(x)\subset W$ for every $x\in K$. We record this uniformity explicitly because it will be used when averaging a family of locally integral cycles.
\begin{lemma}
    \label{Lemma: Uniform monotonicity for semicalibrated cycles}
    Fix a smooth semicalibration $\varphi$ on $V$ and a compact set $K\subset W\Subset V$, where $W$ is open. There exist $r_0,C,A>0$, depending only on $K,W,g_V,\varphi$, and $q$, such that every $\varphi$-semicalibrated locally integral cycle $T\in\mathbf I_{q,\mathrm{loc}}(V)$ satisfies
    \begin{align*}
        r\mapsto e^{Cr}r^{-q}\|T\|(B_r(x)) \quad\text{ is nondecreasing on }(0,r_0)
    \end{align*}
    for every $x\in K$. In particular, the limit
    \begin{align*}
        \theta_T(x):=\lim_{r\rightarrow0}r^{-q}\|T\|(B_r(x))
    \end{align*}
    exists and is finite at every point in $V$ and vanishes outside $\spt T$. Moreover, the mass ratios satisfy
    \begin{align} \label{Equation: Uniform domination of current densities}
        r^{-q}\|T\|(B_r(x))\leq A\|T\|(W) \qquad\forall\,x\in K,\,\forall\,0<r<r_0.
    \end{align}
\end{lemma}

For a locally integral $q$-cycle $T$, we call a point $x\in\spt T$ \textit{regular} if there exist a neighbourhood $U$ of $x$, a connected embedded oriented $C^1$ $q$-dimensional submanifold $\Sigma\subset U$ without boundary in $U$, and an integer $Q\geq1$ such that
\begin{align*}
    \Sigma=\spt T\cap U, \qquad T=Q\llbracket\Sigma\rrbracket\quad\text{ as currents on }U.
\end{align*}
We denote the set of regular points by $\Reg\,T$ and write $\Sing T:=\spt T\smallsetminus\Reg\,T$. The regular set of $T$ is relatively open in its support; hence, $\Sing T$ is closed in $V$.

We use the following regularity theorem for semicalibrated locally integral cycles.
\begin{theorem}
    \label{Theorem: Regularity of semicalibrated cycles}
    Let $q\geq2$ and let $T\in\mathbf I_{q,\mathrm{loc}}(V)$ be a locally integral cycle semicalibrated by a smooth form on the smooth Riemannian manifold $(V,g_V)$. Then
    \begin{align*}
        \dim_{\Hh}\Sing T\leq q-2.
    \end{align*}
    If $q=2$, then $\Sing T$ is locally finite in $V$.
\end{theorem}
The dimension estimate is due to Spolaor \cite[Theorem~0.2]{Spolaor}. For a positive $(1,1)$ integral cycle in an almost complex manifold admitting locally compatible symplectic forms, the 2-dimensional assertion is the theorem of Rivi\`ere--Tian \cite[p.~743]{RiviereTianCurrents}. For an arbitrary semicalibrated integral $2$-cycle, the same conclusion follows from De Lellis--Spadaro--Spolaor \cite[Theorem~0.2, case~(b)]{DeLellisSpadaroSpolaor}. The Riemannian formulation follows by a smooth local isometric embedding and the extension of semicalibrations in \cite[Lemma~0.3]{Spolaor}.
\section{Averaged semicalibrated cycles}
\label{Section: Averaged semicalibrated cycles}

In this section we establish an averaging principle for families of
semicalibrated cycles. The current dimension is denoted by $q$;
in the application to pseudoholomorphic maps, we will take
$q=2m-2$. The argument uses the regularity of the sum of two
cycles and the relation between the density of an averaged mass
measure and the probability that a point belongs to the support
of a member of the family.

Let $(V,g_V)$ be a smooth Riemannian manifold. Let
$(A,\mathcal A)$ be a standard Borel space, meaning that
$A$ admits a complete separable metric whose Borel
$\sigma$-algebra is $\mathcal A$, and let $\lambda$ be a
probability measure on $(A,\mathcal A)$.

Consider a family $(T_a)_{a\in A}$ of locally integral
currents on $V$. For each $a\in A$, the mass measure
$\|T_a\|$ is a locally finite Borel measure on $V$.
We say that these mass measures form a
\emph{measurable kernel} if, for every fixed Borel set
$B\subset V$, the function
\begin{align*}
    A\rightarrow[0,+\infty],
    \qquad a\mapsto\|T_a\|(B),
\end{align*}
is $\mathcal A$-measurable.

We call the family $(T_a)_{a\in A}$ \emph{measurable}
if its mass measures form a measurable kernel and,
for every compactly supported smooth test form $\zeta$
of the appropriate degree, the function
$a\mapsto T_a(\zeta)$ is $\mathcal A$-measurable.

All balls below are geodesic balls for $g_V$, and all
densities are unnormalised as in
Section~\ref{Section: Preliminaries}.

\begin{theorem}
    \label{Theorem: Averaged cycles}
    Let $2\leq q\leq\dim V$, and let $(T_a)_{a\in A}$ be a
    measurable family of locally integral $q$-cycles in $V$,
    semicalibrated by a fixed smooth $q$-form $\varphi$.
    Assume that
    \begin{align}
        \int_A\|T_a\|(W)\,d\lambda(a)<\infty
        \qquad\text{for every open }W\Subset V,
        \label{Equation: Integrable local masses}
    \end{align}
    and that
    \begin{align}
        \|T_a\|\perp\|T_b\|
        \qquad\text{for }
        (\lambda\otimes\lambda)\text{-a.e. }(a,b)\in A\times A.
        \label{Equation: Almost every pair is mutually singular}
    \end{align}
    Define the averaged measure by
    \begin{align}
        \mu(B):=\int_A\|T_a\|(B)\,d\lambda(a)
        \qquad\text{for every Borel set }B\subset V.
        \label{Equation: Averaged mass measure}
    \end{align}
    Then $\mu$ is a Radon measure and the limit
    \begin{align*}
        \theta_\mu(x):=
        \lim_{r\rightarrow0}\frac{\mu(B_r(x))}{r^q}
    \end{align*}
    exists and is finite at every $x\in V$. Its positive-density
    set satisfies
    \begin{align}
        \dim_{\Hh}E\leq q-2,
        \qquad
        E:=\{x\in V\mbox{ : }\theta_\mu(x)>0\}.
        \label{Equation: Dimension of the positive density set}
    \end{align}
    Moreover, if $q=2$, then for every compact set $K\subset V$
    and every $\varepsilon>0$, the set
    \begin{align}
        E_\varepsilon\cap K,
        \qquad
        E_\varepsilon:=\{x\in V\mbox{ : }\theta_\mu(x)\geq\varepsilon\},
        \label{Equation: Positive density superlevel sets}
    \end{align}
    is finite.
\end{theorem}

The assumptions that the currents are locally integral cycles and have the
specified semicalibration may equivalently be imposed outside a fixed
parameter set of measure zero: redefine the currents to be zero on
that set. We work with this representative of the family throughout.
In particular, the density of every individual current exists at
every point by
Lemma~\ref{Lemma: Uniform monotonicity for semicalibrated cycles}.

\subsection{Intersections of supports}
\label{Subsection: Intersections of supports}

\begin{lemma}
    \label{Lemma: Pairwise intersection}
    Let $T_1,T_2\in\mathbf I_{q,\mathrm{loc}}(V)$ be cycles
    semicalibrated by the same smooth form $\varphi$.
    If $\|T_1\|\perp\|T_2\|$, then
    \begin{align}
        \spt T_1\cap\spt T_2
        \subset\Sing(T_1+T_2).
        \label{Equation: Intersections are singular for the sum}
    \end{align}
    Consequently, if $q\geq2$, then
    \begin{align}
        \dim_{\Hh}(\spt T_1\cap\spt T_2)\leq q-2.
        \label{Equation: Dimension of a pairwise intersection}
    \end{align}
    If $q=2$, then $\spt T_1\cap\spt T_2$ is locally finite in $V$.
\end{lemma}
\begin{proof}
    The common semicalibration gives
    \begin{align*}
        \|T_1+T_2\|=\|T_1\|+\|T_2\|,
    \end{align*}
    by \eqref{Equation: Additivity of semicalibrated mass}.
    Moreover, $T_1+T_2$ is a locally integral cycle semicalibrated by
    $\varphi$.

    Suppose that $x\in\spt T_1\cap\spt T_2$ is a regular point
    of $T_1+T_2$. By regularity, after shrinking a neighbourhood
    $U$ of $x$, there is a connected embedded oriented
    $C^1$ $q$-dimensional submanifold
    \begin{align*}
        \Sigma=\spt(T_1+T_2)\cap U,
        \qquad x\in\Sigma,
    \end{align*}
    such that $T_1+T_2=Q\llbracket\Sigma\rrbracket$ as currents in $U$, for some
    integer $Q\geq1$.
    Mass additivity implies that both $T_1$ and $T_2$ are
    supported on $\Sigma$ in $U$. Since their boundaries vanish
    in $U$, the constancy theorem on the connected sheet gives
    \begin{align*}
        T_1=Q_1\llbracket\Sigma\rrbracket,
        \qquad T_2=Q_2\llbracket\Sigma\rrbracket
        \qquad\text{as currents in }U,
    \end{align*}
    where $Q_1,Q_2$ are integers. Since $\varphi$ semicalibrates both
    currents and agrees with the volume form on the oriented sheet $\Sigma$, these
    multiplicities are nonnegative.

    Both multiplicities are positive. Indeed, if $Q_i=0$, then
    $T_i$ vanishes on the neighbourhood $U$ of $x$, contradicting
    $x\in\spt T_i$. Thus the restrictions of $\|T_1\|$ and
    $\|T_2\|$ to $U$ are positive multiples of the same nonzero
    measure $\Hh^q\!\llcorner\Sigma$.
    This contradicts their mutual singularity and proves
    \eqref{Equation: Intersections are singular for the sum}.

    Apply Theorem~\ref{Theorem: Regularity of semicalibrated cycles}
    to $T_1+T_2$ to obtain
    \eqref{Equation: Dimension of a pairwise intersection}.
    When $q=2$, the singular set of the sum is closed and
    discrete in $V$, and hence finite on every compact subset
    of $V$. The same conclusion holds for the support
    intersection by
    \eqref{Equation: Intersections are singular for the sum}.
\end{proof}

\subsection{Density averaging and incidence probabilities}
\label{Subsection: Density averaging and incidence}

We first record the measurability facts needed below. Write
$\rho$ for the Riemannian distance on $V$. For each fixed $r>0$,
the kernel assumption implies that
\begin{align*}
    (a,x)\mapsto\|T_a\|(B_r(x))
    =\int_V\mathbf{1}_{\{\rho(x,y)<r\}}\,d\|T_a\|(y)
\end{align*}
is jointly measurable. Indeed, the integrand is Borel measurable
in $(x,y)$, so the measurability theorem for integration against a
measure kernel applies. In particular, the incidence set
\begin{align}
    \mathcal I:=\{(a,x)\in A\times V\mbox{ : }x\in\spt T_a\}=\bigcap_{j=1}^\infty
    \{(a,x)\in A\times V\mbox{ : }\|T_a\|(B_{1/j}(x))>0\}
    \label{Equation: Measurability of support incidence}
\end{align}
is measurable. We may therefore define the measurable function
\begin{align}
    \mathfrak p(x)
    :=\lambda\big(\{a\in A\mbox{ : }x\in\spt T_a\}\big)
    =\int_A\mathbf{1}_{\mathcal I}(a,x)\,d\lambda(a),
    \label{Equation: Incidence probability}
\end{align}
which takes values in $[0,1]$.

The definition \eqref{Equation: Averaged mass measure} gives a
countably additive Borel measure by monotone convergence.
Assumption~\eqref{Equation: Integrable local masses} makes it
locally finite, and hence Radon on the manifold $V$.

\begin{lemma}
    \label{Lemma: Density and incidence}
    Under the common-semicalibration and local-mass assumptions of
    Theorem~\ref{Theorem: Averaged cycles}, let
    \begin{align*}
        \theta_a(x):=
        \lim_{r\rightarrow0}r^{-q}\|T_a\|(B_r(x)).
    \end{align*}
    Then the density of $\mu$ exists and is finite at every point,
    and
    \begin{align}
        \theta_\mu(x)=\int_A\theta_a(x)\,d\lambda(a)
        \qquad\text{for every }x\in V.
        \label{Equation: Density disintegration}
    \end{align}
    Moreover,
    \begin{align}
        \theta_\mu(x)>0\quad\Rightarrow\quad
        \mathfrak p(x)>0.
        \label{Equation: Positive density implies positive incidence}
    \end{align}
    For every compact $K\subset V$ and every $\varepsilon>0$,
    there is a constant $\delta_{K,\varepsilon}>0$ such that
    \begin{align}
        \mathfrak p(x)\geq\delta_{K,\varepsilon}
        \qquad\text{for every }x\in E_\varepsilon\cap K.
        \label{Equation: Uniform positive incidence}
    \end{align}
    The constant may depend on the family of currents.
\end{lemma}
\begin{proof}
    Fix a compact set $K\subset V$ and choose an open set
    $W\Subset V$ containing $K$.
    Lemma~\ref{Lemma: Uniform monotonicity for semicalibrated cycles}
    gives constants $r_0,C>0$, independent of $a$, such that
    \begin{align}
        r^{-q}\|T_a\|(B_r(x))
        \leq C\|T_a\|(W)=:G_K(a)
        \qquad\forall\,x\in K,\,\forall\,0<r<r_0.
        \label{Equation: Integrable density domination}
    \end{align}
    By \eqref{Equation: Integrable local masses}, $G_K\in L^1(A)$.
    The same monotonicity lemma gives existence of $\theta_a(x)$
    for every $a$ and every $x$, and
    $0\leq\theta_a(x)\leq G_K(a)$ for $x\in K$.
    Thus, for every fixed $x\in K$, dominated convergence gives
    \begin{align*}
        \lim_{r\rightarrow0}\frac{\mu(B_r(x))}{r^q}=\lim_{r\rightarrow0}\int_A
        \frac{\|T_a\|(B_r(x))}{r^q}\,d\lambda(a)=\int_A\theta_a(x)\,d\lambda(a)<\infty.
    \end{align*}
    Since the same function $G_K$ dominates the mass ratios for every
    $x\in K$, a compact exhaustion of $V$ gives
    \eqref{Equation: Density disintegration} at every point of $V$.

    The functions $(a,x)\mapsto\theta_a(x)$ and
    $x\mapsto\theta_\mu(x)$ are measurable because they are pointwise
    limits of the corresponding ball-mass ratios along $r=1/j$.
    Consequently, $E$ and $E_\varepsilon$ are Borel sets.

    If $x\notin\spt T_a$, then $T_a$ has zero mass on some
    neighbourhood of $x$, so $\theta_a(x)=0$.
    Therefore \eqref{Equation: Density disintegration} gives
    \begin{align}
        \theta_\mu(x)
        =\int_{\{a\mbox{ : }x\in\spt T_a\}}
        \theta_a(x)\,d\lambda(a).
        \label{Equation: Density is supported on incident parameters}
    \end{align}
    If $\mathfrak p(x)=0$, the integration domain is $\lambda$-null, so the integral vanishes. This proves
    \eqref{Equation: Positive density implies positive incidence}.

    Finally, fix $\varepsilon>0$ and choose $L>\varepsilon$ so
    large that
    \begin{align*}
        \int_{\{G_K>L\}}G_K\,d\lambda<\frac{\varepsilon}{2}.
    \end{align*}
    Such a choice is possible because $G_K\in L^1(A)$.
    For $x\in E_\varepsilon\cap K$, split the integral in
    \eqref{Equation: Density is supported on incident parameters}
    according to whether $G_K\leq L$ or $G_K>L$. Using
    $\theta_a(x)\leq G_K(a)$, we obtain
    \begin{align*}
        \varepsilon
        \leq\theta_\mu(x)
        \leq L\mathfrak p(x)
        +\int_{\{G_K>L\}}G_K\,d\lambda
        <L\mathfrak p(x)+\frac{\varepsilon}{2}.
    \end{align*}
    Thus \eqref{Equation: Uniform positive incidence} holds with
    $\delta_{K,\varepsilon}:=\varepsilon/(2L)$.
\end{proof}

\subsection{Proof of the averaging theorem}
\label{Subsection: Proof of the averaging theorem}

\begin{proof}[Proof of Theorem~\ref{Theorem: Averaged cycles}]
    The Radon property of $\mu$ and the existence and finiteness
    of its density have already been established. We prove the
    dimension estimate first, and then the 2-dimensional
    finiteness assertion.

    Suppose, towards a contradiction, that $\dim_{\Hh}E>q-2$.
    By a countable compact exhaustion and a countable coordinate
    cover, there are a compact set $K$ contained in a coordinate
    neighbourhood and a number $s$ such that
    \begin{align*}
        q-2<s<\dim_{\Hh}(E\cap K).
    \end{align*}
    Apply Frostman's lemma in its version for Borel subsets of
    Euclidean space, using a coordinate chart on a neighbourhood
    of $K$; see \cite{Mattila}. Since the Riemannian and Euclidean
    distances are locally comparable there, we obtain a nonzero
    finite Radon measure $\nu$, supported on a compact subset of
    $E\cap K$, and constants $C_\nu,r_\nu>0$ such that
    \begin{align}
        \nu(B_r(x))\leq C_\nu r^s
        \qquad\forall\,x\in V,\,\forall\, 0<r<r_\nu.
        \label{Equation: Frostman measure for the density set}
    \end{align}
    In particular, $\nu$ vanishes on every Borel set of zero
    $\Hh^s$-measure. Indeed, covering such a set by sufficiently
    small balls and applying
    \eqref{Equation: Frostman measure for the density set}
    bounds its $\nu$-measure by an arbitrarily small multiple
    of $C_\nu$.

    For a.e. $(a,b)$, assumption
    \eqref{Equation: Almost every pair is mutually singular}
    and Lemma~\ref{Lemma: Pairwise intersection} imply
    \begin{align*}
        \dim_{\Hh}(\spt T_a\cap\spt T_b)\leq q-2<s.
    \end{align*}
    Consequently,
    \begin{align*}
        \nu(\spt T_a\cap\spt T_b)=0
        \qquad\text{for a.e. }(a,b).
    \end{align*}
    By \eqref{Equation: Measurability of support incidence},
    the function
    $(a,b,x)\mapsto\mathbf{1}_{\mathcal I}(a,x)
    \mathbf{1}_{\mathcal I}(b,x)$ is jointly measurable.
    Tonelli's theorem therefore gives
    \begin{align*}
        0
        &=\int_{A\times A}
        \nu(\spt T_a\cap\spt T_b)
        \,d(\lambda\otimes\lambda)(a,b)\\
        &=\int_V\int_A\int_A
        \mathbf{1}_{\mathcal I}(a,x)
        \mathbf{1}_{\mathcal I}(b,x)
        \,d\lambda(a)\,d\lambda(b)\,d\nu(x)\\
        &=\int_V\mathfrak p(x)^2\,d\nu(x).
    \end{align*}
    On the other hand,
    \eqref{Equation: Positive density implies positive incidence}
    gives $\mathfrak p(x)>0$ at every point of $E$, and
    $0<\nu(E)=\nu(V)$. Hence the last integral is strictly
    positive, a contradiction. This proves
    \eqref{Equation: Dimension of the positive density set}.

    Now assume $q=2$, fix a compact set $K\subset V$ and
    $\varepsilon>0$, and suppose that $E_\varepsilon\cap K$
    contains a sequence of pairwise distinct points $(x_j)$.
    By Lemma~\ref{Lemma: Density and incidence},
    \begin{align*}
        \mathfrak p(x_j)\geq\delta_{K,\varepsilon}>0
        \qquad\text{for every }j.
    \end{align*}
    Define measurable functions on $A\times A$ by
    \begin{align*}
        f_j(a,b):=
        \mathbf{1}_{\mathcal I}(a,x_j)
        \mathbf{1}_{\mathcal I}(b,x_j).
    \end{align*}
    For a.e. $(a,b)$,
    Lemma~\ref{Lemma: Pairwise intersection} gives
    \begin{align*}
        \#(\spt T_a\cap\spt T_b\cap K)<\infty.
    \end{align*}
    For each such pair, only finitely many of the distinct points
    $x_j$ belong to the intersection. Hence $f_j\to0$ almost
    everywhere. Since $0\leq f_j\leq1$ and
    $\lambda\otimes\lambda$ is a probability measure, dominated
    convergence yields
    \begin{align*}
        \int_{A\times A}f_j\,d(\lambda\otimes\lambda)
        \rightarrow0.
    \end{align*}
    However, the product structure of the parameter measure gives
    \begin{align*}
        \int_{A\times A}f_j\,d(\lambda\otimes\lambda)
        =\mathfrak p(x_j)^2
        \geq\delta_{K,\varepsilon}^2>0,
    \end{align*}
    which is a contradiction. Hence $E_\varepsilon\cap K$ is finite.
\end{proof}

\begin{remark}
    \label{Remark: Scope of the averaging theorem}
    For $q=2$, the finiteness conclusion uses a fixed positive
    density threshold. Thus, if a set $S\subset V$ satisfies
    $\theta_\mu\geq\varepsilon_K>0$ on $S\cap K$ for every compact
    $K\subset V$, then $S$ is locally finite. Proposition~\ref{Proposition: Initial regularity}
    supplies this density bound for the singular set of a map.
    For $q>2$, the theorem gives a Hausdorff-dimension estimate;
    the argument supplies no local bound for
    $\Hh^{q-2}(E)$.
\end{remark}
\section{Inverse-hyperplane cycles}
\label{Section: Inverse hyperplane cycles}

We now construct the family of semicalibrated integral cycles to which
Theorem~\ref{Theorem: Averaged cycles} will be applied. We retain the
notation and assumptions of Section~\ref{Section: Preliminaries}.
In particular, the Fubini--Study form $\omega$ has the
projective-line normalisation fixed in the introduction, and
\begin{align*}
    \Phi=\frac{\Omega^{m-1}}{(m-1)!},
    \qquad
    \mu_u=\frac12|du|^2\,d\vol_g.
\end{align*}
The construction starts with transverse inverse images of hyperplanes
on $\Reg\,u$. We will prove that for almost every hyperplane
its inverse image defines an integral cycle on all of $M$.
For this purpose, we need an integral-geometric identity with an
arbitrary smooth weight on the space of hyperplanes.

\subsection{Weighted incidence forms}
\label{Subsection: Weighted incidence forms}

Let $A:=(\cp^n)^*$ be the dual projective space, whose points
parametrise the complex projective hyperplanes in $\cp^n$.
We equip $A$ with its Borel $\sigma$-algebra and its unitary-invariant
probability measure $\lambda$. This measure is induced by a smooth
positive volume form, which we also denote by $d\lambda$.
Consider the incidence manifold
\begin{align*}
    \mathcal Z:=\{(y,H)\in\cp^n\times A \mbox{ : } y\in H\},
\end{align*}
with the projections
\begin{align*}
    \pi:\mathcal Z\rightarrow\cp^n,
    \qquad
    \sigma:\mathcal Z\rightarrow A.
\end{align*}
Both projections are proper holomorphic submersions. Their fibres
are projective spaces of complex dimension $n-1$: the fibre of
$\sigma$ over $H$ is $H$ itself, whereas the fibre of $\pi$ over
$y$ is the space of hyperplanes containing $y$.
We give all these fibres their complex orientations.

We write $\pi_*$ for integration of differential forms along the
fibres of $\pi$. We give its pointwise definition in the degree
needed below. For $y\in\cp^n$, set
\begin{align*}
    \mathcal Z_y:=\pi^{-1}(y)
    =\{(y,H)\in\mathcal Z \mbox{ : } y\in H\}.
\end{align*}
This is a compact manifold of real dimension $2n-2$, equipped
with its complex orientation. If $\alpha$ is a smooth real
$2n$-form on $\mathcal Z$, then $\pi_*\alpha$ is the smooth
real 2-form on $\cp^n$ defined as follows. Fix $y\in\cp^n$ and $X,Y\in T_y\cp^n$. For each
$p\in\mathcal Z_y$, choose vectors
$\widetilde X_p,\widetilde Y_p\in T_p\mathcal Z$ satisfying
\begin{align*}
    d\pi_p(\widetilde X_p)=X,
    \qquad
    d\pi_p(\widetilde Y_p)=Y.
\end{align*}
Such lifts exist because $\pi$ is a submersion, and they can
be chosen smoothly along $\mathcal Z_y$. For instance, after
choosing an auxiliary Riemannian metric on $\mathcal Z$, one
can take the unique lifts orthogonal to $\ker d\pi_p$.
Define a smooth $(2n-2)$-form $\eta_{\alpha,y;X,Y}$ on
$\mathcal Z_y$ by
\begin{align*}
    (\eta_{\alpha,y;X,Y})_p(V_1,\ldots,V_{2n-2}):=
    \alpha_p(\widetilde X_p,\widetilde Y_p,
    V_1,\ldots,V_{2n-2}),
    \qquad
    \mbox{ for }\,V_1,\ldots,V_{2n-2}\in T_p\mathcal Z_y.
\end{align*}
Since this is a top-degree form on the oriented fibre, its
integral is well defined. We set
\begin{align}
    (\pi_*\alpha)_y(X,Y)
    :=\int_{\mathcal Z_y}\eta_{\alpha,y;X,Y}.
    \label{Equation: Pointwise fibre integration}
\end{align}

This definition is independent of the chosen lifts. Indeed,
two lifts of $X$ differ by a vector $W\in\ker d\pi_p
=T_p\mathcal Z_y$. Replacing $\widetilde X_p$ by
$\widetilde X_p+W$ changes the integrand by
\begin{align*}
    \alpha_p(W,\widetilde Y_p,V_1,\ldots,V_{2n-2}).
\end{align*}
The $2n-1$ vectors $W,V_1,\ldots,V_{2n-2}$ belong to the
$(2n-2)$-dimensional space $T_p\mathcal Z_y$, so this expression
vanishes by alternation. The same argument applies to the lift
of $Y$. In particular, the auxiliary metric has no effect on
$\pi_*\alpha$. The construction is alternating and bilinear
in $X,Y$, and linear in $\alpha$.

To see that $\pi_*\alpha$ is smooth, take an oriented product chart with coordinates
\begin{align*}
    (y^1,\ldots,y^{2n},t^1,\ldots,t^{2n-2}),
    \qquad
    \pi(y,t)=y,
\end{align*}
and write
\begin{align*}
    \alpha
    &=
    \sum_{1\leq i<j\leq2n}
    a_{ij}(y,t)\,dy^i\wedge dy^j
    \wedge dt^1\wedge\cdots\wedge dt^{2n-2}
    +\alpha_{\mathrm{rem}},
\end{align*}
where every term of $\alpha_{\mathrm{rem}}$ contains fewer
than $2n-2$ fibre differentials. If $\alpha$ is compactly
supported in this chart, then
\begin{align*}
    \pi_*\alpha
    =
    \sum_{1\leq i<j\leq2n}
    \left(
        \int_{\mathbb R^{2n-2}}
        a_{ij}(y,t)\,dt^1\cdots dt^{2n-2}
    \right)dy^i\wedge dy^j,
\end{align*}
where the coefficients are extended by zero outside the
coordinate chart. The terms in $\alpha_{\mathrm{rem}}$
contribute zero. A partition of unity reduces the general
case to this formula. Differentiating the coefficients under
the integral sign now gives the desired smoothness.

The complex orientation of $\mathcal Z$ agrees with the
orientation obtained from the complex orientations of the
base and the fibres. Fubini's theorem in these product charts
therefore gives
\begin{align}
    \int_{\cp^n}(\pi_*\alpha)\wedge\gamma
    =
    \int_{\mathcal Z}\alpha\wedge\pi^*\gamma
    \qquad\forall\,\gamma\in\Omega^{2n-2}(\cp^n).
    \label{Equation: fibre integration duality}
\end{align}
The same construction applies in other degrees, inserting
one lifted vector for each argument of the resulting form;
thus $\pi_*$ lowers the degree by $2n-2$. Since the fibres
have no boundary and their dimension is even, Stokes'
theorem gives
\begin{align}\label{Equation: commutation pi_* and differential}
    d(\pi_*\alpha)=\pi_*(d\alpha).
\end{align}
In particular, for $f\in C^\infty(A;\mathbb R)$ and
$\alpha=\sigma^*(f\,d\lambda)$, the fibre form appearing in
\eqref{Equation: Pointwise fibre integration} is explicitly
\begin{align*}
    &(\eta_{\alpha,y;X,Y})_p(V_1,\ldots,V_{2n-2})\\
    &\qquad=
    f(H)\,(d\lambda)_H\bigl(
        d\sigma_p(\widetilde X_p),
        d\sigma_p(\widetilde Y_p),
        d\sigma_p(V_1),\ldots,d\sigma_p(V_{2n-2})
    \bigr),
    \qquad p=(y,H)\in\mathcal Z_y.
\end{align*}
Integrating this top-degree form over $\mathcal Z_y$ defines
the value of $\pi_*\sigma^*(f\,d\lambda)$ on $X,Y$.
When $n=1$, each fibre is a positively oriented point, and
the fibre integral is simply evaluation at that point.

\begin{lemma}
    \label{Lemma: Weighted incidence forms}
    For $f\in C^\infty(A;\mathbb R)$, define
    \begin{align}
        \beta_f:=\pi_*\sigma^*(f\,d\lambda)\in\Omega^2(\cp^n).
        \label{Equation: Definition of weighted incidence forms}
    \end{align}
    Then
    \begin{align}
        \int_{\cp^n}\beta_f\wedge\gamma
        =\int_A f(H)\left(\int_H\gamma\right)d\lambda(H) \qquad\forall\gamma\in\Omega^{2n-2}(\cp^n).
        \label{Equation: Weighted incidence identity on the target}
    \end{align}
    Moreover,
    \begin{align}
        \beta_1=\omega
        \qquad\mbox{ and }\qquad
        \beta_f=c_f\omega+d\eta_f
        \quad\mbox{ with }\, 
        c_f:=\int_A f\,d\lambda,
        \label{Equation: Cohomology of weighted incidence forms}
    \end{align}
    for some smooth real 1-form $\eta_f$ on $\cp^n$.
\end{lemma}
\begin{proof}
    Since $f\,d\lambda$ has top degree on $A$, it is closed.
    By \eqref{Equation: commutation pi_* and differential}, we therefore have
    \begin{align*}
        d\beta_f
        =\pi_*d\sigma^*(f\,d\lambda)
        =\pi_*\sigma^*d(f\,d\lambda)=0.
    \end{align*}
    Fix $\gamma\in\Omega^{2n-2}(\cp^n)$. By the definition
    of $\beta_f$ and \eqref{Equation: fibre integration duality},
    \begin{align*}
        \int_{\cp^n}\beta_f\wedge\gamma
        =
        \int_{\mathcal Z}
        \sigma^*(f\,d\lambda)\wedge\pi^*\gamma.
    \end{align*}
    To evaluate the last integral, denote by $\sigma_*$
    integration along the complex-oriented fibres of $\sigma$,
    defined by the same construction as $\pi_*$.
    Since $\pi^*\gamma$ has degree $2n-2$, equal to the real
    dimension of these fibres, $\sigma_*(\pi^*\gamma)$ is a
    smooth real-valued function on $A$. For every $H\in A$,
    the restriction
    \begin{align*}
        \pi|_{\sigma^{-1}(H)}:\sigma^{-1}(H)&\rightarrow H\\
        (y,H)&\mapsto y
    \end{align*}
    is a holomorphic diffeomorphism and therefore preserves
    the complex orientations. Consequently,
    \begin{align*}
        \bigl(\sigma_*(\pi^*\gamma)\bigr)(H)
        =
        \int_{\sigma^{-1}(H)}\pi^*\gamma
        =
        \int_H\gamma.
    \end{align*}
    Applying Fubini's theorem for integration along the fibres
    of $\sigma$, with the same orientation convention as in
    \eqref{Equation: fibre integration duality}, now gives
    \begin{align*}
        \int_{\mathcal Z}
        \sigma^*(f\,d\lambda)\wedge\pi^*\gamma
        &=
        \int_A f(H)\,
        \bigl(\sigma_*(\pi^*\gamma)\bigr)(H)\,d\lambda(H)\\
        &=
        \int_A f(H)\left(\int_H\gamma\right)d\lambda(H).
    \end{align*}
    Indeed, in product charts for $\sigma$, the pulled-back
    top-degree form $\sigma^*(f\,d\lambda)$ contains all the
    base differentials. Its wedge product with $\pi^*\gamma$
    therefore retains only the component of $\pi^*\gamma$
    of top degree in the fibre variables. Ordinary Fubini
    integration in these coordinates, followed by a partition
    of unity, gives the first equality. The complex
    orientations of the total space, base, and fibres are
    compatible, so no additional sign occurs.
    This proves
    \eqref{Equation: Weighted incidence identity on the target}.

    The incidence relation and $\lambda$ are invariant under
    projective unitary transformations. Hence $\beta_1$ is an
    invariant real 2-form on $\cp^n$. The invariant real
    2-forms are precisely the constant multiples of $\omega$,
    as follows by considering the isotropy action on a tangent
    space. Our normalisation of $\omega$ gives
    \begin{align*}
        \int_{\cp^n}\omega^n=1,
        \qquad
        \int_H\omega^{n-1}=1.
    \end{align*}
    Taking $f=1$ and $\gamma=\omega^{n-1}$ in
    \eqref{Equation: Weighted incidence identity on the target}
    therefore shows that this constant multiple is one.

    Finally, since $H^2_{\mathrm{dR}}(\cp^n;\mathbb R)$ is
    generated by $[\omega]$, the closed form $\beta_f$ represents
    a scalar multiple of $[\omega]$. Pairing with
    $[\omega]^{n-1}$ determines this scalar:
    \begin{align*}
        \int_{\cp^n}\beta_f\wedge\omega^{n-1}
        =\int_A f(H)\,d\lambda(H)=c_f.
    \end{align*}
    Thus $\beta_f-c_f\omega$ is exact, which gives the asserted
    smooth 1-form $\eta_f$.
\end{proof}

\subsection{The identity on the regular set}
\label{Subsection: Smooth inverse hyperplanes}

We first establish the pullback identity for a smooth
pseudoholomorphic map. This also supplies the measurability
needed to construct the currents associated with $u$.

\begin{lemma}
    \label{Lemma: Smooth inverse hyperplanes}
    Let $V\subset M$ be open and let $v:V\rightarrow\cp^n$ be smooth
    and pseudoholomorphic. There is a Borel set $A_v\subset A$
    with $\lambda(A\smallsetminus A_v)=0$ such that, for every
    $H\in A_v$, the inverse image
    \begin{align*}
        F_H^v:=\{x\in V:v(x)\in H\}
    \end{align*}
    is either empty or a smooth embedded submanifold of real
    dimension $2m-2$, closed in $V$ and without boundary in $V$.
    Its tangent spaces are $J$-invariant, and we give it the
    complex orientation.

    Set $\llbracket F_H^v\rrbracket=0$ for $H\notin A_v$. The current evaluations
    $H\mapsto\llbracket F_H^v\rrbracket(\zeta)$ are Borel measurable for every
    smooth compactly supported $(2m-2)$-form $\zeta$ on $V$.
    The measures
    \begin{align*}
        \nu_H^v:=
        \begin{cases}
            \Hh^{2m-2}\!\llcorner F_H^v & \mbox{ for } H\in A_v,\\
            0 & \mbox{ for } H\notin A_v 
        \end{cases}
    \end{align*}
    form a Borel kernel on $V$.
    Moreover, for every such $\zeta$ and every
    $f\in C^\infty(A;\mathbb R)$,
    \begin{align}
        \int_A f(H)\llbracket F_H^v\rrbracket(\zeta)\,d\lambda(H)
        =\int_V v^*\beta_f\wedge\zeta.
        \label{Equation: Weighted identity for smooth maps}
    \end{align}
    The integral on the left is absolutely convergent.
\end{lemma}
\begin{proof}
    Form the pullback incidence manifold
    \begin{align*}
        \mathcal Z_v
        :=\{(x,H)\in V\times A:v(x)\in H\},
    \end{align*}
    and denote its projections by $\pi_v:\mathcal Z_v\rightarrow V$
    and $\sigma_v:\mathcal Z_v\rightarrow A$.
    Since $\pi:\mathcal Z\rightarrow\cp^n$ is a submersion,
    $\mathcal Z_v$ is smooth, has real dimension $2m+2n-2$,
    and $\pi_v$ is a submersion with fibres $\cp^{n-1}$.
    Moreover, $\pi_v$ is proper, since the inverse image of any compact
    $K\subset V$ is a closed subset of $K\times A$.

    Let $j_A$ be the complex structure on $A$. The tangent
    spaces of $\mathcal Z_v$ are invariant under $J\oplus j_A$.
    Indeed, the map $(x,H)\mapsto(v(x),H)$ is
    pseudoholomorphic, and the tangent spaces of the incidence
    manifold $\mathcal Z$ are complex subspaces.
    Thus $\mathcal Z_v$ has an induced almost complex structure,
    and both $\pi_v$ and $\sigma_v$ are pseudoholomorphic.
    We use the corresponding orientation on $\mathcal Z_v$.

    Let $A_v$ be the set of regular values of $\sigma_v$.
    Sard's theorem shows that its complement has
    $\lambda$-measure zero. To see that $A_v$ is Borel, observe that
    the critical locus of $\sigma_v$
    is closed in the second-countable manifold $\mathcal Z_v$,
    and hence is a countable union of compact sets. Its image
    under $\sigma_v$ is therefore a countable union of compact
    sets in $A$. Thus its complement $A_v$ is Borel.

    For $H\in A_v$, the projection $\pi_v$ identifies
    $\sigma_v^{-1}(H)$ with $F_H^v$.
    Equivalently, regularity of $H$ means that $v$ is transverse
    to $H$. At every $x\in F_H^v$,
    \begin{align*}
        T_xF_H^v
        =\{X\in T_xV:dv_x(X)\in T_{v(x)}H\}.
    \end{align*}
    This space is $J$-invariant and has real dimension $2m-2$.
    Since $H$ is closed, $F_H^v$ is closed in $V$.
    The regular value theorem gives smoothness and the absence
    of boundary. The fibre orientation induced by $\sigma_v$
    agrees with its complex orientation: the tangent and
    quotient spaces in the submersion exact sequence are
    complex vector spaces, and their real dimensions are even.
    The same observation applies to the fibres of $\pi_v$.

    Fibre integration commutes with this pullback construction:
    \begin{align}
        v^*\beta_f=(\pi_v)_*\sigma_v^*(f\,d\lambda).
        \label{Equation: Base change for incidence forms}
    \end{align}
    To verify the identity, trivialise $\pi$ locally and pull
    that trivialisation back by $v$. The fibre variables and
    their orientations are unchanged. Pullback in the base
    variables consequently commutes with integration in the
    fibre variables. A partition of unity proves
    \eqref{Equation: Base change for incidence forms} globally.

    We can now integrate against $\zeta$ to obtain
    \begin{align*}
        \int_V v^*\beta_f\wedge\zeta=\int_{\mathcal Z_v}
        \sigma_v^*(f\,d\lambda)\wedge\pi_v^*\zeta=\int_A f(H)\left(\int_{F_H^v}\zeta\right)d\lambda(H).
    \end{align*}
    To justify the second equality, observe that at a critical point of $\sigma_v$
    the form $\sigma_v^*(d\lambda)$ vanishes, since $d\lambda$
    has degree $2n$. On the open set where $\sigma_v$ is a
    submersion, use submersion coordinates and a countable
    partition of unity. Fubini's theorem in these coordinates
    gives the displayed equality; the critical values can be
    discarded by Sard's theorem. Moreover, $\pi_v^*\zeta$ has
    compact support by properness of $\pi_v$. The absolute
    value of the smooth top-degree form
    $\sigma_v^*(f\,d\lambda)\wedge\pi_v^*\zeta$ is therefore
    integrable. Applying the same coordinate argument to its
    absolute value proves absolute convergence of the fibre
    integral. The orientation agreements proved above ensure
    that there is no additional sign.

    To prove measurability, cover the
    submersion locus of $\sigma_v$ by countably many submersion
    charts and choose a nonnegative smooth partition of unity
    there. In each such chart, $H$ is a base variable and the
    remaining variables parametrise a portion of $F_H^v$.
    The metric induced from $g$ on that portion has a smooth
    positive volume density. Consequently, for every Borel
    set $B\subset V$, the integral of the corresponding
    partition function over $F_H^v\cap B$ is a parameter
    integral of a nonnegative Borel function.
    It is Borel measurable by Tonelli's theorem.
    Summing over the charts and restricting to the Borel set
    $A_v$ proves measurability of $H\mapsto\nu_H^v(B)$.
    The same local expressions, with the coefficients of
    $\zeta$ evaluated on the oriented tangent vectors, prove
    measurability of $H\mapsto\llbracket F_H^v\rrbracket(\zeta)$.
\end{proof}

\subsection{Extension as integral cycles}
\label{Subsection: Extension of inverse hyperplanes}

We now apply the preceding construction to the smooth representative
of $u$ on $\Reg\,u$. The following proposition verifies all the
hypotheses on the family of currents in
Theorem~\ref{Theorem: Averaged cycles}.

\begin{proposition}
    \label{Proposition: Hyperplane representation}
    Let $u$ satisfy \eqref{Equation: Assumptions on the map}.
    There is a measurable family
    \begin{align*}
        (T_H)_{H\in A}
        \subset\mathbf I_{2m-2,\mathrm{loc}}(M)
    \end{align*}
    of cycles semicalibrated by $\Phi$ such that, for
    $\lambda$-a.e. $H$, $T_H$ is integration with
    multiplicity one and complex orientation over
    \begin{align*}
        F_H:=\{x\in\Reg\,u:u(x)\in H\},
    \end{align*}
    regarded as a current on $M$.
    In particular, $\|T_H\|(\Sing u)=0$ for every $H$.
    For every $f\in C^\infty(A;\mathbb R)$ and every smooth
    compactly supported $(2m-2)$-form $\zeta$ on $M$,
    \begin{align}
        \int_A f(H)T_H(\zeta)\,d\lambda(H)
        =\int_M u^*\beta_f\wedge\zeta.
        \label{Equation: Weighted hyperplane representation}
    \end{align}
    As Radon measures on $M$,
    \begin{align}
        \mu_u=\int_A\|T_H\|\,d\lambda(H).
        \label{Equation: Energy represented by hyperplane cycles}
    \end{align}
    Finally,
    \begin{align}
        \|T_H\|\perp\|T_{H'}\|
        \qquad\text{ for }\,(\lambda\otimes\lambda)\text{-a.e. }
        (H,H')\in A\times A.
        \label{Equation: Mutual singularity of hyperplane cycles}
    \end{align}
\end{proposition}
\begin{proof}
    Write $S:=\Sing u$. Proposition~\ref{Proposition: Initial regularity}
    shows that $S$ is closed and $\Hh^{2m-2}(S)=0$.
    In particular, $S$ has zero Riemannian volume and
    $\mu_u(S)=0$. The smooth representatives of $u$ on its
    regular neighbourhoods agree on overlaps, since they agree
    almost everywhere. They therefore define a single smooth
    pseudoholomorphic map on $\Reg\,u$.

    \medskip
    \noindent\emph{Step 1: mass representation on $\Reg\,u$.}
    Apply Lemma~\ref{Lemma: Smooth inverse hyperplanes} to this
    map, and denote the resulting Borel set of regular
    parameters by $A_0$.
    For $H\in A_0$, let
    $\nu_H:=\Hh^{2m-2}\!\llcorner F_H$ on $\Reg\,u$; set $\nu_H=0$ otherwise.
    These measures form a Borel kernel.
    By Lemma~\ref{Lemma: Energy identity}, the complex
    orientation of $F_H$ satisfies
    \begin{align*}
        \llbracket F_H\rrbracket(\chi\Phi)=\int_{F_H}\chi\,d\Hh^{2m-2}
    \end{align*}
    for every $\chi\in C_c^\infty(\Reg\,u)$.
    Taking $f=1$ and $\zeta=\chi\Phi$ in
    \eqref{Equation: Weighted identity for smooth maps}, for
    $\chi\geq0$, gives
    \begin{align}
        \int_A\left(\int_{\Reg\,u}\chi\,d\nu_H\right)d\lambda(H)=\int_{\Reg\,u}\chi\,u^*\omega\wedge\Phi=\int_{\Reg\,u}\chi\,d\mu_u.
        \label{Equation: Mass representation on the regular set}
    \end{align}
    Tonelli's theorem shows that $B\mapsto\int_A\nu_H(B)d\lambda(H)$
    is a measure on $\Reg\,u$. It is locally finite, by testing
    \eqref{Equation: Mass representation on the regular set}
    with a smooth function equal to one on any prescribed
    compact subset of $\Reg\,u$.
    Equality on the smooth nonnegative test functions thus
    gives equality of Radon measures. Consequently,
    \begin{align}
        \int_A\nu_H(B)\,d\lambda(H)=\mu_u(B)
        \qquad\text{for every Borel set } B\subset\Reg\,u.
        \label{Equation: Borel mass representation on the regular set}
    \end{align}
    In particular, this identity also holds for Borel sets
    whose closures meet $S$.

    \medskip
    \noindent\emph{Step 2: extension with locally finite mass.}
    Extend each $\nu_H$ to $M$ by setting
    $\nu_H(B):=\nu_H(B\cap\Reg\,u)$ for Borel $B\subset M$.
    This preserves the Borel kernel property.
    Choose an increasing exhaustion $W_j\Subset M$ by
    relatively compact open sets, with
    $\overline{W_j}\subset W_{j+1}$.
    By \eqref{Equation: Borel mass representation on the regular set},
    \begin{align*}
        \int_A\nu_H(W_j)\,d\lambda(H)
        =\mu_u(W_j\cap\Reg\,u)=\mu_u(W_j)<\infty.
    \end{align*}
    Therefore
    \begin{align*}
        A_1:=A_0\cap\bigcap_{j=1}^{\infty}
        \{H:\nu_H(W_j)<\infty\}
    \end{align*}
    is Borel and has full $\lambda$-measure.
    For $H\in A_1$, define
    \begin{align*}
        T_H(\zeta):=\int_{F_H}\zeta
    \end{align*}
    for every smooth compactly supported $(2m-2)$-form $\zeta$ on $M$,
    and set $T_H=0$ for $H\notin A_1$.
    The local mass bound makes this integral absolutely
    convergent. The submanifold $F_H$ is countably rectifiable
    in $M$ and has multiplicity one, so $T_H$ is a locally
    integer-rectifiable current, with
    \begin{align*}
        \|T_H\|=\nu_H
        \qquad\forall\,H\in A_1.
    \end{align*}
    Since $F_H$ carries its complex orientation, $T_H$ is semicalibrated
    by $\Phi$. Moreover, $\nu_H(S)=0$ by construction.
    After setting the measures to zero on $A\smallsetminus A_1$,
    \eqref{Equation: Borel mass representation on the regular set}
    and $\mu_u(S)=0$ give
    \eqref{Equation: Energy represented by hyperplane cycles}
    on all Borel subsets of $M$.

    \medskip
    \noindent\emph{Step 3: the weighted identity on $M$.}
    Choose functions $\chi_j\in C_c^\infty(\Reg\,u)$ such that
    $0\leq\chi_j\leq1$ and $\chi_j\to1$ pointwise on $\Reg\,u$.
    Fix a smooth compactly supported $(2m-2)$-form $\zeta$ on
    $M$ and a smooth real function $f$ on $A$.
    Since $\chi_j\zeta$ has compact support in $\Reg\,u$, the smooth identity gives
    \begin{align*}
        \int_A f(H)T_H(\chi_j\zeta)\,d\lambda(H)
        =\int_{\Reg\,u}\chi_j\,u^*\beta_f\wedge\zeta.
    \end{align*}
    For each $H\in A_1$, dominated convergence on $F_H$ gives
    $T_H(\chi_j\zeta)\rightarrow T_H(\zeta)$.
    Moreover,
    \begin{align*}
        |T_H(\chi_j\zeta)|
        \leq\sup_M|\zeta|_{\mathrm{comass}}
        \|T_H\|(\spt\zeta),
    \end{align*}
    and the right-hand side is integrable in $H$ by
    \eqref{Equation: Energy represented by hyperplane cycles}.
    Each evaluation on $\chi_j\zeta$ is Borel measurable by
    Lemma~\ref{Lemma: Smooth inverse hyperplanes} and the
    definition of $A_1$. Passing to the pointwise limit also
    proves measurability of $H\mapsto T_H(\zeta)$.
    Together with the Borel mass kernel, this is precisely
    the measurability required in
    Section~\ref{Section: Averaged semicalibrated cycles}.

    Since $\beta_f$ is smooth on the compact target,
    \begin{align*}
        |u^*\beta_f|\leq C_f|du|^2
        \qquad\text{a.e. on }M.
    \end{align*}
    The right-hand side is locally integrable on $M$, and $S$ has
    zero volume. Dominated convergence on both sides therefore
    proves \eqref{Equation: Weighted hyperplane representation}.
    This passage to the limit uses only pointwise convergence and the
    bounds above, so it requires no estimate on $d\chi_j$.

    \medskip
    \noindent\emph{Step 4: vanishing of the individual boundaries.}
    Since $d(u^*\omega)=0$, Lemmas~\ref{Lemma: Weighted incidence forms}
    and~\ref{Lemma: Sobolev pullbacks} give
    \begin{align*}
        d(u^*\beta_f)
        =c_f\,d(u^*\omega)+d^2(u^*\eta_f)=0
        \qquad\text{in distributions on }M
    \end{align*}
    for every smooth real $f$ on $A$.
    Let $\xi$ be a smooth compactly supported $(2m-3)$-form
    on $M$. We may now apply
    \eqref{Equation: Weighted hyperplane representation} to
    the fixed test form $d\xi$, obtaining
    \begin{align}
        \int_A f(H)\,\partial T_H(\xi)\,d\lambda(H)
        =\int_Mu^*\beta_f\wedge d\xi=0
        \qquad\text{for every }f\in C^\infty(A;\mathbb R).
        \label{Equation: Weighted vanishing of the boundaries}
    \end{align}
    The function $H\mapsto\partial T_H(\xi)=T_H(d\xi)$
    is Borel measurable and integrable, since
    \begin{align*}
        |\partial T_H(\xi)|
        \leq\sup_M|d\xi|_{\mathrm{comass}}
        \|T_H\|(\spt\xi).
    \end{align*}
    Smooth functions are uniformly dense in $C^0(A)$.
    Equation~\eqref{Equation: Weighted vanishing of the boundaries}
    therefore says that the finite signed measure
    $\partial T_H(\xi)\,d\lambda(H)$ is zero.
    It follows that $\partial T_H(\xi)=0$ for a.e. $H$.

    We must make the exceptional parameter set independent of
    $\xi$. To this end, for each $j$ choose a countable family of smooth
    $(2m-3)$-forms supported in $\overline{W_j}$ which is dense
    in the $C^1$ norm among such forms. These spaces are
    separable in that norm, as can be seen in a finite
    collection of coordinate charts covering $\overline{W_j}$.
    Apply the preceding conclusion to the union of these
    countable families and discard the resulting Borel null
    set of parameters.
    For every remaining $H$, every test form $\xi$ is supported
    in some $W_j$ and can be approximated in $C^1$ by the
    corresponding dense family. The estimate
    \begin{align*}
        |T_H(d\xi-d\xi_\ell)|
        \leq\sup_M|d\xi-d\xi_\ell|_{\mathrm{comass}}
        \|T_H\|(\overline{W_j})
        \rightarrow0
    \end{align*}
    proves $\partial T_H(\xi)=0$ for every test form.
    Thus $\partial T_H=0$ simultaneously for all these parameters.
    Redefine $T_H=0$ on the discarded Borel null set.
    All preceding identities and measurability properties
    remain valid. The currents are now locally
    integer-rectifiable with zero boundary, and hence belong to
    $\mathbf I_{2m-2,\mathrm{loc}}(M)$.

    \medskip
    \noindent\emph{Step 5: mutual singularity.}
    For this step, use the set-theoretic definition
    $F_{H'}=\{x\in\Reg\,u:u(x)\in H'\}$ for every $H'\in A$,
    including those for which we have set $T_{H'}=0$. The incidence set
    \begin{align*}
        \{(H',x)\in A\times\Reg\,u \mbox{ : } u(x)\in H'\}
    \end{align*}
    is Borel. Together with the Borel mass kernel, this shows
    that $(H,H')\mapsto\|T_H\|(F_{H'}\cap W_j)$ is measurable.
    For each fixed $H$, Tonelli's theorem gives
    \begin{align*}
        \int_A\|T_H\|(F_{H'}\cap W_j)\,d\lambda(H')=\int_{\Reg\,u\cap W_j}
        \lambda\big(\{H'\in A\mbox{ : }u(x)\in H'\}\big)\,d\|T_H\|(x)=0.
    \end{align*}
    Indeed, the hyperplanes containing a fixed point of
    $\cp^n$ form a copy of $\cp^{n-1}$ in $A=(\cp^n)^*$,
    and consequently have zero measure for the smooth
    $2n$-dimensional volume $\lambda$.
    For $n=1$, this set consists of one point and is still null.
    Applying Tonelli once more and then taking the countable
    union of the exceptional sets associated with the $W_j$
    gives
    \begin{align*}
        \|T_H\|(F_{H'})=0
        \qquad\text{ for }\,(\lambda\otimes\lambda)\text{-a.e. }
        (H,H').
    \end{align*}
    On the other hand, $\|T_{H'}\|$ is concentrated on $F_{H'}$
    for every $H'$, including those for which $T_{H'}=0$.
    Since $F_{H'}$ is Borel, this proves the mutual singularity in
    \eqref{Equation: Mutual singularity of hyperplane cycles}.
\end{proof}

In particular, for every relatively compact open set $W\Subset M$,
\begin{align*}
    \int_A\|T_H\|(W)\,d\lambda(H)=\mu_u(W)<\infty.
\end{align*}
This local mass bound, the common semicalibration $\Phi$, and mutual
singularity of the mass measures for almost every pair verify the
hypotheses of Theorem~\ref{Theorem: Averaged cycles} for the family
in Proposition~\ref{Proposition: Hyperplane representation}, with
$V=M$, $q=2m-2$, and averaged measure $\mu_u$.

\begin{remark}
    \label{Remark: Supports of the extended hyperplane cycles}
    The assertion $\|T_H\|(\Sing u)=0$ does not imply that
    $\spt T_H$ avoids $\Sing u$. The closure in $M$ of $F_H$
    may meet the singular set of the map. The extension argument
    proves that this accumulation creates no boundary for almost
    every hyperplane. It is the intersections of the supports
    of these extended cycles that enter the averaging theorem.
\end{remark}
\section{Partial regularity and the structure of the singular set}
\label{Section: Partial regularity and sharpness}

We combine the averaging theorem with the representation by
inverse-hyperplane cycles to obtain the regularity statement for $u$.
The rectifiability theorem for semicalibrated currents, together
with a measure-theoretic result of De Lellis--Fleschler, also gives
countable rectifiability of the singular set. We then recover local finiteness in real dimension four using \cite{RiviereTianCurrents} and exhibit a map whose singular
set has Hausdorff dimension exactly $2m-4$. Finally, we give a sufficient
condition for local finiteness of the critical measure in higher dimensions
in terms of the same family of currents.

\subsection{The partial regularity theorem}
\label{Subsection: Proof of partial regularity}

We can now give the proof of Theorem \ref{Theorem: Intro main regularity}. Recall that a set is \emph{countably $(2m-4)$-rectifiable} if, in local coordinates, it is covered up to a set of zero $\Hh^{2m-4}$-measure by countably many Lipschitz images of subsets of $\mathbb R^{2m-4}$. For $m=2$, this means that the set is countable. This convention does not include local finiteness of its Hausdorff measure.

\begin{proof}[Proof of Theorem \ref{Theorem: Intro main regularity}]
    Fix a smooth $J$-compatible metric $g$, and use the notation
    of Section~\ref{Section: Preliminaries}.

    \medskip
    \noindent\emph{Step 1: the dimension estimate and regularity.}
    Proposition~\ref{Proposition: Hyperplane representation}
    provides a measurable family $(T_H)_{H\in A}$ of locally
    integral $(2m-2)$-cycles, semicalibrated by the common
    smooth form $\Phi$, such that
    \begin{align*}
        \mu_u=\int_A\|T_H\|\,d\lambda(H).
    \end{align*}
    Their local masses are integrable in $H$, and their mass
    measures are mutually singular for almost every pair of
    parameters. Theorem~\ref{Theorem: Averaged cycles}, with
    $V=M$ and $q=2m-2$, therefore gives
    \begin{align*}
        \dim_{\Hh}\{x\in M\mbox{ : }\theta_{\mu_u}(x)>0\}\leq2m-4.
    \end{align*}
    The density used in that theorem agrees with the one in
    Proposition~\ref{Proposition: Initial regularity}, since
    \begin{align*}
        \theta_{\mu_u}(x)
        =\lim_{r\rightarrow0}r^{2-2m}\mu_u(B_r(x))
        =\lim_{r\rightarrow0}\mathscr{E}_u(x,r)
        =\theta_u(x).
    \end{align*}
    Hence \eqref{Equation: Initial singular set estimate} identifies the positive-density set of $\mu_u$ with $\Sing u$, giving $\dim_{\Hh}\Sing u\leq2m-4$.
    Closedness of $\Sing u$ and smoothness on its complement
    follow from Proposition~\ref{Proposition: Initial regularity}.
    The smooth local representatives agree on overlaps because
    they agree almost everywhere. They consequently define a
    smooth representative on $\Reg\,u$, and its
    pseudoholomorphic equation holds everywhere by continuity.

    \medskip
    \noindent\emph{Step 2: rectifiability of the pairwise intersections.}
    We first treat $m>2$. For $H,H'\in A$, set
    \begin{align}
        C_{H,H'}:=\spt T_H\cap\spt T_{H'}.
        \label{Equation: Pairwise intersections for rectifiability}
    \end{align}
    For a.e. pair, mutual singularity and
    Lemma~\ref{Lemma: Pairwise intersection} give
    \begin{align*}
        C_{H,H'}\subset\Sing(T_H+T_{H'}).
    \end{align*}
    The sum is a locally integral $(2m-2)$-cycle semicalibrated
    by the same smooth form $\Phi$. The rectifiability theorem
    of Minter--Parise--Skorobogatova--Spolaor
    \cite[Theorem~1.2]{MinterPariseSkorobogatovaSpolaor}
    therefore implies that $\Sing(T_H+T_{H'})$, and hence its
    Borel subset $C_{H,H'}$, is countably $(2m-4)$-rectifiable.
    The theorem applies locally in a smooth Riemannian ambient
    manifold by the same isometric embedding and extension of the
    semicalibration as in
    Section~\ref{Section: Preliminaries}. Localisation gives the
    conclusion for locally integral currents, since the added
    boundary of a restriction lies outside the region under
    consideration.

    Retain the incidence probability from
    \eqref{Equation: Incidence probability}, namely
    \begin{align*}
        \mathfrak p(x)
        =\lambda\big(\{H\in A\mbox{ : }x\in\spt T_H\}\big).
    \end{align*}
    Lemma~\ref{Lemma: Density and incidence} and the identity
    $\Sing u=\{\theta_u>0\}$ give
    \begin{align*}
        \mathfrak p(x)>0
        \qquad\text{for every }x\in\Sing u.
    \end{align*}

    \medskip
    \noindent\emph{Step 3: exclusion of purely unrectifiable subsets.}
    Let $F\subset\Sing u$ be a Borel, purely
    $(2m-4)$-unrectifiable set with finite $\Hh^{2m-4}$-measure.
    Pure unrectifiability means that its intersection with every
    countably $(2m-4)$-rectifiable set has zero
    $\Hh^{2m-4}$-measure. By Step~2,
    \begin{align*}
        \Hh^{2m-4}(F\cap C_{H,H'})=0
        \qquad\text{for a.e. }(H,H').
    \end{align*}
    The incidence set is jointly measurable by
    \eqref{Equation: Measurability of support incidence}.
    Since $\Hh^{2m-4}\!\llcorner F$ is a finite measure,
    Tonelli's theorem yields
    \begin{align}
        0
        &=\int_{A\times A}
        \Hh^{2m-4}(F\cap C_{H,H'})
        \,d(\lambda\otimes\lambda)(H,H')\notag\\
        &=\int_F\left(\int_A
        \mathbf{1}_{\spt T_H}(x)\,d\lambda(H)\right)^2
        \,d\Hh^{2m-4}(x)\notag\\
        &=\int_F\mathfrak p(x)^2\,d\Hh^{2m-4}(x).
        \label{Equation: Exclusion of purely unrectifiable subsets}
    \end{align}
    The integrand is strictly positive on $F$. Consequently,
    $\Hh^{2m-4}(F)=0$.
    Since the measure used in Tonelli's theorem is restricted to $F$,
    this argument applies before $\sigma$-finiteness on $\Sing u$ has been established.

    \medskip
    \noindent\emph{Step 4: $\sigma$-finiteness and rectifiability.}
    We use the following consequence of the theorem of
    De Lellis--Fleschler \cite[Theorem~1.1]{DeLellisFleschler}: a Borel subset
    of $\mathbb R^{2m}$ which is not $\sigma$-finite with respect
    to $\Hh^{2m-4}$ contains a closed, purely
    $(2m-4)$-unrectifiable subset of positive finite
    $\Hh^{2m-4}$-measure. The restriction $m>2$ ensures that
    the exponent $2m-4$ in this statement is a positive
    integer strictly smaller than the ambient dimension.

    Cover $M$ by countably many compact sets contained in
    coordinate neighbourhoods on which the coordinates are
    bi-Lipschitz. If $\Hh^{2m-4}\!\llcorner\Sing u$ were not $\sigma$-finite,
    there would be such a compact set $K$ for which $\Sing u\cap K$
    is not $\sigma$-finite with respect to $\Hh^{2m-4}$. Apply the
    De Lellis--Fleschler theorem to the coordinate image of $\Sing u\cap K$.
    Pulling the resulting subset back to $M$ would give a Borel,
    purely $(2m-4)$-unrectifiable subset of $\Sing u$ with
    positive finite $\Hh^{2m-4}$-measure, contrary to Step~3.
    Here we have used that bi-Lipschitz changes of coordinates
    preserve pure unrectifiability, finiteness and positivity
    of Hausdorff measure, and $\sigma$-finiteness. Thus
    $\Hh^{2m-4}\!\llcorner\Sing u$ is $\sigma$-finite.

    We can therefore write $\Sing u$ as a countable union of
    Borel sets of finite $\Hh^{2m-4}$-measure, each contained
    in a coordinate neighbourhood. Apply the usual
    rectifiable--purely unrectifiable decomposition to each
    of these sets; see \cite{Mattila}. The purely unrectifiable
    part has zero measure by Step~3. Each piece is consequently
    countably $(2m-4)$-rectifiable, and so is their countable
    union.

    Finally, let $m=2$ and fix a compact set $K\subset M$. Proposition~\ref{Proposition: Initial regularity} supplies $\varepsilon_K>0$ such that
    \begin{align*}
    	\theta_u(x)\geq\varepsilon_K
    	\qquad\text{for every }x\in\Sing u\cap K.
    \end{align*}
    As $m=2$, the cycles in Proposition~\ref{Proposition: Hyperplane representation} are positive $(1,1)$ integral cycles. By Rivi\`ere--Tian
    \cite{RiviereTianCurrents}, the singular set of each pairwise
    sum is locally finite. The inclusion
    \eqref{Equation: Intersections are singular for the sum} thus
    gives the pairwise finiteness used in the $q=2$ case of
    Theorem~\ref{Theorem: Averaged cycles}. We conclude that
    \begin{align*}
    	\{x\in K:\theta_u(x)\geq\varepsilon_K\}
    \end{align*}
    is finite. Since this set contains $\Sing u\cap K$, the singular set is locally finite. A countable compact exhaustion of $M$ then shows that $\Sing u$ is countable and hence countably $0$-rectifiable. If $M$ is compact, taking $K=M$ gives finiteness.
\end{proof}
\subsection{Optimality of the dimension exponent}
\label{Subsection: Sharpness of codimension four}

The exponent $2m-4$ in \eqref{Equation: Intro Haussdorf dimension bound} cannot be replaced by a smaller one, even when the complex structure of the domain
is integrable and the image is contained in a projective line.

\begin{proposition}
    \label{Proposition: Sharpness of codimension four}
    Let $m\geq2$ and $n\geq1$, and equip $\mathbb C^m$ with its
    standard complex structure and Euclidean metric.
    Set
    \begin{align*}
        L:=\{z\in\mathbb C^m\mbox{ : }z_1=z_2=0\}.
    \end{align*}
    The map
    \begin{align}
        u(z_1,\ldots,z_m):=[z_1:z_2:0:\cdots:0]
        \qquad\text{on }\mathbb C^m\smallsetminus L,
        \label{Equation: Sharpness example}
    \end{align}
    assigned any fixed value in $\cp^n$ on $L$, belongs to
    $W^{1,2}_{\mathrm{loc}}(\mathbb C^m,\cp^n)$ and satisfies
    \eqref{Equation: Assumptions on the map}. Moreover,
    \begin{align*}
        \Sing u=L,
        \qquad
        \dim_{\Hh}\Sing u=2m-4.
    \end{align*}
    In fact,
    \begin{align}
        0<\Hh^{2m-4}(\Sing u\cap B_r(p))<\infty
        \qquad(p\in L,\ r>0).
        \label{Equation: Critical measure of the dimension example}
    \end{align}
\end{proposition}
\begin{proof}
    The expression in \eqref{Equation: Sharpness example} is
    well defined and holomorphic away from $L$.
    Write
    \begin{align*}
        \rho(z):=(|z_1|^2+|z_2|^2)^{1/2}.
    \end{align*}
    The map is independent of $z_3,\ldots,z_m$, and its
    dependence on $(z_1,z_2)$ is homogeneous of degree zero
    under positive real dilations. Its restriction to the unit
    sphere in $\mathbb C^2$ is smooth. Scaling the differential
    therefore gives
    \begin{align}
        |du(z)|\leq\frac{C}{\rho(z)}
        \qquad(z\notin L).
        \label{Equation: Differential bound for the sharpness example}
    \end{align}
    Here the target metric is the normalised Fubini--Study
    metric used throughout the paper.

    \medskip
    \noindent\emph{Step 1: the Sobolev extension.}
    Fix a smooth isometric embedding
    $\iota:\cp^n\hookrightarrow\mathbb R^N$, and put
    $U:=\iota\circ u$. The map $U$ is bounded because the
    target is compact. Its classical differential off $L$
    satisfies the bound in
    \eqref{Equation: Differential bound for the sharpness example}.
    For every compact $K\subset\mathbb C^m$, integration in
    the four real variables corresponding to $(z_1,z_2)$ gives,
    for $0<\varepsilon<1$,
    \begin{align}
        \int_{K\cap\{0<\rho<2\varepsilon\}}|dU|^2\,d\vol
        \leq C_K\int_0^{2\varepsilon}t^{3}t^{-2}\,dt
        \leq C_K\varepsilon^2.
        \label{Equation: Transverse energy integrability}
    \end{align}
    We use $d\vol$ for Euclidean volume in this proof.
    The constant absorbs the volume of a bounded set containing
    the projection of $K$ onto the remaining variables.
    When $m=2$, there are no remaining variables, and the same
    estimate holds. Since $U$ is smooth away from $L$,
    \eqref{Equation: Transverse energy integrability} shows that
    its classical differential, extended by zero on $L$, is
    locally square integrable.

    We verify that this is the distributional differential.
    Choose $\chi\in C^\infty([0,\infty))$ such that
    $0\leq\chi\leq1$, $\chi=0$ on $[0,1]$, and $\chi=1$ on
    $[2,\infty)$, and set
    \begin{align*}
        \chi_\varepsilon(z):=\chi(\rho(z)/\varepsilon).
    \end{align*}
    This function is smooth on $\mathbb C^m$, is zero near $L$,
    and satisfies
    \begin{align*}
        |d\chi_\varepsilon|
        \leq\frac{C}{\varepsilon}
        \mathbf{1}_{\{\varepsilon<\rho<2\varepsilon\}}.
    \end{align*}
    For a compactly supported smooth real test function $\varphi$
    and a real coordinate direction $i\in\{1,\ldots,2m\}$,
    integration by parts componentwise gives
    \begin{align*}
        -\int_{\mathbb C^m}\chi_\varepsilon U\,\partial_i\varphi\,d\vol
        &=\int_{\mathbb C^m}
        \chi_\varepsilon\partial_iU\,\varphi\,d\vol
        +\int_{\mathbb C^m}
        (\partial_i\chi_\varepsilon)U\varphi\,d\vol.
    \end{align*}
    The volume of the transverse shell
    $\{\varepsilon<\rho<2\varepsilon\}$ over a fixed bounded
    set of tangential variables is at most $C\varepsilon^4$.
    Boundedness of $U$ therefore gives
    \begin{align*}
        \left|\int_{\mathbb C^m}
        (\partial_i\chi_\varepsilon)U\varphi\,d\vol\right|
        \leq C_\varphi\varepsilon^{-1}\varepsilon^4
        =C_\varphi\varepsilon^3\rightarrow0.
    \end{align*}
    The other terms converge by dominated convergence and
    local integrability of the classical differential.
    Hence $\partial_iU$ off $L$, extended by zero, is the weak
    derivative of $U$. This proves
    $u\in W^{1,2}_{\mathrm{loc}}(\mathbb C^m,\cp^n)$.
    Since $L$ has zero volume and $u$ is holomorphic on its
    complement, $u$ is weakly pseudoholomorphic.

    \medskip
    \noindent\emph{Step 2: distributional closedness of the pullback.}
    Let $\alpha:=u^*\omega$. This is a locally integrable
    2-form, smooth and closed on $\mathbb C^m\smallsetminus L$.
    The pointwise estimate for pullbacks and
    \eqref{Equation: Differential bound for the sharpness example}
    give
    \begin{align*}
        |\alpha(z)|\leq C|du(z)|^2
        \leq\frac{C}{\rho(z)^2}
        \qquad(z\notin L).
    \end{align*}
    Fix a smooth compactly supported $(2m-3)$-form $\xi$ on
    $\mathbb C^m$. Since $\chi_\varepsilon\xi$ is supported
    away from $L$, we can use $d\alpha=0$ on its support to obtain
    \begin{align*}
        \int_{\mathbb C^m}\chi_\varepsilon\alpha\wedge d\xi
        =-\int_{\mathbb C^m}
        d\chi_\varepsilon\wedge\alpha\wedge\xi.
    \end{align*}
    On the support of $d\chi_\varepsilon$ we have
    $|\alpha|\leq C\varepsilon^{-2}$.
    Using again the transverse shell volume, we obtain
    \begin{align*}
        \left|\int_{\mathbb C^m}
        d\chi_\varepsilon\wedge\alpha\wedge\xi\right|
        \leq C_\xi\varepsilon^{-1}\varepsilon^{-2}\varepsilon^4
        =C_\xi\varepsilon\rightarrow0.
    \end{align*}
    Since $\alpha$ is locally integrable, dominated convergence
    on the left yields
    \begin{align*}
        \int_{\mathbb C^m}\alpha\wedge d\xi=0.
    \end{align*}
    Thus $d(u^*\omega)=0$ in distributions on all of
    $\mathbb C^m$.

    \medskip
    \noindent\emph{Step 3: identification of the singular set.}
    Smoothness away from $L$ gives $\Sing u\subset L$.
    Fix $p=(0,0,p_3,\ldots,p_m)\in L$ and suppose that the
    Sobolev class of $u$ has a continuous representative on
    a neighbourhood of $p$.
    On the complement of $L$ in that neighbourhood, this
    representative and the smooth map in
    \eqref{Equation: Sharpness example} agree almost everywhere.
    Since both are continuous there, they agree everywhere.
    However, for positive $t$ tending to zero,
    \begin{align*}
        u(t,0,p_3,\ldots,p_m)&=[1:0:0:\cdots:0],\\
        u(0,t,p_3,\ldots,p_m)&=[0:1:0:\cdots:0].
    \end{align*}
    These are distinct points of $\cp^n$, contradicting
    continuity at $p$. Thus every point of $L$ is singular,
    independently of the value assigned to $u$ on $L$.
    Finally, $L$ is a complex linear subspace of dimension
    $m-2$, so its Hausdorff dimension is $2m-4$.
    For $m>2$, its intersection with a ball centred on $L$
    is a Euclidean ball of dimension $2m-4$, which has positive
    finite $\Hh^{2m-4}$-measure. For $m=2$, $L=\{0\}$ and
    the corresponding $\Hh^0$-measure equals one. This proves
    \eqref{Equation: Critical measure of the dimension example}.
\end{proof}

Since the singular plane has positive $\Hh^{2m-4}$-measure,
the example also rules out a general conclusion of the form
$\Hh^{2m-4}(\Sing u)=0$.

\subsection{The critical Hausdorff measure}
\label{Subsection: Scope of partial regularity}

For $m>2$, Theorem~\ref{Theorem: Intro main regularity} gives countable
$(2m-4)$-rectifiability and $\sigma$-finiteness of the critical
Hausdorff measure. The argument does not establish its local
finiteness. The rectifiability theorem for semicalibrated currents
used above does not supply this additional conclusion, and
Proposition~\ref{Proposition: Sharpness of codimension four}
provides no counterexample to it.

The following conditional statement identifies a sufficient
additional input. In particular, one would not need an integrable
bound for the critical measures of the pairwise intersections.

\begin{proposition}
    \label{Proposition: Conditional local finiteness}
    Under the assumptions of Theorem~\ref{Theorem: Intro main regularity},
    let $m>2$ and let $(T_H)_{H\in A}$ be the family from
    Proposition~\ref{Proposition: Hyperplane representation}.
    Fix a compact set $K\subset M$. Assume, in addition, that
    \begin{align}
        \Hh^{2m-4}(\spt T_H\cap\spt T_{H'}\cap K)<\infty
        \qquad\text{ for }\,(\lambda\otimes\lambda)
        \text{-a.e. }(H,H').
        \label{Equation: Additional pairwise finiteness assumption}
    \end{align}
    Then
    \begin{align*}
        \Hh^{2m-4}(\Sing u\cap K)<\infty.
    \end{align*}
\end{proposition}
\begin{proof}
    Proposition~\ref{Proposition: Initial regularity} gives
    $\varepsilon_K>0$ such that $\theta_u\geq\varepsilon_K$
    on $\Sing u\cap K$. By
    Lemma~\ref{Lemma: Density and incidence}, there is
    $\delta_K>0$ such that
    \begin{align*}
        \mathfrak p(x)\geq\delta_K
        \qquad\text{for every }x\in\Sing u\cap K.
    \end{align*}
    Suppose that $\Hh^{2m-4}(\Sing u\cap K)=\infty$.
    The $\sigma$-finiteness proved in
    Theorem~\ref{Theorem: Intro main regularity} then supplies Borel
    sets $F_j\subset\Sing u\cap K$ satisfying
    \begin{align*}
        0<\Hh^{2m-4}(F_j)<\infty,
        \qquad
        \Hh^{2m-4}(F_j)\rightarrow\infty.
    \end{align*}
    Indeed, take increasing finite unions of a countable cover of
    $\Sing u\cap K$ by Borel sets of finite measure, discarding any
    initial unions of zero measure. Define probability measures
    \begin{align*}
        \nu_j:=
        \frac{\Hh^{2m-4}\!\llcorner F_j}
             {\Hh^{2m-4}(F_j)}.
    \end{align*}
    Recall the notation $C_{H,H'}$ from
    \eqref{Equation: Pairwise intersections for rectifiability}.
    Since $F_j\subset K$, \eqref{Equation: Additional pairwise finiteness assumption}
    gives, for almost every pair $(H,H')$,
    \begin{align*}
        0\leq\nu_j(C_{H,H'})
        \leq
        \frac{\Hh^{2m-4}(C_{H,H'}\cap K)}
             {\Hh^{2m-4}(F_j)}
        \rightarrow0.
    \end{align*}
    Each function $(H,H')\mapsto\nu_j(C_{H,H'})$ is measurable
    by the joint measurability of support incidence.
    It takes values in $[0,1]$, so dominated convergence yields
    \begin{align*}
        \int_{A\times A}\nu_j(C_{H,H'})
        \,d(\lambda\otimes\lambda)(H,H')
        \rightarrow0.
    \end{align*}
    On the other hand, Tonelli's theorem gives
    \begin{align*}
        \int_{A\times A}\nu_j(C_{H,H'})
        \,d(\lambda\otimes\lambda)(H,H')=\int_M\mathfrak p(x)^2\,d\nu_j(x)\geq\delta_K^2,
    \end{align*}
    a contradiction.
\end{proof}
\section{A conservation law on real-analytic almost complex manifolds}
\label{Section: Automatic closedness}

We now prove Theorem~\ref{Theorem: Automatic closedness}. Throughout
this section the domain $(M^{2m},J)$ is a real analytic almost complex manifold and $g$ is a smooth and $J$-compatible metric on $M$. We retain the notation of Section~\ref{Section: Preliminaries} and we assume that $u\in W^{1,2}_{\mathrm{loc}}(M,\cp^n)$ satisfies the pseudoholomorphic equation \eqref{Equation: Intro pseudoholomorphic equation}.

An almost complex submanifold is an embedded
submanifold $S$ whose tangent bundle is $J$-invariant. We equip $S$
with the induced structure $J_S=J|_{TS}$, which is real analytic
whenever $S$ is a real-analytic submanifold and $J$ is real analytic.

Here, $\omega$ denotes the Fubini--Study form with the normalisation fixed in the introduction. Recall that, for $m\geq2$, distributional closedness means that
\begin{align}\label{Equation: Distributional closedness}
    \int_M u^*\omega\wedge d\eta=0
    \qquad\forall\,\eta\in C_c^\infty(M;\Lambda^{2m-3}T^*M).
\end{align}
For $m=1$, closedness is automatic since $u^*\omega$ has top degree.
We first prove closedness for maps into $\cp^1$ by induction on the
domain dimension. When $J$ is integrable, we use the holomorphic-chain
theorems of King and Harvey--Shiffman, in the codimension-one form
proved by Shiffman (Subsection~\ref{Subsection: Integrable sphere valued case}).
To pass from $\cp^1$ to $\cp^n$, we justify composition of an arbitrary
$W^{1,2}_{\mathrm{loc}}$ map with almost every singular projective
projection of \cite[Appendix~A]{CaniatoRiviere}. Averaging over these
projections then gives the induction step in the target dimension.

For the induction on the domain dimension, an analytic vector
field $X$ satisfying $du(X)=0$ and $X(p)\ne0$ allows us to factor
$u$ locally through an analytic almost complex domain of dimension
$2m-2$ (Lemma~\ref{Lemma: Analytic dimension reduction}). In the nonintegrable case we choose finitely many
analytic fields annihilated by $du$ and let $Z$ be their common
zero set. Since at least one field is nonzero at each point outside
$Z$, the induction hypothesis gives closedness there. The strata of
codimension two are almost complex, and each stratum of codimension
three is locally contained in an almost complex hypersurface
(Lemma~\ref{Lemma: Geometry of the common zero set}). Positivity
therefore gives closedness across both kinds of strata. An elementary
support argument then gives closedness across the remaining set of
codimension at least four
(Subsections~\ref{Subsection: Removability} and~\ref{Subsection: Sphere in every dimension}). This argument requires no mass bound on the distributional derivative of $u^*\omega$.

\subsection{Currents and positivity}
\label{Subsection: Currents and positivity}

If $T$ is a $q$-current of locally finite mass on an open set $V$, its action on a
test form is given by
\begin{align*}
    T(\eta)=\int_V\langle\eta_x,\vec T(x)\rangle\,d\|T\|(x),
\end{align*}
where $\vec T$ has mass norm one $\|T\|$-almost everywhere, the mass
norm being dual to the comass norm. This representation extends $T$
to compactly supported bounded Borel forms and defines its restriction
$T\llcorner E:=\mathbf{1}_E T$ to any Borel set $E$. We call $T$
locally normal if both $T$ and $\partial T$ have locally finite mass.

If $\alpha$ is a smooth $r$-form with $0\leq r\leq q$, we define
the contraction $T\llcorner\alpha$, a current of dimension $q-r$,
by
\begin{align*}
    (T\llcorner\alpha)(\eta):=T(\alpha\wedge\eta).
\end{align*}
This convention gives, for $r\leq q-1$,
\begin{align*}
    \partial(T\llcorner\alpha)
    =(-1)^r\bigl((\partial T)\llcorner\alpha
                 -T\llcorner d\alpha\bigr).
\end{align*}
Indeed, evaluate both sides on a test form $\eta$ and use
$d(\alpha\wedge\eta)=d\alpha\wedge\eta
 +(-1)^r\alpha\wedge d\eta$. In particular,
\begin{align*}
    \partial(fT)=f\partial T-T\llcorner df
\end{align*}
for smooth $f$ and $q\geq1$. Contractions with continuous forms are also defined
by the integral representation whenever $T$ has locally finite mass.

An $L^1_{\mathrm{loc}}$ $k$-form $\alpha$ determines a current
of dimension $2m-k$ by
\begin{align*}
    T_\alpha(\zeta):=\int_M\alpha\wedge\zeta.
\end{align*}
Its mass is locally finite. For smooth compactly supported
$(2m-k-1)$-forms $\eta$, the distributional convention is
\begin{align*}
    \langle d\alpha,\eta\rangle
    =(-1)^{k+1}\int_M\alpha\wedge d\eta.
\end{align*}
Consequently $d\alpha=0$ if and only if $\partial T_\alpha=0$.
No mass bound for $\partial T_\alpha$ follows from local
integrability of $\alpha$ alone. Distributional supports are
understood in the same way as current supports.

On the complexified cotangent bundle, the action induced by $J$ has
eigenspaces for $i$ and $-i$, consisting respectively of the $(1,0)$
and $(0,1)$ covectors. Their exterior products
define forms of type $(a,b)$, even when $J$ is not integrable.
A real 2-form $\alpha$ is of type $(1,1)$ exactly when
$\alpha(JX,JY)=\alpha(X,Y)$, and is positive when in addition
$\alpha(X,JX)\geq0$ for every $X$.

A real current $T$ of dimension $2p$ has \emph{bidimension $(p,p)$}
if its complex-linear extension vanishes on all test forms of type
$(a,b)$ with $a+b=2p$ and $(a,b)\ne(p,p)$. Its \emph{bidegree} is
then $(m-p,m-p)$. A real $(p,p)$-form is \emph{strongly positive}
if, at each point, it is a finite sum with nonnegative real
coefficients of forms of the shape
\begin{align*}
    \left(\frac{i}{2}\right)^p
    \alpha_1\wedge\overline{\alpha_1}\wedge\cdots\wedge
    \alpha_p\wedge\overline{\alpha_p},
\end{align*}
where the $\alpha_j$ are complex covectors of type $(1,0)$.
A current of bidimension $(p,p)$ is \emph{positive} if it takes
nonnegative values on every smooth compactly supported strongly
positive $(p,p)$-form.

In particular, a positive real $(1,1)$-form $\alpha$ in
$L^1_{\mathrm{loc}}$ defines a positive current $T_\alpha$ of
dimension $2m-2$ and bidimension $(m-1,m-1)$. We write $T_\alpha$
when taking the boundary or a contraction of this current.

The weak pseudoholomorphic equation and compatibility of $j$ and
$\omega$ imply
\begin{align}\label{Equation: Positivity of the pullback}
    (u^*\omega)(JX,JY)&=(u^*\omega)(X,Y),\\
    (u^*\omega)(X,JX)&=h(duX,duX)\geq0.\notag
\end{align}
Thus $u^*\omega$ is a positive $(1,1)$-form in $L^1_{\mathrm{loc}}$.
Since the energy identity \eqref{Equation: Energy identity} uses only
the complex linearity of $du$, it applies here as well.

\subsection{Sobolev identities and analytic dimension reduction}
\label{Subsection: Sobolev identities and reduction}

Let $(N,j_N)$ be a compact complex manifold. For
$v\in W^{1,2}_{\mathrm{loc}}(M,N)$ and a smooth $1$-form $\beta$
on $N$, the identity $d(v^*\beta)=v^*(d\beta)$ holds in
distributions by Lemma~\ref{Lemma: Sobolev pullbacks}.

We recall that the Nijenhuis tensor $N_J$ associated with the almost complex structure $J$ is given by
\begin{align*}
    N_J(X,Y)=[JX,JY]-J[JX,Y]-J[X,JY]-[X,Y] \qquad\forall\, X,Y\in\mathfrak{X}(M).
\end{align*}
Moreover, we define the first-order differential operator $d_J^c$ on $C^{\infty}(M;\mathbb R)$ by
\begin{align*}
    d_J^c\psi:=-d\psi\circ J \qquad\forall\,\psi\in C^{\infty}(M;\mathbb R).
\end{align*}
\begin{lemma}[Nijenhuis annihilation]\label{Lemma: Nijenhuis annihilation}
    Let $(N,j_N)$ be a compact complex manifold. If
    $u\in W^{1,2}_{\mathrm{loc}}(M,N)$ is weakly pseudoholomorphic, then
    \begin{align*}
        du\bigl(N_J(X,Y)\bigr)=0 \qquad\mbox{ a.e. on } M,
    \end{align*}
    for all smooth vector fields $X,Y\in\mathfrak{X}(M)$.
\end{lemma}
\begin{proof}
    For a smooth real-valued function $\psi\in C^{\infty}(M;\mathbb R)$, the formula for the exterior derivative of a 1-form gives
    \begin{align*}
        dd_J^c \psi(X,Y)
        =-X(JY\psi)+Y(JX\psi)+d\psi(J[X,Y]).
    \end{align*}
    Since $J^2=-\mathrm{Id}$, evaluating the same formula at $JX,JY$
    gives
    \begin{align*}
        dd_J^c \psi(JX,JY)
        =JX(Y\psi)-JY(X\psi)+d\psi(J[JX,JY]).
    \end{align*}
    Subtracting the first identity from the second and grouping the
    differences of second derivatives into Lie brackets gives
    \begin{align}\label{Equation: Scalar Nijenhuis identity}
    \begin{split}
        dd_J^c \psi(JX,JY)-dd_J^c \psi(X,Y)
        &=d\psi\bigl([JX,Y]+[X,JY]+J[JX,JY]-J[X,Y]\bigr)\\
        &=d\psi\bigl(JN_J(X,Y)\bigr).
    \end{split}
    \end{align}
    For fixed smooth $J,X,Y$, every operation in \eqref{Equation: Scalar Nijenhuis identity} is a
    linear differential operator with smooth coefficients applied to $\psi$.
    In particular, evaluating the distributional 2-form $dd_J^c\psi$ on
    $X,Y$ means contracting its coefficients with the smooth
    coefficients of those fields. Consequently the identity extends
    from smooth functions to scalar distributions, and in particular to
    $W^{1,2}_{\mathrm{loc}}$ functions, by local mollification.

    Now let $f\in C^\infty(N;\mathbb R)$ and set $\psi=f\circ u$. The ordinary
    Sobolev chain rule gives $\psi\in W^{1,2}_{\mathrm{loc}}$. For almost every
    domain point and every tangent vector $V$, pseudoholomorphicity gives
    \begin{align*}
        d_J^c \psi(V)
        =-df_u\bigl(du(JV)\bigr)
        =-df_u\bigl(j_N(u)du(V)\bigr)
        =\bigl(u^*d_{j_N}^cf\bigr)(V).
    \end{align*}
    Both sides are locally $L^2$ 1-forms. Applying
    the Sobolev pullback identity above to the smooth target 1-form $d_{j_N}^cf$
    therefore gives the distributional identity
    \begin{align*}
        dd_J^c \psi=u^*(dd_{j_N}^cf).
    \end{align*}
    Its right-hand side is a locally $L^1$ 2-form. Since $j_N$ is integrable,
    $dd_{j_N}^cf$ is of type $(1,1)$, and complex linearity of $du$ makes
    its pullback of type $(1,1)$ almost everywhere. The left-hand side of
    \eqref{Equation: Scalar Nijenhuis identity} is therefore the zero distribution.
    We conclude that
    \begin{align*}
        df_u\bigl(du(JN_J(X,Y))\bigr)=0
    \end{align*}
    as a distribution. This expression is an $L^2_{\mathrm{loc}}$ function,
    so it vanishes almost everywhere.

    Choose a smooth embedding $\iota:N\hookrightarrow\mathbb R^L$ and apply
    the preceding argument to its finitely many coordinate functions
    $f_1,\ldots,f_L$. Outside the union of the resulting null sets,
    \begin{align*}
        d\iota_u\bigl(du(JN_J(X,Y))\bigr)=0.
    \end{align*}
    The differential of an embedding is injective on each tangent
    space, so $du(JN_J(X,Y))=0$ almost everywhere. Finally,
    \begin{align*}
        du(JN_J(X,Y))=j_N(u)\,du(N_J(X,Y)),
    \end{align*}
    and $j_N(u)$ is invertible, proving the assertion.

    On a coordinate neighbourhood we can choose a single exceptional
    null set for all $X,Y$. Indeed, apply the argument first to the
    finitely many pairs of coordinate vector fields. Tensoriality of
    $N_J$ then gives the conclusion for every pair of smooth fields
    outside that same null set.
\end{proof}

For a fixed weakly pseudoholomorphic map $u$, call a smooth vector field $X$ an
\emph{annihilator} if $du(X)=0$ almost everywhere on its domain.

We write $\mathcal L_XJ$ for the Lie derivative of $J$ along $X$;
explicitly, $(\mathcal L_XJ)Y=[X,JY]-J[X,Y]$.

\begin{lemma}[Annihilators and Lie derivatives]\label{Lemma: Annihilators and Lie derivatives}
    Annihilators form a module over smooth real functions and are
    preserved by $J$. If $X$ is an annihilator and $Y$ is any smooth
    vector field, then $(\mathcal L_XJ)Y$ is an annihilator.
\end{lemma}
\begin{proof}
    The module property follows from $du(aX)=a\,du(X)$ for a smooth
    real function $a$. Pseudoholomorphicity gives
    $du(JX)=j_N(u)du(X)$, proving preservation by $J$.

    For the remaining assertion, work in a coordinate neighbourhood
    $U$ on which $X$ is smooth and $Xu=0$ almost everywhere. Fix
    $U_0\Subset U$ and let $\delta>0$ be such that the local flow $F_t$ of
    $X$ is defined on $U_0$ and satisfies $F_t(U_0)\Subset U$ for
    every $|t|<\delta$. After decreasing $\delta$ if necessary, all the sets $F_t(U_0)$ lie in a fixed relatively compact neighbourhood $U_1\Subset U$, and the Jacobians and first derivatives of $F_t$ and $F_t^{-1}$ are uniformly bounded on the relevant domains.

    Regard $u$ as Euclidean-valued through a smooth target embedding.
    For a smooth ambient-valued function $v$, the fundamental theorem
    of calculus along flow lines gives
    \begin{align*}
        v\circ F_t-v=\int_0^t(Xv)\circ F_s\,ds.
    \end{align*}
    For each fixed $t$, this identity passes to a $W^{1,2}$ map by
    strong Sobolev approximation. Indeed, the uniform Jacobian bounds
    make composition with $F_s$ bounded on $L^2$, uniformly in $s$,
    and the integral on the right converges in $L^2(U_0)$ by
    Minkowski's inequality. Applied to $u$, whose derivative $Xu$
    vanishes, it yields
    \begin{align}\label{Equation: Invariance under the flow}
        u\circ F_t=u \qquad\mbox{ in } L^2(U_0)
    \end{align}
    for every fixed $|t|<\delta$.

    The pullback of the endomorphism $J$ is
    \begin{align*}
        (F_t^*J)_x
        =DF_t(x)^{-1}\,J_{F_t(x)}\,DF_t(x) \qquad\forall\,x\in U_0.
    \end{align*}
    The Sobolev chain rule for a smooth domain diffeomorphism, together
    with the equation for $u$ at $F_t(x)$, gives
    \begin{align*}
        d(u\circ F_t)\circ F_t^*J
        =j_N(u\circ F_t)\circ d(u\circ F_t)
        \qquad\mbox{ a.e. on } U_0.
    \end{align*}
    Thus, by \eqref{Equation: Invariance under the flow}, for each fixed $t$,
    \begin{align*}
        du\circ F_t^*J=j_N(u)\circ du
        \qquad\mbox{ a.e. on } U_0.
    \end{align*}
    Choose nonzero $t_\ell\to0$ with $|t_\ell|<\delta$. Discard the
    union of the null sets for these countably many equations and for
    the original pseudoholomorphicity equation. At every remaining
    point, subtraction and division by $t_\ell$ give
    \begin{align*}
        du\circ\frac{F_{t_\ell}^*J-J}{t_\ell}=0.
    \end{align*}
    The smooth tensor quotients converge locally uniformly to
    $\mathcal L_XJ$. Passing to the limit at such a point yields
    $du\circ\mathcal L_XJ=0$. Evaluating this tensor identity on any smooth vector field $Y$ gives
    \begin{align*}
        du\bigl((\mathcal L_XJ)Y\bigr)=0
        \qquad\mbox{ a.e. on } U_0.
    \end{align*}
    Since $U_0\Subset U$ was arbitrary, the same identity holds almost everywhere on $U$. Thus, $(\mathcal L_XJ)Y$ is an annihilator.
\end{proof}

\begin{lemma}[Analytic reduction by two real dimensions]
    \label{Lemma: Analytic dimension reduction}
    Let $(M^{2m},J)$ be real analytic, with $m\geq2$, let $(N,j_N)$
    be a compact complex manifold, and let
    $u\in W^{1,2}_{\mathrm{loc}}(M,N)$ be weakly pseudoholomorphic.
    Suppose a real-analytic annihilator $X$ satisfies $X(p)\ne0$.
    Then, near $p$, there are a bounded open cube $Q\subset\mathbb R^{2m-2}$ containing $0$, real-analytic coordinates $(t,s,q)$ with $q\in Q$, a real-analytic almost complex structure $I$ on $Q$, and $f\in W^{1,2}_{\mathrm{loc}}(Q,N)$
    such that
    \begin{align}\label{Equation: Reduced pseudoholomorphic map}
        u(t,s,q)=f(q) \quad\mbox{ a.e.},
        \qquad df\circ I=j_N\circ df\quad\text{a.e. on }Q.
    \end{align}
    If a smooth target 2-form $\beta$ satisfies $d(f^*\beta)=0$,
    then $d(u^*\beta)=0$ in the sense of distributions on this neighbourhood.
\end{lemma}
\begin{proof}
    Set $d=2m-2$, and fix a smooth embedding of $N$ into
    $\mathbb R^L$. We identify $N$ with its image and take weak
    derivatives componentwise in this embedding. Throughout the
    proof, we use the same notation for maps and tensors and their
    expressions in local coordinates.

    We work in an open neighbourhood of $p$ in which $X$ never vanishes.
    The analytic flow-box theorem gives coordinates $(t,y)$ on a
    product $I_0\times Y$, where $I_0\subset\mathbb R$ is a bounded
    open interval containing zero and $Y\subset\mathbb R^{d+1}$ is
    a bounded open neighbourhood of zero, such that $p=(0,0)$ and
    $X=\partial_t$. Choose this product with closure contained in
    a larger flow-box chart. The Sobolev chain rule for smooth
    coordinate diffeomorphisms then gives
    $u\in W^{1,2}(I_0\times Y,N)$ after shrinking the product if
    necessary. Since $X$ is an annihilator,
    \begin{align*}
        \partial_tu=du(X)=0
        \qquad\mbox{ a.e. on } I_0\times Y.
    \end{align*}
    By Sobolev slicing, for almost every $y\in Y$ the slice
    $t\mapsto u(t,y)$ has an absolutely continuous representative
    whose derivative equals $\partial_tu(t,y)$ almost everywhere.
    For this representative and every $t,t'\in I_0$,
    \begin{align*}
        u(t,y)-u(t',y)
        =\int_{t'}^t\partial_tu(\tau,y)\,d\tau=0.
    \end{align*}
    Thus almost every slice has a constant representative.

    Choose $\rho\in C_c^\infty(I_0)$ with $\rho\geq0$ and
    $\int_{I_0}\rho(t)\,dt=1$, and define
    \begin{align*}
        v(y):=\int_{I_0}\rho(t)u(t,y)\,dt.
    \end{align*}
    For almost every $y$, this average equals the constant value
    of the slice. Consequently, $v(y)\in N$ almost everywhere and
    $u(t,y)=v(y)$ almost everywhere on $I_0\times Y$.
    For each $r\in\{1,\ldots,d+1\}$, set
    \begin{align*}
        G_r(y)=\int_{I_0}\rho(t)\partial_{y_r}u(t,y)\,dt.
    \end{align*}
    For every $\varphi\in C_c^\infty(Y)$, the function
    $(t,y)\mapsto\rho(t)\varphi(y)$ is an admissible test function
    on $I_0\times Y$. Therefore, componentwise,
    \begin{align*}
        \int_Yv\,\partial_{y_r}\varphi\,dy
        &=\int_{I_0\times Y}
        \rho(t)u(t,y)\partial_{y_r}\varphi(y)\,dt\,dy\\
        &=-\int_{I_0\times Y}
        \rho(t)\partial_{y_r}u(t,y)\varphi(y)\,dt\,dy\\
        &=-\int_YG_r\varphi\,dy.
    \end{align*}
    Moreover, Jensen's inequality gives
    \begin{align*}
        \int_Y\left(|v|^2+\sum_{r=1}^{d+1}|G_r|^2\right)dy
        &\leq\int_{I_0\times Y}\rho(t)
        \left(|u|^2+\sum_{r=1}^{d+1}
        |\partial_{y_r}u|^2\right)dt\,dy\\
        &\leq\|\rho\|_{L^\infty(I_0)}
        \|u\|_{W^{1,2}(I_0\times Y)}^2.
    \end{align*}
    It follows that $v\in W^{1,2}(Y,N)$ and
    $\partial_{y_r}v=G_r$. To identify the derivatives of the
    factorisation $u(t,y)=v(y)$, take $\Psi\in C_c^\infty(I_0\times Y)$.
    Fubini's theorem and the weak derivative identity for $v$ give
    \begin{align*}
        \int_{I_0\times Y}u\,\partial_{y_r}\Psi\,dt\,dy
        &=\int_{I_0}\int_Yv(y)
        \partial_{y_r}\Psi(t,y)\,dy\,dt\\
        &=-\int_{I_0\times Y}
        \partial_{y_r}v(y)\Psi(t,y)\,dt\,dy.
    \end{align*}
    By uniqueness of weak derivatives,
    \begin{align*}
        \partial_{y_r}u(t,y)=\partial_{y_r}v(y)
        \qquad\mbox{ for a.e. }(t,y).
    \end{align*}
    Write
    \begin{align*}
        J\partial_t=a(t,y)\partial_t+
        \sum_{r=1}^{d+1}b_r(t,y)\partial_{y_r},
    \end{align*}
    where all coefficients are real analytic. The vector
    $b(t,y)=(b_1(t,y),\ldots,b_{d+1}(t,y))$ never vanishes.
    Indeed, if it vanished at a point, then at that point
    \begin{align*}
        -\partial_t=J^2\partial_t
        =J(a\partial_t)=a^2\partial_t,
    \end{align*}
    which is impossible because $a$ is real. Pseudoholomorphicity
    and $\partial_tu=0$ imply
    \begin{align*}
        0=j_N(u)\partial_tu
        =du(J\partial_t)
        =\sum_{r=1}^{d+1}b_r(t,y)\partial_{y_r}v(y)
        \qquad\mbox{ for a.e. }(t,y).
    \end{align*}
    By Fubini's theorem, for almost every $y$ all the vectors
    $\partial_{y_r}v(y)$ are defined and the last identity holds
    for almost every $t$. For each such $y$, the function
    \begin{align*}
        t\mapsto
        \sum_{r=1}^{d+1}b_r(t,y)\partial_{y_r}v(y)
    \end{align*}
    is continuous and vanishes almost everywhere, so it vanishes
    everywhere on $I_0$. In particular, the real-analytic field
    \begin{align*}
        B(y)=\sum_{r=1}^{d+1}b_r(0,y)\partial_{y_r}
    \end{align*}
    is nonvanishing and satisfies
    \begin{align*}
        dv_y(B(y))
        =\sum_{r=1}^{d+1}b_r(0,y)\partial_{y_r}v(y)=0
        \qquad\mbox{ for a.e. }y\in Y.
    \end{align*}

    Apply the analytic flow-box theorem to $B$ near $y=0$.
    It gives an analytic coordinate parametrisation
    $\psi:I_1\times Q\rightarrow Y'\subset Y$, where
    $I_1\subset\mathbb R$ is a bounded open interval containing
    zero, $Q\subset\mathbb R^d$ is a bounded open cube containing
    zero, and
    \begin{align*}
        \psi(0,0)=0,
        \qquad\partial_s\psi(s,q)=B(\psi(s,q)).
    \end{align*}
    Choose this product with closure contained in a larger
    flow-box chart whose image lies in $Y$. Then
    $w:=v\circ\psi\in W^{1,2}(I_1\times Q,N)$, and the Sobolev
    chain rule gives
    \begin{align*}
        \partial_sw(s,q)
        =dv_{\psi(s,q)}\bigl(\partial_s\psi(s,q)\bigr)
        =dv_{\psi(s,q)}\bigl(B(\psi(s,q))\bigr)=0
        \qquad\mbox{ a.e.}
    \end{align*}
    Applying the same slicing and averaging argument in the $s$
    variable gives $f\in W^{1,2}(Q,N)$ with $w(s,q)=f(q)$
    almost everywhere. In the analytic coordinates $(t,s,q)$ on
    $I_0\times I_1\times Q$ we thus have $X=\partial_t$ and
    $u(t,s,q)=f(q)$ almost everywhere. With $\pi(t,s,q)=q$, the
    weak derivative identities give
    \begin{align*}
        du=df\circ d\pi\qquad\text{a.e. on } I_0\times I_1\times Q.
    \end{align*}

    We next construct the almost complex structure on $Q$.
    Let $\sigma(q):=(0,0,q)$ and
    $a_0(q):=a(0,\psi(0,q))$. At $t=0$, the $y$ component of
    $J\partial_t$ is $B$, which becomes $\partial_s$ in the
    second flow box. Therefore, along $\sigma(Q)$,
    \begin{align*}
        J\partial_t=a_0\partial_t+\partial_s.
    \end{align*}
    Applying $J$ pointwise to
    $\partial_s=J\partial_t-a_0\partial_t$ gives, along the
    same section,
    \begin{align*}
        J\partial_s
        &=J^2\partial_t-a_0J\partial_t\\
        &=-\partial_t-a_0(a_0\partial_t+\partial_s)\\
        &=-(1+a_0^2)\partial_t-a_0\partial_s.
    \end{align*}
    In particular, $J$ preserves
    $\ker d\pi=\operatorname{span}\{\partial_t,\partial_s\}$
    along $\sigma(Q)$.

    On the product $I_0\times I_1\times Q$, write
    \begin{align*}
        J\partial_{q_i}
        =c_i(t,s,q)\partial_t+d_i(t,s,q)\partial_s
        +\sum_{k=1}^d A_i^k(t,s,q)\partial_{q_k},
        \qquad 1\leq i\leq d,
    \end{align*}
    with real-analytic coefficients. Define
    \begin{align*}
        I_q(\partial_{q_i})
        =\sum_{k=1}^d A_i^k(0,0,q)\partial_{q_k}.
    \end{align*}
    Since the coefficients $A_i^k$ are real analytic, this defines a
    real-analytic endomorphism of $TQ$. Equivalently,
    \begin{align*}
        I_q=d\pi_{\sigma(q)}\circ J_{\sigma(q)}\circ d\sigma_q.
    \end{align*}
    To compute its square, apply $d\pi$ to
    $J^2\partial_{q_i}=-\partial_{q_i}$ at $\sigma(q)$.
    The two preceding formulas for $J\partial_t$ and
    $J\partial_s$ show that their $q$ components vanish there.
    Consequently,
    \begin{align*}
        -\partial_{q_i}
        &=d\pi_{\sigma(q)}
        \left(J^2\partial_{q_i}\big|_{\sigma(q)}\right)\\
        &=d\pi_{\sigma(q)}
        \left(\left.c_iJ\partial_t+d_iJ\partial_s
        +\sum_{k=1}^d A_i^kJ\partial_{q_k}
        \right|_{\sigma(q)}\right)\\
        &=\sum_{\ell=1}^d\left(
        \sum_{k=1}^d A_i^k(0,0,q)A_k^\ell(0,0,q)
        \right)\partial_{q_\ell}.
    \end{align*}
    Thus $\sum_k A_i^k(0,0,q)A_k^\ell(0,0,q)
    =-\delta_i^\ell$, and
    \begin{align*}
        I_q^2(\partial_{q_i})
        &=\sum_{k,\ell=1}^d
        A_i^k(0,0,q)A_k^\ell(0,0,q)\partial_{q_\ell}=-\partial_{q_i}.
    \end{align*}
    Hence $I^2=-\mathrm{Id}$, as required.

    We now verify pseudoholomorphicity of $f$. Evaluating the
    equation for $u$ on $\partial_{q_i}$, and using the weak
    derivative identities established above, gives
    \begin{align*}
        \sum_{k=1}^d A_i^k(t,s,q)\partial_{q_k}f(q)
        &=du_{(t,s,q)}\bigl(J\partial_{q_i}\bigr)\\
        &=j_N(u(t,s,q))\partial_{q_i}u(t,s,q)\\
        &=j_N(f(q))\partial_{q_i}f(q)
        \qquad\mbox{ for a.e. }(t,s,q).
    \end{align*}
    By Fubini's theorem, there is a set $Q_*\subset Q$ of full
    measure such that, for each $q\in Q_*$, all the vectors
    $\partial_{q_k}f(q)$ are defined and the preceding identities
    hold for almost every $(t,s)$, simultaneously for
    $i=1,\ldots,d$. Fix such a $q$. For each $i$, the function
    \begin{align*}
        (t,s)\mapsto
        \sum_{k=1}^d A_i^k(t,s,q)\partial_{q_k}f(q)
        -j_N(f(q))\partial_{q_i}f(q)
    \end{align*}
    is continuous and vanishes almost everywhere on
    $I_0\times I_1$, so it vanishes everywhere. Evaluating at
    $(t,s)=(0,0)$ therefore gives
    \begin{align*}
        df_q\bigl(I_q\partial_{q_i}\bigr)
        &=\sum_{k=1}^d A_i^k(0,0,q)\partial_{q_k}f(q)\\
        &=j_N(f(q))\partial_{q_i}f(q)
        =j_N(f(q))df_q(\partial_{q_i}).
    \end{align*}
    Since this holds for every coordinate vector,
    $df\circ I=j_N(f)\circ df$ almost everywhere on $Q$.
    Thus evaluation at $(0,0)$ uses only the continuity of the
    coefficients in $(t,s)$ after factorisation.

    Finally, for a smooth target $2$-form $\beta$,
    the identity $du=df\circ d\pi$ yields
    $u^*\beta=\pi^*(f^*\beta)$ almost everywhere in the product
    coordinates. The coefficients of $f^*\beta$ belong to
    $L^1_{\mathrm{loc}}(Q)$. For such a form $\gamma$, pullback
    by the coordinate projection commutes with $d$ in distributions.
    Indeed, if $g(q)$ is any coefficient of $\gamma$, then
    $g\circ\pi$ has zero derivatives in $t,s$ and
    \begin{align*}
        \langle\partial_{q_i}(g\circ\pi),\varphi\rangle
        =\left\langle\partial_{q_i}g,\,
           \int_{I_0\times I_1}\varphi(t,s,\cdot)\,dt\,ds
         \right\rangle
    \end{align*}
    for every compactly supported smooth test function $\varphi$,
    by Fubini's theorem. Hence $d(f^*\beta)=0$ implies
    $d(u^*\beta)=0$ in these coordinates and, by change of
    variables, on the corresponding neighbourhood of $p$.
\end{proof}

\begin{remark}
    The coordinate plane field
    $\ker d\pi=\spn\{\partial_t,\partial_s\}$ agrees with
    $\spn\{X,JX\}$ along the section $\sigma(Q)$, although the
    two plane fields may differ away from this section.
    The invariance along $\sigma(Q)$ suffices to define $I$, so the
    construction requires no involutivity assumption on
    $\spn\{X,JX\}$.
\end{remark}

\subsection{The integrable sphere-valued case}
\label{Subsection: Integrable sphere valued case}

The integrable case belongs to the classical theory of holomorphic
chains developed by King \cite{King} and Harvey--Shiffman
\cite{HarveyShiffman}. We recall the precise consequence needed
for our domain induction and explain its deduction from that theory.

A holomorphic $r$-chain is a locally finite sum
$\sum_\nu n_\nu\llbracket V_\nu\rrbracket$, where
$n_\nu\in\mathbb Z$ and the $V_\nu$ are irreducible complex
analytic subvarieties of pure complex dimension $r$.
The current $\llbracket V_\nu\rrbracket$ integrates over the
regular locus with its complex orientation. Such chains are closed;
their positive and negative parts, obtained by separating the signs
of the coefficients, are also closed.

King's theorem \cite[Theorem~5.2.1]{King} characterises positive
closed locally integer rectifiable currents of pure complex
bidimension as positive holomorphic chains. Harvey--Shiffman
\cite{HarveyShiffman} allow signed currents, with the
additional support hypothesis
$\Hh^{2r+1}(\spt R)=0$ for a current of bidimension $(r,r)$.
In the hypersurface case, Shiffman \cite{ShiffmanCodimensionOne}
removed the support hypothesis. Thus a closed
locally integer rectifiable current of bidimension $(m-1,m-1)$ on
a complex $m$-manifold is a holomorphic $(m-1)$-chain, without a
positivity or support assumption.

\begin{proposition}[The classical integrable case]
    \label{Proposition: Integrable sphere valued case}
    On a complex manifold of any complex dimension, every weakly
    holomorphic map $u\in W^{1,2}_{\mathrm{loc}}(M,\cp^1)$
    satisfies $d(u^*\omega)=0$ in the sense of distributions.
\end{proposition}
\begin{proof}[Deduction from the classical holomorphic-chain theorem]
    Closedness is automatic in complex dimension one. Otherwise,
    work in a relatively compact holomorphic coordinate ball $B$
    on which $u\in W^{1,2}(B,\cp^1)$, and on $B$ let
    \begin{align*}
        T_u:=T_{u^*\omega},
        \qquad \Gamma_u:=\partial T_u.
    \end{align*}
    We use the Sobolev level-current theorem of
    Alberti--Baldo--Orlandi
    \cite[Theorem~3.8 and Section~7.5]{AlbertiBaldoOrlandi},
    identifying $\cp^1$ with the oriented 2-sphere. With a fixed
    Borel representative of $u$ and its countable Lipschitz
    decomposition, almost every $a\in\cp^1$ determines a locally
    integer rectifiable $(2m-2)$-current $C_a$ carried by that
    representative's level set. The decomposition covers almost
    all of $B$ by Borel pieces on which $u$ agrees with a Lipschitz
    map and the approximate differentials agree almost everywhere.
    The level currents use the coarea orientation, and satisfy
    \begin{align}\label{Equation: Common fibre boundary}
        T_u=\int_{\cp^1}C_a\,\omega(a),\qquad
        \partial C_a=\Gamma_u\quad\text{for almost every }a.
    \end{align}
    The first equality holds on every compactly supported test form.
    The integral is absolutely convergent, since the coarea estimate gives
    $\int_{\cp^1}\|C_a\|(K)\,\omega(a)
    \leq C\int_K|du|^2\,d\vol_g$ for compact $K\subset B$.
    The common boundary is a distribution; the cited construction
    imposes no mass assumption on it. The normalisation
    $\int_{\cp^1}\omega=1$ identifies that boundary with
    $\partial T_u$ in the displayed formula.

    By the coarea formula, for almost every $a$ the approximate
    differential is complex linear and has rank two at
    $\|C_a\|$-almost every point. The approximate tangent plane of
    the level at such a point is therefore $\ker du$, a complex hyperplane. On its
    Hermitian orthogonal complement, $du$ is a complex-linear
    isomorphism to the target tangent line. Consequently the coarea
    orientation is the positive complex orientation, and $C_a$ is
    positive of bidimension $(m-1,m-1)$. For distinct admissible
    values $a,b$, the measures $\|C_a\|$ and $\|C_b\|$ are mutually
    singular, since they are carried by disjoint Borel level sets.

    Choose two such values. The difference $R:=C_a-C_b$ is locally
    integer rectifiable, has pure bidimension $(m-1,m-1)$, and is
    closed by \eqref{Equation: Common fibre boundary}. Shiffman's
    theorem therefore makes $R$ a holomorphic $(m-1)$-chain.
    Write $R=R^+-R^-$ by separating its positive and negative
    coefficients. Relative to the canonical complex orientation,
    mutual singularity of the two fibre measures identifies these
    parts uniquely as $R^+=C_a$ and $R^-=C_b$: the multiplicity is
    positive on the carrier of $C_a$ and negative on that of $C_b$,
    with no cancellation. Since both holomorphic chains $R^+$ and
    $R^-$ are closed, $\Gamma_u=\partial C_a=0$.
    This proves $d(u^*\omega)=0$ on $B$, hence on $M$ by locality.
\end{proof}

\subsection{Finite analytic annihilators and their common zero set}
\label{Subsection: Analytic annihilators}

A \emph{real-analytic subset} of a real-analytic manifold is a set
locally defined by the simultaneous vanishing of finitely many
real-analytic functions. A \emph{semianalytic set} is locally a finite
union of sets defined by finitely many real-analytic equations and
strict inequalities. A set is \emph{subanalytic} if locally it is the
image, under a coordinate projection, of a relatively compact
semianalytic set in a product with a Euclidean space. We use these
notions only on real-analytic manifolds. A point of a subanalytic
set is called \emph{regular} if the set agrees, in an ambient
neighbourhood of that point, with an embedded real-analytic
submanifold. The dimension of a subanalytic set is the maximum of
the dimensions of its strata in a real-analytic stratification;
the stratification facts needed below are recalled in
Remark~\ref{Remark: Analytic stratification}. We set $\dim\varnothing=-\infty$.
For these definitions and their local properties, see
\cite{BierstoneMilman}.

\begin{lemma}\label{Lemma: Finite analytic annihilators}
    Let $(M^{2m},J)$ be a real-analytic almost complex manifold, let
    $(N,j_N)$ be a compact complex manifold, and let
    $u\in W^{1,2}_{\mathrm{loc}}(M,N)$ be weakly pseudoholomorphic.
    Every $p\in M$ has a connected analytic coordinate neighbourhood
    $V$ on which at least one of the following holds.
    \begin{enumerate}
        \item $J$ is integrable;
        \item there are finitely many real-analytic annihilators
        $X_1,\ldots,X_r$, not all identically zero, and real-analytic
        functions $c_{iab}$ on $V$ satisfying
        \begin{align}\label{Equation: Finite annihilator relations}
            (\mathcal L_{X_i}J)\partial_{x^a}
            =\sum_{b=1}^r c_{iab}X_b
            \qquad\forall\,1\leq i\leq r,\,\forall\, 1\leq a\leq2m.
        \end{align}
        In this case the common zero set of $X_1,\ldots,X_r$
        \begin{align*}
            Z:=\{x\in V:X_1(x)=\cdots=X_r(x)=0\}
        \end{align*}
        is a proper relatively closed real-analytic subset of $V$.
    \end{enumerate}
\end{lemma}
\begin{proof}
    Let $\mathcal O_p^{\mathbb R}$ denote the ring of germs at $p$ of
    real-valued real-analytic functions. Two representatives define the
    same germ if they agree on a neighbourhood of $p$. In fixed analytic
    coordinates, germs of analytic vector fields form the finite free module
    $(\mathcal O_p^{\mathbb R})^{2m}$. Let $\mathcal A_p$ be its subset
    consisting of germs that admit a representative $X$ with $du(X)=0$
    almost everywhere on some neighbourhood of $p$.

    This condition is independent of the representative. Since sums
    and analytic multiples of annihilators still annihilate $du$ on
    their common domain, $\mathcal A_p$ is a submodule. The ring
    $\mathcal O_p^{\mathbb R}$ is Noetherian because it is isomorphic
    to the ring of convergent real power series in $2m$ variables. Since a submodule of a finite free module
    over a Noetherian ring is finitely generated, $\mathcal A_p$ has
    finitely many generators.

    If $\mathcal A_p=0$, Lemma~\ref{Lemma: Nijenhuis annihilation} shows that each of the
    analytic fields $N_J(\partial_{x^a},\partial_{x^b})$ has zero germ.
    Each field therefore vanishes identically on a neighbourhood of $p$.
    There are finitely many pairs $(a,b)$, so all of them vanish on one
    common neighbourhood. Hence $N_J=0$ there, and the
    Newlander--Nirenberg theorem gives the first alternative
    \cite{NewlanderNirenberg}.

    Suppose instead that $\mathcal A_p\ne0$. Choose analytic
    representatives $X_1,\ldots,X_r$ of a finite generating family.
    By the definition of $\mathcal A_p$, each representative annihilates
    $du$ on a sufficiently small neighbourhood. Intersecting these
    finitely many neighbourhoods gives a common domain for all fields
    and annihilation identities. Lemma~\ref{Lemma: Annihilators and Lie derivatives} then implies that
    each germ $(\mathcal L_{X_i}J)\partial_{x^a}$ belongs to
    $\mathcal A_p$. Express it in the chosen generators, obtaining
    coefficient germs $c_{iab}$. Choose analytic representatives of
    these coefficients and shrink the neighbourhood until every identity
    \eqref{Equation: Finite annihilator relations} holds pointwise. Only finitely many
    fields, coefficients, and identities occur, so all these choices
    are valid on one connected coordinate neighbourhood $V$ of $p$.

    The set $Z$ is the zero set of the $2mr$ analytic scalar components
    of the fields. At least one component has a nonzero germ. Since
    $V$ is connected, the analytic identity principle implies that
    this component cannot vanish on a nonempty open subset of $V$. Consequently
    $Z$ is relatively closed, proper, and has empty interior. The
    argument requires only generation at the single stalk $p$ and
    extension of the finitely many identities above to a common neighbourhood.
\end{proof}

The next lemma gives analytic propagation of the common zero set,
with the parameter dependence needed in higher dimension.

\begin{lemma}[Analytic propagation of the common zero set]\label{Lemma: Parametric analytic propagation}
    Let $X_i,J,Z$ satisfy the second alternative of
    Lemma~\ref{Lemma: Finite analytic annihilators}. Let
    $c(s_1,\ldots,s_k)$ be a real-analytic map from a neighbourhood of $0\in\mathbb R^k$ into $Z$, where $k\geq1$. The Cauchy problem
    \begin{align}\label{Equation: Analytic Cauchy problem}
    \begin{cases}
        \partial_t\phi(s,t)=J(\phi(s,t))\partial_{s_1}\phi(s,t),\\
        \phi(s,0)=c(s)
    \end{cases}
    \end{align}
    has a unique local analytic solution $\phi$ from an open neighbourhood of $(0,0)$ in $\mathbb R^{k+1}$ into $Z$. If $c$ is an immersion at $0$ and
    $J_{c(0)}\partial_{s_1}c(0)\notin\operatorname{im}dc_0$, then $d\phi$ has rank $k+1$ throughout a neighbourhood of
    $(0,0)$. In particular, $\phi$, viewed as a map into the
    ambient manifold, is an immersion on that neighbourhood.
\end{lemma}
\begin{proof}
    In the coordinates of $V$, equation~\eqref{Equation: Analytic Cauchy problem} is an analytic
    first-order system solved for every component of $\partial_t\phi$.
    The hypersurface $t=0$ is therefore noncharacteristic. The
    Cauchy--Kowalevski theorem, applied with all the variables
    $s_1,\ldots,s_k$ tangential to this hypersurface, gives a unique
    solution jointly analytic in $(s,t)$ on a neighbourhood of $(0,0)$.
    In particular, the solution depends analytically on
    $s_2,\ldots,s_k$, although no derivatives in these variables
    appear in the equation. Shrink the neighbourhood so that
    $\phi$ takes values in $V$.

    Regard $J$ as an analytic matrix and $X_i$ as analytic maps into
    $\mathbb R^{2m}$ in this coordinate trivialisation. Write $DX_i$
    for the coordinate Jacobian and
    $D_XJ=\sum_aX^a\partial_{x^a}J$. Put
    \begin{align*}
        G_i(s,t):=X_i(\phi(s,t)),\qquad
        v(s,t):=\partial_{s_1}\phi(s,t).
    \end{align*}
    The coordinate formula for the Lie derivative is
    \begin{align*}
        (\mathcal L_XJ)v=(D_XJ)v-DX(Jv)+JDX(v).
    \end{align*}
    The chain rule and \eqref{Equation: Analytic Cauchy problem} give
    \begin{align*}
        \partial_tG_i
        &=DX_i(\phi)\partial_t\phi
          =DX_i(\phi)J(\phi)v,\\
        \partial_{s_1}G_i
        &=DX_i(\phi)\partial_{s_1}\phi
          =DX_i(\phi)v.
    \end{align*}
    The formula for the Lie derivative and \eqref{Equation: Finite annihilator relations} imply
    \begin{align}\label{Equation: Homogeneous propagation system}
        \partial_tG_i-J(\phi)\partial_{s_1}G_i
        &=DX_i(\phi)J(\phi)v
          -J(\phi)DX_i(\phi)v\notag\\
        &=(D_{X_i}J)(\phi)v
          -(\mathcal L_{X_i}J)(\phi)v\notag\\
        &=\sum_{a=1}^{2m}G_i^a
          (\partial_{x^a}J)(\phi)v
          -\sum_{b=1}^r
          \left(
              \sum_{a=1}^{2m}c_{iab}(\phi)v^a
          \right)G_b\\
        &=\sum_{b=1}^r B_{ib}(s,t)G_b.\notag
    \end{align}
    To specify the coefficient matrices, define the linear map
    $A(s,t):\mathbb R^{2m}\rightarrow\mathbb R^{2m}$ by
    \begin{align*}
        A(s,t)w&:=\sum_{a=1}^{2m}w^a
        (\partial_{x^a}J)(\phi(s,t))v(s,t),\\
        B_{ib}(s,t)&:=\delta_{ib}A(s,t)
        -\left(\sum_{a=1}^{2m}c_{iab}(\phi(s,t))v^a(s,t)\right)
        \mathrm{Id}.
    \end{align*}
    Since $J$, the functions $c_{iab}$, and $\phi$ are analytic,
    $J(\phi)$ and $B_{ib}$ are prescribed analytic coefficients.
    With these coefficients, consider the Cauchy problem for functions
    $H_i$ with values in $\mathbb R^{2m}$:
    \begin{align*}
        \begin{cases}
            \displaystyle
            \partial_tH_i
            =J(\phi)\partial_{s_1}H_i
             +\sum_{b=1}^r B_{ib}H_b,
             &1\leq i\leq r,\\
            H_i(s,0)=0,
             &1\leq i\leq r.
        \end{cases}
    \end{align*}
    Taken together, these equations form a linear system of $2mr$
    scalar equations with analytic coefficients, solved for every
    component of the derivative in $t$. Cauchy--Kowalevski
    uniqueness therefore applies near $(0,0)$.

    The functions $H_i=G_i$ solve this system by
    \eqref{Equation: Homogeneous propagation system}, and their
    initial values satisfy
    \begin{align*}
        G_i(s,0)=X_i(\phi(s,0))=X_i(c(s))=0,
    \end{align*}
    because $c(s)\in Z$. The functions $H_i\equiv0$ solve the
    same system with the same initial data. Uniqueness implies
    $G_i\equiv0$ for every $i$, on a smaller neighbourhood of
    $(0,0)$. Since
    \begin{align*}
        Z=\{x\in V:X_i(x)=0\text{ for every }1\leq i\leq r\},
    \end{align*}
    it follows that $\phi$ takes values in $Z$. Uniqueness among
    analytic solutions taking values in $Z$ follows from uniqueness
    for the original Cauchy problem in the ambient coordinate
    neighbourhood $V$.

    At $(0,0)$, the first $k$ columns of $d\phi$ are the columns of
    $dc_0$, and the additional column is
    $J_{c(0)}\partial_{s_1}c(0)$. The stated hypotheses make these
    $k+1$ columns linearly independent. Hence some $(k+1)$-minor of
    $d\phi$ is nonzero at $(0,0)$ and, by continuity, remains nonzero
    on a neighbourhood. This proves the immersion assertion.
\end{proof}

\begin{remark}\label{Remark: Analytic stratification}
    Let $U$ be an open subset of a real-analytic manifold, and let
    $E\subset U$ be relatively closed and subanalytic. By a
    \emph{real-analytic stratification} of $E$ we mean a partition
    \begin{align*}
        E=\bigsqcup_{S\in\mathcal S}S
    \end{align*}
    with the following properties:
    \begin{enumerate}
        \item Each member $S\in\mathcal S$, called a
        \emph{stratum}, is a nonempty connected embedded
        real-analytic submanifold of $U$ without boundary, and
        is a subanalytic subset of $U$.

        \item The family $\mathcal S$ is \emph{locally finite
        in $U$}: every point of $U$ has an open neighbourhood
        meeting only finitely many strata.

        \item The \emph{frontier condition} holds: for any
        distinct strata $A,B\in\mathcal S$,
        \begin{align*}
            A\cap\overline B\ne\varnothing
            \quad\Rightarrow\quad
            A\subset\overline B\smallsetminus B
            \quad\mbox{and}\quad
            \dim A<\dim B.
        \end{align*}
        Here and below, closures are taken relative to $U$.
    \end{enumerate}
    Such stratifications exist for relatively closed subanalytic
    sets; see \cite{BierstoneMilman}. Since $E$ is relatively closed,
    $\overline B\subset E$ for every stratum $B$. The frontier
    condition therefore says that
    $\overline B\smallsetminus B$ is a union of entire strata,
    each of dimension strictly smaller than $\dim B$.

    For an integer $r\geq0$, the \emph{$r$-skeleton} of this
    stratification is
    \begin{align*}
        E^{(r)}:=\bigcup_{\substack{S\in\mathcal S\\\dim S\leq r}}S.
    \end{align*}
    This set is relatively closed in $U$. Indeed, let
    $x_\nu\in E^{(r)}$ converge to $x\in U$. Since $E$ is
    relatively closed, $x$ belongs to a unique stratum $A$.
    By local finiteness, after passing to a subsequence all
    $x_\nu$ belong to one stratum $B$ with $\dim B\leq r$.
    Hence $x\in A\cap\overline B$. If $A=B$, then
    $\dim A\leq r$; otherwise, the frontier condition gives
    $\dim A<\dim B\leq r$. In either case, $x\in E^{(r)}$.
    Moreover, $E^{(r)}$ is subanalytic because, near each point
    of $U$, it is a finite union of subanalytic strata.

    We also record why strata of maximal dimension consist
    of regular points of $E$. Let $B$ be such a stratum and
    let $x\in B$. If $x\in\overline A$ for a distinct stratum
    $A$, the frontier condition would imply
    $\dim B<\dim A$, contradicting maximality. Choose an
    ambient neighbourhood $W$ of $x$ meeting only finitely
    many strata. Since $x$ belongs to the closure of none of
    these strata other than $B$, we can shrink $W$ to an open
    neighbourhood $V$ of $x$ avoiding all the other strata.
    Consequently,
    \begin{align*}
        E\cap V=B\cap V.
    \end{align*}
    Thus $E$ is an embedded real-analytic submanifold near
    $x$, which is the meaning of regularity used here.

    Finally, only finitely many strata meet a relatively compact
    open set $V\Subset U$, since $\overline V$ can be covered by
    finitely many neighbourhoods, each meeting only finitely many
    strata.
    If $U$ carries an almost complex structure $J$, the definition
    of the stratification does not require $T_xS$ to be $J_x$-invariant.
    Any such invariance must therefore be proved separately.
\end{remark}

\begin{lemma}
    \label{Lemma: Geometry of the common zero set}
    Let $m\geq2$ and let $Z$ be the common zero set in the second
    alternative of Lemma~\ref{Lemma: Finite analytic annihilators}.
    Then the following hold.
    \begin{enumerate}
        \item at every regular point of $Z$, its tangent space is
        $J$-invariant; in particular $\dim Z\leq2m-2$;
        \item every real-analytic embedded submanifold
        $K^{2m-3}\subset Z$ is, locally near each of its points,
        contained in a real-analytic almost complex submanifold $S^{2m-2}\subset Z$.
    \end{enumerate}
\end{lemma}
\begin{proof}
    Let $q$ be a regular point and $v\in T_qZ$ be nonzero. In an
    analytic submanifold chart for $Z$, an appropriate straight
    parameter line gives an analytic arc $c$ in $Z$ with
    $c(0)=q$ and $c'(0)=v$. Apply Lemma~\ref{Lemma: Parametric analytic propagation} with $k=1$.
    The resulting map $\phi$ takes values in $Z$, so the curve
    $t\mapsto\phi(0,t)$ has tangent at $t=0$ in $T_qZ$. Equation~\eqref{Equation: Analytic Cauchy problem}
    identifies this tangent vector with $J_qv$. Thus $J_q(T_qZ)\subset T_qZ$;
    since $J_q^2=-\mathrm{Id}$, this inclusion is an equality.

    Every $J$-invariant real vector space has even dimension. Since
    $Z$ is a proper analytic subset of a connected open set $V$, it has empty
    interior and hence no $2m$-dimensional stratum. If it had dimension
    $2m-1$, every point of a stratum of that dimension would be regular
    by Remark~\ref{Remark: Analytic stratification}. The tangent space
    of $Z$ at such a point would then be $J$-invariant and odd-dimensional,
    a contradiction. Therefore $\dim Z\leq2m-2$.

    Fix a locally finite real-analytic stratification of $Z$ satisfying the frontier condition, and write $Z_2$ for the
    union of its strata of dimension $2m-2$. By Remark~\ref{Remark: Analytic stratification},
    $Z_2$ is relatively open in $Z$, each of its points is regular,
    and
    \begin{align*}
        Z\smallsetminus Z_2
        \text{ is relatively closed and subanalytic},\qquad
        \dim(Z\smallsetminus Z_2)\leq2m-3.
    \end{align*}
    If no stratum has dimension $2m-2$, then $Z_2=\varnothing$.

    For the second assertion, fix $q\in K$. Since
    $\dim_{\mathbb R}T_qK=2m-3$ is odd, there is a vector
    $v\in T_qK$ such that $J_qv\notin T_qK$. Indeed, otherwise
    $T_qK$ would be $J_q$-invariant and would therefore have
    even real dimension.

    Choose a real-analytic parametrisation
    $c:U\rightarrow K$ of a neighbourhood of $q$, where
    $U\subset\mathbb R^{2m-3}$ is an open neighbourhood of zero
    and $c(0)=q$. Since
    $dc_0:\mathbb R^{2m-3}\rightarrow T_qK$ is an isomorphism,
    precomposing $c$ with a suitable invertible linear map
    allows us to arrange that
    \begin{align*}
        \partial_{s_1}c(0)=v.
    \end{align*}
    The hypotheses of
    Lemma~\ref{Lemma: Parametric analytic propagation} are then
    satisfied with $k=2m-3$. It gives a real-analytic map
    $\phi$ into $Z$ whose differential has rank $2m-2$
    throughout a neighbourhood of $(0,0)$.

    By the real-analytic constant-rank theorem, after restricting
    to a smaller open neighbourhood $W$ of $(0,0)$, the map
    $\phi|_W$ is an embedding and its image
    $S:=\phi(W)$ is a real-analytic embedded submanifold of
    dimension $2m-2$. Moreover, $S\subset Z$. Since
    $\phi(s,0)=c(s)$ and $W$ contains $(s,0)$ for every $s$
    sufficiently close to zero, $S$ contains a neighbourhood
    of $q$ in $K$.

    It remains to prove that $J_x(T_xS)=T_xS$ for every $x\in S$. Put $F:=Z\smallsetminus Z_2$, so that $F$ is subanalytic
    and $\dim F\leq2m-3$.

    We first show that $S\cap Z_2$ is dense in $S$.
    Otherwise there would be a nonempty relatively open subset
    $O\subset S$ contained in $F$. Choose $x\in O$ and a
    sufficiently small ambient real-analytic coordinate
    neighbourhood $U_x$ such that $S\cap U_x\subset O$ and
    $S\cap U_x$ is an analytic coordinate slice, hence subanalytic
    in $U_x$. The inclusion
    $S\cap U_x\subset F\cap U_x$ and monotonicity of
    subanalytic dimension would then give
    \begin{align*}
        2m-2
        =\dim(S\cap U_x)
        \leq\dim(F\cap U_x)
        \leq2m-3,
    \end{align*}
    a contradiction.

    Now let $x\in S\cap Z_2$. Since $x$ is a regular point
    of $Z$, the set $Z$ agrees near $x$ with a real-analytic
    submanifold of dimension $2m-2$. The inclusion $S\subset Z$
    gives $T_xS\subset T_xZ$. Both spaces have dimension
    $2m-2$, so $T_xS=T_xZ$. Assertion~(1) therefore gives
    \begin{align*}
        J_x(T_xS)=T_xS
        \qquad\mbox{for every }x\in S\cap Z_2.
    \end{align*}

    To extend this identity to all of $S$, use the coordinates
    of $V$ to identify its tangent spaces with $\mathbb R^{2m}$,
    and let $P_x$ be the Euclidean orthogonal projection onto
    $T_xS$. Since $S$ is smooth, $P_x$ depends continuously on
    $x\in S$. The preceding invariance identity is equivalent to
    \begin{align*}
        (\mathrm{Id}-P_x)J_xP_x=0
        \qquad\mbox{for every }x\in S\cap Z_2.
    \end{align*}
    The left-hand side depends continuously on $x$, and
    $S\cap Z_2$ is dense in $S$. Hence this matrix identity
    holds for every $x\in S$. For $w\in T_xS$, we have
    $P_xw=w$, and therefore
    \begin{align*}
        (\mathrm{Id}-P_x)J_xw=0.
    \end{align*}
    Thus $J_xw\in T_xS$, proving $J_x(T_xS)\subset T_xS$.
    Since $J_x^2=-\mathrm{Id}$, this inclusion is an equality.
    Consequently $S$ is a real-analytic almost complex
    submanifold, as required.
\end{proof}

\subsection{Removability across the analytic strata}
\label{Subsection: Removability}

We first recall the version of Elkhadhra's extension theorem for
positive currents that we use below \cite{ElkhadhraPluripolar}.
Let $(V^{2m},J)$ be a smooth almost complex manifold. A \emph{$J$-holomorphic disc} is a smooth map $\gamma:D\rightarrow V$, where $D$ is an open disc
in $\mathbb C$, satisfying $d\gamma\circ i=J\circ d\gamma$.
An upper semicontinuous function
$\rho:V\rightarrow[-\infty,\infty)$ is \emph{$J$-plurisubharmonic}
if $\rho\circ\gamma$ is subharmonic on every such disc, allowing
the composition to be identically $-\infty$.

If $f$ is $C^2$ and $J$-plurisubharmonic, then
$dd_J^cf(v,Jv)\geq0$ for every tangent vector $v$. This is the
positivity of the $(1,1)$-part of $dd_J^cf$; the entire form need
not be of type $(1,1)$ when $J$ is nonintegrable. The relation
with discs follows from $\gamma^*d_J^cf=d_i^c(f\circ\gamma)$ and
the existence of local $J$-holomorphic discs with a prescribed
tangent vector.

A relatively closed set $S\subset V$ is \emph{locally regular complete
$J$-pluripolar} if near each of its points there is a
$J$-plurisubharmonic function $\rho$ whose $-\infty$ locus is
exactly $S$ in that neighbourhood and which is finite and $C^2$
off $S$. Elkhadhra proved that every smooth almost complex
submanifold of $V$ of positive codimension has this property \cite{ElkhadhraPluripolar}; the
precise statement is also recalled in
\cite[Section~3.1, Remark~1(3)]{ElkhadhraPoincareLelong}.

Let $T$ be a current of dimension $q$ with locally finite mass on
$V\smallsetminus S$, where $S$ is relatively closed in $V$.
We say that its mass is \emph{locally finite across $S$} if
\begin{align*}
    \|T\|(K\smallsetminus S)<\infty
    \qquad\mbox{for every compact }K\subset V.
\end{align*}
Under this assumption, its \emph{trivial extension} $\widetilde T$
to $V$ is the current defined by
\begin{align*}
    \widetilde T(\eta)
    :=\int_{V\smallsetminus S}
       \langle\eta_x,\vec T(x)\rangle\,d\|T\|(x),
    \qquad \eta\in C_c^\infty(V;\Lambda^qT^*V).
\end{align*}
It has locally finite mass, agrees with $T$ on
$V\smallsetminus S$, and satisfies
\begin{align*}
    \|\widetilde T\|(E)=\|T\|(E\smallsetminus S)
    \qquad\mbox{for every Borel set }E\subset V.
\end{align*}
In particular, $\|\widetilde T\|(S)=0$. If $T$ is positive,
then $\widetilde T$ is positive as well.

\begin{theorem}[Elkhadhra, removability across almost complex submanifolds]\label{Theorem: Positive current extension}
    Let $(V^{2m},J)$ be smooth, with $m\geq2$, and let
    $S\subset V$ be a relatively closed smooth almost complex
    submanifold of $V$ of positive codimension. Let $T$ be a positive
    cycle of dimension $2m-2$ and bidimension $(m-1,m-1)$ on
    $V\smallsetminus S$, whose mass is locally finite across $S$.
    Then its trivial extension $\widetilde T$ is a positive cycle
    with locally finite mass on $V$.
\end{theorem}
\begin{proof}
    Positivity and local finiteness of the mass follow from the
    definition of the trivial extension. It remains to prove
    $\partial\widetilde T=0$. This is a local assertion, and it
    already holds on $V\smallsetminus S$.

    \medskip
    \noindent\emph{Step 1: the potential and auxiliary currents.}

    Near a point of $S$, choose the $J$-plurisubharmonic potential
    $\rho$ described above: it is upper semicontinuous, its
    $-\infty$ locus is precisely $S$, and it is finite and $C^2$
    off $S$. Choose a smaller relatively compact neighbourhood
    $W$ whose closure is contained in the domain of $\rho$.
    Upper semicontinuity bounds $\rho$ from above on $\overline W$.
    Subtracting a constant, we may assume $\rho\leq-1$ on $W$.
    This preserves its $-\infty$ locus and
    $J$-plurisubharmonicity.

    Let $g$ be a smooth $J$-compatible metric and let
    $\Omega(v,w):=g(Jv,w)$. All masses and pointwise norms in
    the proof are computed using $g$. All currents below are
    restricted to $W$. Define
    \begin{align*}
        \sigma&:=\widetilde T\llcorner\Omega^{m-2},\\
        \mu&:=\widetilde T\llcorner\Omega^{m-1},\\
        R&:=\widetilde T\llcorner d(\Omega^{m-2}).
    \end{align*}
    Thus $\sigma$ is a positive current of dimension two and
    bidimension $(1,1)$, $\mu$ is a positive current of dimension
    zero, and $R$ is a current of dimension one. We identify the
    positive zero-dimensional current $\mu$ with its positive
    Radon measure. All three have locally finite mass, and
    $\mu(S)=0$.

    On $W\smallsetminus S$, the contraction formula and
    $\partial T=0$ give
    \begin{align*}
        \partial\sigma=-R.
    \end{align*}
    Since $R$ has locally finite mass, $\sigma$ is locally normal
    on $W\smallsetminus S$. This calculation uses no closedness
    assumption on $\Omega$.

    \medskip
    \noindent\emph{Step 2: a pointwise estimate for current pairings.}
    We record the linear algebra underlying the estimates.
    At $\|\widetilde T\|$-almost every point, positivity and
    bidimension $(m-1,m-1)$ allow us to choose a $g$-orthonormal
    frame $e_1,Je_1,\ldots,e_m,Je_m$ in which
    \begin{align*}
        \vec{\widetilde T}
        =\sum_{i=1}^m a_iE_{\widehat i},
        \qquad a_i\geq0.
    \end{align*}
    Here $E_{\widehat i}$ is the unit simple $(2m-2)$-vector
    obtained from
    $e_1\wedge Je_1\wedge\cdots\wedge e_m\wedge Je_m$
    by omitting $e_i\wedge Je_i$. This follows from the spectral
    theorem for the positive semidefinite Hermitian form on
    $\Lambda^{m-1,0}T_x^*W$ determined by
    $\vec{\widetilde T}$. Here we use that, in complex dimension
    $m$, every $(m-1,0)$-covector is decomposable. The resulting
    diagonalisation is used only pointwise; no regularity of the
    chosen frame is needed.

    The form $\Omega^{m-1}/(m-1)!$ has comass one and takes
    the value one on every $E_{\widehat i}$. The triangle
    inequality and this evaluation show that the mass norm of
    $\sum_i a_iE_{\widehat i}$ is exactly $\sum_i a_i$.
    Since $\vec{\widetilde T}$ has mass norm one, $\sum_i a_i=1$.
    It follows that
    \begin{align*}
        \mu=(m-1)!\,\|\widetilde T\|.
    \end{align*}
    We may therefore write
    \begin{align*}
        \widetilde T(\eta)
        &=\int_W\langle\eta,\tau\rangle\,d\mu,\\
        \tau&=\sum_{i=1}^m\lambda_iE_{\widehat i},
        \qquad \lambda_i\geq0,
        \qquad \sum_{i=1}^m\lambda_i=\frac1{(m-1)!}.
    \end{align*}
    All these pointwise identities hold $\mu$-almost everywhere.
    For $\widetilde T=T_{u^*\omega}$, take $g$ to be the metric
    used to define $\mu_u$. The energy identity then gives $\mu=(m-1)!\mu_u$.

    For a continuous real covector field $\xi$, put
    $\xi^c:=-\xi\circ J$. The zero-dimensional current
    $\sigma\llcorner(\xi\wedge\xi^c)$ is a positive measure
    absolutely continuous with respect to $\mu$. Write it as
    $q_\xi\mu$. If
    $|\xi_j|^2:=\xi(e_j)^2+\xi(Je_j)^2$, expansion in the
    preceding frame gives
    \begin{align*}
        q_\xi
        &=(m-2)!\sum_{i=1}^m\lambda_i
                        \sum_{j\ne i}|\xi_j|^2,\\
        q_{\xi^c}&=q_\xi.
    \end{align*}
    Similarly, for a continuous real $(2m-3)$-form $\Theta$,
    the zero-dimensional current
    $\widetilde T\llcorner(\Theta\wedge\xi)$ is a signed
    measure with density $\langle\Theta\wedge\xi,\tau\rangle$
    relative to $\mu$. Since $E_{\widehat i}$ contains only
    the complex lines with $j\ne i$, the Cauchy--Schwarz
    inequality gives
    \begin{align*}
        |\langle\Theta\wedge\xi,\tau\rangle|
        &\leq C_m|\Theta|_g
             \sum_i\lambda_i
             \left(\sum_{j\ne i}|\xi_j|^2\right)^{1/2}\\
        &\leq C_m|\Theta|_g
             \left(\sum_i\lambda_i\right)^{1/2}
             \left(\sum_i\lambda_i
                         \sum_{j\ne i}|\xi_j|^2\right)^{1/2}\\
        &\leq C_m|\Theta|_g\,q_\xi^{1/2}.
    \end{align*}
    The same estimate holds with $\xi\wedge\Theta$ in place
    of $\Theta\wedge\xi$, since these differ only by a sign.
    The formula for $q_\xi$ also shows that the bilinear pairing
    \begin{align*}
        (\xi,\zeta)\mapsto
        \langle\Omega^{m-2}\wedge\xi\wedge\zeta^c,\tau\rangle
    \end{align*}
    is symmetric, since it equals
    $(m-2)!\sum_i\lambda_i\sum_{j\ne i}
       [\xi(e_j)\zeta(e_j)+\xi(Je_j)\zeta(Je_j)]$.

    \medskip
    \noindent\emph{Step 3: cutoffs near the submanifold.}
    Choose a smooth convex nondecreasing function $\chi$ on
    $(-\infty,0]$ such that $\chi=0$ on $(-\infty,-1]$,
    $\chi(0)=1$, and $(\chi')^2\leq C\chi''$.
    For example, set $\chi(t)=\exp(4-4/(t+1))$ for
    $-1<t\leq0$ and extend it by zero to the left. The
    extension is smooth at $-1$. With $z=t+1\in(0,1]$,
    \begin{align*}
        \chi'=\frac{4\chi}{z^2},\qquad
        \chi''=\frac{8\chi(2-z)}{z^4},\qquad
        \frac{(\chi')^2}{\chi''}
        =\frac{2\chi}{2-z}\leq2.
    \end{align*}
    For each integer $\ell\geq1$, define
    \begin{align*}
        f_\ell(x):=
        \begin{cases}
            \chi(\rho(x)/\ell),&x\in W\smallsetminus S,\\
            0,&x\in S.
        \end{cases}
    \end{align*}
    The open set $\{\rho<-\ell\}$ contains $S$, and $f_\ell$
    vanishes there. Thus $f_\ell\in C^2(W)$,
    $0\leq f_\ell\leq1$, and $f_\ell\to1$ pointwise off $S$.
    If a product involving $f_\ell$ or its derivatives is
    localised to a compact subset of $W$, its support is
    contained in the intersection of that compact set with
    $\{\rho\geq-\ell\}$. This intersection is compact and
    disjoint from $S$.

    Fix $\eta\in C_c^\infty(W;\Lambda^{2m-3}T^*W)$.
    Choose $\psi\in C_c^\infty(W)$ with $0\leq\psi\leq1$
    and $\psi=1$ on a neighbourhood of $\spt\eta$.
    Let $K\subset W$ be a compact set containing $\spt\psi$.
    Put
    \begin{align*}
        E_\ell
        &:=\sigma(\psi\,df_\ell\wedge d_J^cf_\ell)
          =\int_W\psi q_{df_\ell}\,d\mu,\\
        I_\ell&:=\sigma(\psi\,dd_J^cf_\ell).
    \end{align*}
    These quantities are finite, because $f_\ell$ is $C^2$
    and the currents involved have locally finite mass.
    On $W\smallsetminus S$, the chain rule gives
    \begin{align*}
        dd_J^cf_\ell
        =\frac{\chi'(\rho/\ell)}\ell\,dd_J^c\rho
         +\frac{\chi''(\rho/\ell)}{\ell^2}
          \,d\rho\wedge d_J^c\rho.
    \end{align*}
    Since $\sigma$ has bidimension $(1,1)$, only the
    $(1,1)$-part of $dd_J^c\rho$ contributes to its pairing.
    This part is positive by $J$-plurisubharmonicity. Positivity
    of $\sigma$, $\chi'\geq0$, and
    $(\chi')^2\leq C\chi''$ therefore imply
    \begin{align*}
        0\leq E_\ell\leq C I_\ell.
    \end{align*}
    There is no contribution near $S$, where $f_\ell$ vanishes.

    \medskip
    \noindent\emph{Step 4: a uniform energy bound.}
    The 1-forms $\psi d_J^cf_\ell$ and $f_\ell d_J^c\psi$
    have compact support in $W\smallsetminus S$. Using
    $\partial\sigma=-R$ there gives
    \begin{align*}
        -R(\psi d_J^cf_\ell)
        &=\sigma\bigl(d(\psi d_J^cf_\ell)\bigr)\\
        &=\sigma(d\psi\wedge d_J^cf_\ell)+I_\ell.
    \end{align*}
    The symmetry proved in Step 2 yields
    \begin{align*}
        \sigma(d\psi\wedge d_J^cf_\ell)
        =\sigma(df_\ell\wedge d_J^c\psi).
    \end{align*}
    Applying the boundary identity to $f_\ell d_J^c\psi$ gives
    \begin{align*}
        -R(f_\ell d_J^c\psi)
        &=\sigma\bigl(d(f_\ell d_J^c\psi)\bigr)\\
        &=\sigma(df_\ell\wedge d_J^c\psi)
          +\sigma(f_\ell dd_J^c\psi).
    \end{align*}
    Combining these equalities, we obtain
    \begin{align*}
        I_\ell
        =\sigma(f_\ell dd_J^c\psi)
         +R(f_\ell d_J^c\psi)
         -R(\psi d_J^cf_\ell).
    \end{align*}
    These boundary identities remain valid for the $C^1$
    1-forms used here: approximate them in $C^1$ by smooth
    1-forms supported in a fixed compact subset of
    $W\smallsetminus S$. Their exterior derivatives converge
    uniformly, and both $\sigma$ and $R$ have locally finite
    mass, so the pairings converge.

    Since $0\leq f_\ell\leq1$ and $\Omega$, $d(\Omega^{m-2})$,
    and the derivatives of $\psi$ are bounded on $K$, the integral
    representation in Step 2 bounds each of the first two terms
    in absolute value by $C\mu(K)$.
    For the last term, apply the estimate of Step 2 with
    $\Theta=d(\Omega^{m-2})$ and $\xi=d_J^cf_\ell$.
    Since $q_{d_J^cf_\ell}=q_{df_\ell}$, we get
    \begin{align*}
        |R(\psi d_J^cf_\ell)|
        &\leq C\int_W\psi q_{df_\ell}^{1/2}\,d\mu\\
        &\leq C\left(\int_W\psi\,d\mu\right)^{1/2}
                 \left(\int_W\psi q_{df_\ell}\,d\mu\right)^{1/2}\\
        &\leq C\mu(K)^{1/2}E_\ell^{1/2}.
    \end{align*}
    Together with $E_\ell\leq C I_\ell$, this implies
    \begin{align*}
        E_\ell\leq C\mu(K)
                      +C\mu(K)^{1/2}E_\ell^{1/2}.
    \end{align*}
    Young's inequality absorbs half of $E_\ell$ into the
    left-hand side and gives $E_\ell\leq C'\mu(K)$ for all
    $\ell$. The constants depend on the fixed compact set and
    cutoff, but not on $\ell$.

    \medskip
    \noindent\emph{Step 5: the boundary of the extension.}
    We claim that
    \begin{align*}
        \widetilde T\llcorner df_\ell\rightarrow0
        \qquad\mbox{as currents on }W.
    \end{align*}
    Evaluating the contraction on the fixed test form $\eta$ gives
    \begin{align*}
        (\widetilde T\llcorner df_\ell)(\eta)
        =\int_W\langle df_\ell\wedge\eta,\tau\rangle\,d\mu.
    \end{align*}
    For $a>0$, split this integral over $\{\rho<-a\}$ and
    $\{\rho\geq-a\}$. On the first set, the estimate of
    Step 2, $\psi=1$ near $\spt\eta$, and Cauchy--Schwarz give
    \begin{align*}
        \left|\int_{\{\rho<-a\}}
            \langle df_\ell\wedge\eta,\tau\rangle\,d\mu\right|
        \leq C_\eta E_\ell^{1/2}
                    \mu(K\cap\{\rho<-a\})^{1/2}.
    \end{align*}
    The sets $K\cap\{\rho<-a\}$ decrease to $K\cap S$ as
    $a\rightarrow\infty$. Since $\mu(K)<\infty$ and $\mu(S)=0$,
    their measures tend to zero. The uniform bound for $E_\ell$
    therefore shows that this first integral tends to zero as
    $a\rightarrow\infty$, uniformly in $\ell$.

    For fixed $a$, the set $K\cap\{\rho\geq-a\}$ is compact
    and disjoint from $S$. The covector field $d\rho$ is bounded
    there, and boundedness of $\chi'$ gives
    $|df_\ell|_g\leq C_a/\ell$. The absolute value of the
    integral over this second set is consequently at most
    $C_{\eta,a}\mu(K)/\ell$. Taking first
    $\limsup_{\ell\rightarrow\infty}$ and then $a\rightarrow\infty$ proves
    the claim.

    Dominated convergence in the integral representation of
    $\widetilde T$ also gives
    \begin{align*}
        f_\ell\widetilde T\rightarrow\widetilde T
        \qquad\mbox{as currents on }W,
    \end{align*}
    since $f_\ell\to1$ off $S$ and
    $\|\widetilde T\|(S)=0$. Finally, for every
    $\eta\in C_c^\infty(W;\Lambda^{2m-3}T^*W)$, the form
    $f_\ell\eta$ has compact support in $W\smallsetminus S$.
    Since $T$ is a cycle there,
    \begin{align*}
        0
        &=T\bigl(d(f_\ell\eta)\bigr)\\
        &=\widetilde T(df_\ell\wedge\eta)
          +\widetilde T(f_\ell d\eta).
    \end{align*}
    As above, this identity for $C^1$ test forms follows by
    approximation and the local finiteness of the mass. Hence
    \begin{align*}
        \partial(f_\ell\widetilde T)
        =-\widetilde T\llcorner df_\ell.
    \end{align*}
    The right-hand side tends to zero by the preceding estimate.
    Since the boundary operator is continuous under convergence of
    currents, the left-hand side tends to $\partial\widetilde T$.
    Thus $\partial\widetilde T=0$ on $W$, and the conclusion follows
    by locality.
\end{proof}

For a related application, see
\cite[Section~3.2, proof of Lemma~1]{ElkhadhraPoincareLelong}.
In our setting, let $S\subset V$ be as in
Theorem~\ref{Theorem: Positive current extension}, and let $\alpha$
be a positive $(1,1)$-form with $L^1_{\mathrm{loc}}$ coefficients
on $V$. Its restriction to $V\smallsetminus S$ has locally finite mass across $S$. Since $S$ has positive codimension, it has zero
Lebesgue measure, so the trivial extension of this restriction
is exactly the current defined by $\alpha$. Hence
\begin{align}\label{Equation: Positive pullback extension}
    d\alpha=0\text{ on }V\smallsetminus S
    \quad\Rightarrow\quad d\alpha=0\text{ on }V.
\end{align}
Only the mass of $\alpha$ is used; no mass assumption on
$d\alpha$ is needed.

For strata of higher codimension, local integrability gives the
following support estimate for the exterior derivative.

\begin{lemma}[Removability in higher codimension]\label{Lemma: High codimension support}
    Let $\alpha$ be a locally integrable differential $k$-form on a smooth $d$-manifold. If $D=d\alpha$ is supported in a relatively closed embedded smooth submanifold $S$ of codimension $c>k+1$, then $D=0$. If the ambient manifold is real analytic, the same conclusion holds when $S$ is relatively closed and subanalytic with $\dim S\leq d-k-2$.
\end{lemma}
\begin{proof}
    The smooth-submanifold assertion is local. Choose coordinates
    $(x,y)\in\mathbb R^{d-c}\times\mathbb R^c$ in which
    $S=\{y=0\}$. Let $\vartheta\in C_c^\infty(\mathbb R^c)$
    equal one near zero, with support in $\{|y|\leq C\}$, and put
    $\vartheta_r(y)=\vartheta(y/r)$ for $r>0$. Fix a smooth
    $(d-k-1)$-form $\eta$ compactly supported in this coordinate
    neighbourhood and a compact $K$ containing its support.
    For each normal coordinate $y^a$, the support assumption gives
    \begin{align*}
        \langle y^aD,\eta\rangle
        =\langle D,y^a\eta\rangle
        =\langle D,\vartheta_r y^a\eta\rangle.
    \end{align*}
    Indeed, $(1-\vartheta_r)y^a\eta$ vanishes on an open
    neighbourhood of $S$. Expanding the exterior derivative gives
    \begin{align*}
        d(\vartheta_r y^a\eta)
        =y^a d\vartheta_r\wedge\eta
        +\vartheta_rdy^a\wedge\eta
        +\vartheta_ry^a d\eta.
    \end{align*}
    All coefficients on the right are uniformly bounded for small
    $r$. The only factor of order $r^{-1}$ occurs in
    $d\vartheta_r$, where it is multiplied by $y^a=O(r)$ on its
    support. All three terms are supported in
    $K\cap\{|y|\leq Cr\}$. By the distributional convention
    for $D=d\alpha$ and local equivalence of coordinate norms,
    \begin{align*}
        |\langle D,\vartheta_r y^a\eta\rangle|
        &=\left|(-1)^{k+1}
        \int\alpha\wedge d(\vartheta_r y^a\eta)\right|\\
        &\leq C_\eta\int_{K\cap\{|y|\leq Cr\}}|\alpha|
        \rightarrow0.
    \end{align*}
    The limit follows from absolute continuity of the integral:
    $\alpha$ is locally integrable and $S$ has zero Lebesgue
    measure. Thus $y^aD=0$ for every $a$.

    Since $dD=d^2\alpha=0$, the distributional product rule now
    implies
    \begin{align*}
        dy^a\wedge D=d(y^aD)-y^a dD=0
        \qquad(1\leq a\leq c).
    \end{align*}
    Expand $D$ in the coordinate basis of $(k+1)$-forms, with
    distributional coefficients. For a fixed $a$, wedging with
    $dy^a$ sends basis monomials not containing $dy^a$ to distinct
    basis monomials, up to sign. Their coefficients must therefore
    vanish separately. Repeating for every $a$ shows that any
    surviving term would contain all $c$ distinct normal factors.
    Its degree is only $k+1<c$, so there are no surviving terms.
    This proves the first assertion.

    For the second assertion, choose a locally finite real-analytic
    stratification of the closed containing set $S$. Its maximal
    stratum dimension is at most $d-k-2$. At a point of any
    top-dimensional stratum, $S$ is locally precisely that embedded
    submanifold by Remark~\ref{Remark: Analytic stratification}. Its codimension
    is at least $k+2$, so the first assertion proves that $D$
    vanishes on a neighbourhood of the point. After removing all
    strata of this dimension, the support of $D$ is contained in
    the closed lower-dimensional skeleton. Apply the same argument
    to that skeleton. Since the dimension drops at each step,
    after finitely many steps $\spt D$ is empty. This argument uses
    subanalyticity only for the skeletons containing $\spt D$.
\end{proof}

\subsection{Closedness for sphere-valued maps}
\label{Subsection: Sphere in every dimension}

\begin{theorem}\label{Theorem: Sphere valued maps}
    For every $m\geq1$, every weakly pseudoholomorphic map
    $u\in W^{1,2}_{\mathrm{loc}}(M^{2m},\cp^1)$ on a
    real-analytic almost complex manifold without boundary satisfies
    $d(u^*\omega)=0$ in the sense of distributions on $M$.
\end{theorem}
\begin{proof}
    We proceed by induction on $m$. If $m=1$, then $u^*\omega$
    has top degree and is therefore closed.
    Assume the assertion in real dimension $2m-2$, and let the
    domain have dimension $2m$, with $m\geq2$. Fix $p\in M$. If $J$ is
    integrable near $p$, Proposition~\ref{Proposition: Integrable sphere valued case}
    gives closedness there. Otherwise, the second alternative in
    Lemma~\ref{Lemma: Finite analytic annihilators} gives a neighbourhood
    $V$ of $p$ and fields $X_1,\ldots,X_r$ with common zero set $Z$.
    Set $\alpha:=(u^*\omega)|_V$ and $D:=d\alpha$ on $V$.

    \medskip
    \noindent\emph{Step 1: reduction away from the common zero set.}
    At every $q\in V\smallsetminus Z$ some $X_i(q)$ is nonzero.
    Lemma~\ref{Lemma: Analytic dimension reduction} then factors $u$ on a neighbourhood
    of $q$ through a weakly pseudoholomorphic
    $W^{1,2}_{\mathrm{loc}}$ map on a real-analytic almost complex
    domain of dimension $2m-2$. By the induction hypothesis, the pullback of $\omega$ by the reduced map is closed.
    The last assertion of the reduction lemma therefore gives $D=0$
    near $q$. Since $q\in V\smallsetminus Z$ was arbitrary,
    \begin{align}\label{Equation: Support in the analytic zero set}
        \spt D\subset Z.
    \end{align}
    \medskip
    \noindent\emph{Step 2: removal of the codimension-two strata.}
    By Lemma~\ref{Lemma: Geometry of the common zero set}, $\dim Z\leq2m-2$.
    Fix a real-analytic stratification of $Z$ satisfying
    Remark~\ref{Remark: Analytic stratification}, and let $Z_2$ be the union
    of its strata of dimension $2m-2$, as in the proof of that
    lemma. Each point of $Z_2$ is regular, so it has a neighbourhood $W$
    in which $Z\cap W$ is a relatively closed embedded analytic submanifold
    of dimension $2m-2$. Its tangent spaces are $J$-invariant by
    Lemma~\ref{Lemma: Geometry of the common zero set}(1).

    By \eqref{Equation: Positivity of the pullback}, $\alpha$ is a positive
    $(1,1)$-form with locally integrable coefficients. Since
    \eqref{Equation: Support in the analytic zero set} gives closedness
    on $W\smallsetminus Z$, \eqref{Equation: Positive pullback extension}
    gives closedness throughout $W$. Thus $D$ vanishes near every point
    of $Z_2$, and
    \begin{align*}
        \spt D\subset F:=Z\smallsetminus Z_2,
        \qquad \dim F\leq2m-3.
    \end{align*}
    The set $F$ is the closed subanalytic $(2m-3)$-skeleton of
    the fixed stratification and inherits its stratification.

    \medskip
    \noindent\emph{Step 3: removal of the codimension-three strata.}
    Let $q$ lie on a stratum of $F$ of dimension $2m-3$.
    By local finiteness and the frontier condition, $F$ agrees near
    $q$ with a relatively closed piece $K$ of that stratum. Lemma~\ref{Lemma: Geometry of the common zero set}(2) gives, after shrinking,
    a real-analytic almost complex submanifold
    $S^{2m-2}\subset Z$ containing $K$ near $q$. Choose a still
    smaller neighbourhood $W$ of $q$ in which $S$ is relatively
    closed and $F\cap W\subset S$. Since $\spt D\subset F$, the form $\alpha$ is closed on
    $W\smallsetminus S$. Moreover, $\alpha$ is locally integrable and
    $S$ has zero Riemannian volume, so the trivial extension of
    $T_\alpha|_{W\smallsetminus S}$ is $T_\alpha|_W$.
    Thus \eqref{Equation: Positive pullback extension} gives $D=0$ on $W$.

    \medskip
    \noindent\emph{Step 4: removal of the remaining support.}
    The support of $D$ is now contained in the closed subanalytic
    $(2m-4)$-skeleton of the same stratification. This set has
    codimension at least four. Lemma~\ref{Lemma: High codimension support}
    with $d=2m$ and $k=2$ gives $D=0$ on $V$.
    Since $p$ was arbitrary, the induction is complete.
\end{proof}

\begin{remark}
    The codimension-two strata are almost complex submanifolds.
    For a codimension-three stratum, parametric analytic propagation
    instead gives a codimension-two almost complex submanifold
    containing it locally. The elementary $L^1$ support lemma removes the
    remaining strata, which have codimension at least four.
    The argument therefore does not require the lower strata
    of $Z$ themselves to be almost complex.
\end{remark}

\subsection{Projective spaces and projective algebraic targets}
\label{Subsection: Projective targets}

To distinguish the target dimensions in the induction, we write $j_k,\omega_k,h_k$ for the standard complex
structure, normalised Fubini--Study form, and metric on $\cp^k$,
respectively. For the original target, $j_n=j$, $\omega_n=\omega$,
and $h_n=h$. We use the singular projections of
\cite[Appendix~A]{CaniatoRiviere}, whose averaging formula
\cite[Lemma~A.2]{CaniatoRiviere} reduces the target dimension one step at a time.
To apply it here, we first justify
composition with these projections for arbitrary $W^{1,2}_{\mathrm{loc}}$ maps.

For $k\geq2$ and $p\in\cp^k$, let $L_p\subset\mathbb C^{k+1}$ be
the complex line represented by $p$, and let $Q_p$ be the
orthogonal projection onto $L_p^\perp$. As in
\cite[Appendix~A]{CaniatoRiviere}, define
\begin{align*}
    F_p:\cp^k\smallsetminus\{p\}\rightarrow\mathbb P(L_p^\perp),
    \qquad F_p([z])=[Q_pz].
\end{align*}
A unitary isomorphism $L_p^\perp\simeq\mathbb C^k$, chosen
measurably in $p$, identifies $\mathbb P(L_p^\perp)$ with $\cp^{k-1}$. Such choices can
be made using smooth local orthonormal frames on a finite
cover and a Borel partition subordinate to that cover. The
pulled-back Fubini--Study form and metric are independent of
this choice. We equip the parameter space with the probability
measure
\begin{align*}
    d\lambda_k(p):=
    \frac{d\vol_{h_k}(p)}{\vol_{h_k}(\cp^k)},
    \qquad h_k(\cdot,\cdot):=\omega_k(\cdot,j_k\cdot).
\end{align*}
All differential norms in the following lemma are Hilbert--Schmidt
norms, computed with the domain metric and these normalised
Fubini--Study metrics on the targets. Parameter integrals of forms
mean integrals of their coefficients in local coordinates. We will
prove their absolute integrability before interchanging them with
domain integration.

\begin{lemma}
    \label{Lemma: Sobolev projective projections}
    Let $u\in W^{1,2}_{\mathrm{loc}}(M,\cp^k)$, where $k\geq2$ and
    $M$ is any smooth manifold equipped with a Riemannian metric $g$. For $\lambda_k$-almost
    every $p$, the composition $v_p:=F_p\circ u$ is defined almost
    everywhere, belongs to $W^{1,2}_{\mathrm{loc}}(M,\cp^{k-1})$, and
    satisfies
    \begin{align*}
        dv_p=dF_p(u)\circ du\qquad\text{almost everywhere}.
    \end{align*}
    Moreover,
    \begin{align}\label{Equation: Averaged Sobolev pullbacks}
        u^*\omega_k=\int_{\cp^k}v_p^*\omega_{k-1}\,d\lambda_k(p)
        \qquad\text{in }L^1_{\mathrm{loc}},
    \end{align}
    and, for every relatively compact coordinate neighbourhood
    $U\subset M$,
    \begin{align}\label{Equation: Averaged projection energy}
        \int_{\cp^k}\int_U|dv_p|^2\,d\vol_g\,d\lambda_k(p)
        =\int_U|du|^2\,d\vol_g.
    \end{align}
    If $M$ is also equipped with an almost complex structure $J$ and
    $u$ is weakly pseudoholomorphic with respect to $J$, then so is
    $v_p$ for almost every $p$.
\end{lemma}
\begin{proof}
    \noindent\emph{Step 1: Sobolev composition.}
    Identify $\cp^k$ with the set of rank-one orthogonal projectors
    $P=zz^*$, with $|z|=1$, in the real vector space
    $\operatorname{Herm}(k+1)$ of Hermitian matrices, equipped with
    the Hilbert--Schmidt norm. Here the superscript $*$ denotes
    conjugate transpose.
    Set
    \begin{align*}
        \rho_p(P):=\tr(Q_pPQ_p)=|Q_pz|^2.
    \end{align*}
    For fixed $P$, unitary invariance of $\lambda_k$ gives
    \begin{align}\label{Equation: Integrable projection majorant}
        \int_{\cp^k}\frac{d\lambda_k(p)}{\rho_p(P)}
        =\int_0^1\frac{1}{t}\,kt^{k-1}\,dt
        =\frac{k}{k-1}.
    \end{align}
    Indeed, if $p=[a]$ with $|a|=1$, then
    $\rho_p(P)=1-|\langle a,z\rangle|^2$ has density
    $kt^{k-1}$ on $(0,1)$. To see this, normalise a vector of
    $k+1$ independent complex Gaussian variables, each with density
    $\pi^{-1}e^{-|z|^2}$ on $\mathbb C$, to obtain
    a uniform unit vector $a$, and use a unitary transformation to
    take $z=e_0$. If $X$ is the squared modulus of the first Gaussian
    variable and $Y$ the sum of the remaining squared moduli, then $X$ and $Y$ are independent, $X$ has density
    $e^{-x}$, and $Y$ is a sum of $k$ such variables. Hence,
    for $0<t<1$,
    \begin{align*}
        \Pr\left(\frac{Y}{X+Y}\leq t\right)
        =\mathbb E\left[e^{-((1-t)/t)Y}\right]=t^k,
    \end{align*}
    which gives the asserted density. We now use
    \eqref{Equation: Integrable projection majorant} to control
    the squared differential of the projection, as in the estimate
    of \cite[Lemma~A.1]{CaniatoRiviere}.

    For $\varepsilon>0$, the matrix-valued map
    \begin{align*}
        G_{p,\varepsilon}(P):=\frac{Q_pPQ_p}{\rho_p(P)+\varepsilon}
    \end{align*}
    is smooth on all of $\cp^k$ and has norm at most one. Where
    $\rho_p(P)>0$, it converges smoothly to
    \begin{align*}
        G_p(P):=\frac{Q_pPQ_p}{\rho_p(P)},
    \end{align*}
    the orthogonal projector representing $F_p(P)$. If $v\perp z$
    and $H=vz^*+zv^*$, write $x=Q_pz$ and $y=Q_pv$.
    Differentiating numerator and denominator gives
    \begin{align*}
        |dG_{p,\varepsilon}(P)[H]|
        &\leq\frac{2|x||y|}{|x|^2+\varepsilon}
        +\frac{2|x|^3|y|}{(|x|^2+\varepsilon)^2}\\
        &\leq\frac{4|y|}{|x|}
        \leq\frac{4|v|}{\sqrt{\rho_p(P)}}.
    \end{align*}
    Consequently, with respect to the fixed smooth metrics,
    \begin{align*}
        |d(G_{p,\varepsilon}\circ u)|^2
        \leq\frac{C_k}{\rho_p(u)}|du|^2
    \end{align*}
    almost everywhere on $\{\rho_p(u)>0\}$.

    Fix a relatively compact coordinate neighbourhood $U$ and choose a
    Borel representative of $u$. For each $x\in U$, the set of centres $p$ for which
    $\rho_p(u(x))=0$ is the singleton $\{u(x)\}$, which has
    zero $\lambda_k$-measure. Fubini's theorem
    therefore gives $\rho_p(u)>0$ almost everywhere on $U$ for
    almost every $p$. By \eqref{Equation: Integrable projection majorant} and
    Tonelli's theorem,
    \begin{align*}
        \int_{\cp^k}\int_U\frac{|du|^2}{\rho_p(u)}
        \,d\vol_g\,d\lambda_k(p)
        =\frac{k}{k-1}\int_U|du|^2\,d\vol_g<\infty.
    \end{align*}
    For almost every $p$, the preceding differential estimate thus
    has an integrable majorant independent of $\varepsilon$. The smooth
    Sobolev chain rule and dominated convergence show that
    $G_{p,\varepsilon}\circ u\rightarrow G_p\circ u$ strongly in $W^{1,2}(U)$,
    with differential $dG_p(u)\circ du$. A countable exhaustion
    gives the asserted local Sobolev regularity and chain rule.
    The pseudoholomorphicity statement follows from the established
    chain rule and the complex linearity of $dF_p$ away from $p$.
    The family of values and differentials is measurable in $(p,x)$:
    in the matrix representation the formulas are smooth off
    $\{(p,P):\rho_p(P)=0\}$, and the chosen target identifications
    are Borel. Choosing a Borel representative of $du$ and extending the pullback and energy integrands by zero on the exceptional set gives jointly Borel representatives.

    \medskip
    \noindent\emph{Step 2: averaging.}
    With the probability measure $\lambda_k$ and our normalisation
    $\int_{\cp^1}\omega_1=1$, the identity of
    \cite[Lemma~A.2]{CaniatoRiviere} reads
    \begin{align}\label{Equation: normalised projection average}
        \omega_k=\int_{\cp^k}F_p^*\omega_{k-1}\,d\lambda_k(p).
    \end{align}
    To verify the normalisation, first observe that the averaged
    2-form is finite by the differential bound in Step~1 and
    \eqref{Equation: Integrable projection majorant}. For each
    $A\in\U(k+1)$, the maps $F_p\circ\mathbb P(A)$ and $F_{A^{-1}p}$
    agree up to a unitary identification of their targets. Since
    $\lambda_k$ is unitary invariant, so is the averaged 2-form. At a fixed point,
    the unitary transformations fixing that point act as the full
    unitary group on its complex tangent space. A Hermitian matrix
    invariant under all these transformations is a scalar multiple of
    the identity. Since the unitary group acts transitively on $\cp^k$,
    the scalar is independent of the point. Thus the average is
    $c\omega_k$ for a constant $c\geq0$. On a projective line
    $L\subset\cp^k$, the restriction $F_p|_L$ is a projective
    isomorphism onto a projective line whenever $p\notin L$.
    The exceptional centres have zero $\lambda_k$-measure. Positivity and Tonelli's theorem allow integration over $L$ first:
    \begin{align*}
        c=c\int_L\omega_k
        =\int_{\cp^k}\left(\int_LF_p^*\omega_{k-1}\right)
        \,d\lambda_k(p)=1.
    \end{align*}
    Indeed, for $p\notin L$ the restriction has degree one, and
    the Fubini--Study form has integral one on its image.

    Since each $F_p$ is holomorphic, the corresponding metric
    identity also holds:
    \begin{align*}
        h_k=\int_{\cp^k}F_p^*h_{k-1}\,d\lambda_k(p).
    \end{align*}
    Apply these pointwise identities to $du$ and use the chain
    rule already proved. Tonelli's theorem gives
    \eqref{Equation: Averaged projection energy}; the estimate
    $|v_p^*\omega_{k-1}|\leq C_k|dv_p|^2$ then gives absolute
    integrability and proves \eqref{Equation: Averaged Sobolev pullbacks} in
    $L^1_{\mathrm{loc}}$.
\end{proof}

\begin{remark}
    In \cite[Lemma~5.3]{CaniatoRiviere} singular projective
    projections are applied to locally strongly approximable maps.
    Lemma~\ref{Lemma: Sobolev projective projections} gives the chain
    rule for arbitrary $W^{1,2}_{\mathrm{loc}}$ maps by regularising
    the projection in the ambient space of Hermitian matrices.
\end{remark}

\begin{proof}[Proof of Theorem~\ref{Theorem: Automatic closedness}]
    The case $m=1$ is automatic by degree, so assume $m\geq2$.
    We argue by induction on the target dimension $n$. The case
    $n=1$ is Theorem~\ref{Theorem: Sphere valued maps}. Suppose that the conclusion
    holds for maps into $\cp^{n-1}$ and let
    $u\in W^{1,2}_{\mathrm{loc}}(M,\cp^n)$ be weakly
    pseudoholomorphic. By Lemma~\ref{Lemma: Sobolev projective projections}, for almost
    every $p\in\cp^n$ the map $v_p=F_p\circ u$ is a weakly
    pseudoholomorphic map in
    $W^{1,2}_{\mathrm{loc}}(M,\cp^{n-1})$. The induction hypothesis
    gives $d(v_p^*\omega_{n-1})=0$.

    Let $\eta$ be a test form as in \eqref{Equation: Distributional closedness}.
    The $L^1$ identity \eqref{Equation: Averaged Sobolev pullbacks}, the energy bound
    \eqref{Equation: Averaged projection energy}, and Fubini's theorem give
    \begin{align*}
        \int_M u^*\omega_n\wedge d\eta
        =\int_{\cp^n}\left(\int_M v_p^*\omega_{n-1}\wedge d\eta\right)
        \,d\lambda_n(p)=0.
    \end{align*}
    This completes the induction and the proof.
\end{proof}

\begin{corollary}[Possibly singular projective targets]
    \label{Corollary: Projective algebraic targets}
    Let $(M^{2m},J)$ be a real-analytic almost complex manifold
    without boundary, with $m\geq 1$, and let $X\subset\cp^n$ be a
    complex projective algebraic variety. Every weakly pseudoholomorphic
    map $u\in W^{1,2}_{\mathrm{loc}}(M,X)$ satisfies $d(u^*\omega)=0$
    in the sense of distributions on $M$.
\end{corollary}
\begin{proof}
    By the ambient convention of Section~\ref{Section: Introduction},
    $u$ is a weakly pseudoholomorphic map into $\cp^n$ taking values
    in $X$ almost everywhere. Theorem~\ref{Theorem: Automatic closedness}
    applies directly.
\end{proof}

\begingroup
\linespread{1.07}\selectfont
\bibliographystyle{amsalpha}
\bibliography{main}
\endgroup
\end{document}